\documentclass[12pt,dvipsnames]{article}

\usepackage[english]{babel}
\usepackage[a4paper]{geometry}

\usepackage{mathtools}
\usepackage{amssymb}
\usepackage{amsthm}
\usepackage{bm}
\usepackage{bbm}
\usepackage{dsfont}
\usepackage{mathrsfs}
\usepackage{tensor}

\usepackage{quiver}
\usepackage[export]{adjustbox}
\usepackage[svgnames]{xcolor}

\usepackage{enumitem}

\usepackage[linktoc=page]{hyperref}
\hypersetup{
    pdfauthor={T. Franke, I. Runkel, and V. Schomerus},
    colorlinks=true,  
    breaklinks=true,
    pdfstartview={FitH}, 
    pdftitle={Braided module categories from small braidings}, 
    linkcolor=blue,
    citecolor=DarkGreen,
    urlcolor=magenta,
    filecolor=magenta
}

\numberwithin{equation}{section}
\tikzset{every picture/.style={line width=0.5mm}}

\def\centerarc[#1](#2)(#3:#4:#5);%
{
	\draw[#1]([shift=(#3:#5)]#2) arc (#3:#4:#5);
}

\theoremstyle{definition}
\newtheorem{definition}{Definition}[section]
\newtheorem*{definition*}{Definition}
\newtheorem{remark}[definition]{Remark}
\newtheorem{theorem}[definition]{Theorem}
\newtheorem*{theorem*}{Theorem}
\newtheorem{corollary}[definition]{Corollary}
\newtheorem{lemma}[definition]{Lemma}
\newtheorem*{lemma*}{Lemma}

\newtheorem{proposition}[definition]{Proposition}
\newtheorem*{proposition*}{Proposition}
\newtheorem{example}[definition]{Example}
\newtheorem*{example*}{Example}

\newcommand{\cat}[1]{\mathcal{#1}}

\newcommand{\C}{\ensuremath{\mathbb{C}}}
\newcommand{\Z}{\ensuremath{\mathbb{Z}}}
\newcommand{\1}{\ensuremath{\mathbbm{1}}}

\DeclareMathOperator{\Id}{id}
\DeclareMathOperator{\Hom}{Hom}

\DeclareMathOperator{\ev}{ev}
\DeclareMathOperator{\coev}{coev}
\newcommand{\lc}{\ensuremath{\overset{\rightarrow}{\sigma}}}
\newcommand{\rc}{\ensuremath{\overset{\leftarrow}{\sigma}}}
\newcommand{\blank}{{-}}
\newcommand{\bblank}{{=}}
\newcommand{\bbblank}{{\equiv}}

\newcommand{\llb}{[\hspace{-1.7pt}[}
\newcommand{\rrb}{]\hspace{-1.7pt}]}

\newcommand\arxiv[2]{\href{https://arXiv.org/abs/#1}{\texttt{arXiv:\allowbreak #1} #2}}
\newcommand\doi[2]{\href{https://doi.org/#1}{#2}}

\def\centerarc[#1](#2)(#3:#4:#5);%
{
	\draw[#1]([shift=(#3:#5)]#2) arc (#3:#4:#5);
}

\tikzset{
    leftstr/.store in = \leftstr,
    leftstr = black,
    rightstr/.store in  = \rightstr,
    rightstr  = black,
    braid/.pic = {
    \centerarc[\leftstr](0.25,0)(90:180:0.25);
    \centerarc[\leftstr](0.25,0.5)(270:360:0.25);
    \centerarc[\rightstr](0.25,0)(0:45:0.25);
    \centerarc[\rightstr](0.25,0.5)(180:225:0.25);
    }; 
}

\tikzset{
    leftstr/.store in = \leftstr,
    leftstr = black,
    rightstr/.store in  = \rightstr,
    rightstr  = black,
    invbraid/.pic = {
    \centerarc[\leftstr](0.25,0)(135:180:0.25);
    \centerarc[\leftstr](0.25,0.5)(315:360:0.25);
    \centerarc[\rightstr](0.25,0)(0:90:0.25);
    \centerarc[\rightstr](0.25,0.5)(180:270:0.25);
    }; 
}

\tikzset{
    leftstr/.store in = \leftstr,
    leftstr = black,
    rightstrA/.store in  = \rightstrA,
    rightstrA  = black,
    rightstrB/.store in  = \rightstrB,
    rightstrB  = black,
    braidots/.pic = {
    \centerarc[\leftstr](0.5,0)(135:180:0.5);
    \draw[\leftstr] (0.145,0.353) -- (0.854,0.646);
    \centerarc[\leftstr](0.5,1)(315:360:0.5);
    \centerarc[\rightstrA](0,0)(0:45:0.5);
    \centerarc[\rightstrA](0.5,1)(180:225:0.5);
    \draw[\rightstrA] (0.146,0.646) -- (0.25,0.5);
    \centerarc[\rightstrB](0.5,0)(0:45:0.5);
    \centerarc[\rightstrB](1,1)(180:225:0.5);
    \draw[\rightstrB] (0.853,0.354) -- (0.75,0.5);
    }; 
}

\tikzset{
    leftstr/.store in = \leftstr,
    leftstr = black,
    rightstrA/.store in  = \rightstrA,
    rightstrA  = black,
    rightstrB/.store in  = \rightstrB,
    rightstrB  = black,
    braiduts/.pic = {
    \centerarc[\leftstr](0.5,0)(135:180:0.5);
    \draw[\leftstr] (0.4,0.453) -- (0.6,0.546);
    \centerarc[\leftstr](0.5,1)(315:360:0.5);
    \centerarc[\rightstrA](0,0)(0:45:0.5);
    \centerarc[\rightstrA](0.5,1)(180:225:0.5);
    \draw[\rightstrA] (0.146,0.646) -- (0.354,0.354);
    \centerarc[\rightstrB](0.5,0)(0:45:0.5);
    \centerarc[\rightstrB](1,1)(180:225:0.5);
    \draw[\rightstrB] (0.853,0.354) -- (0.646,0.646);
    }; 
}

\tikzset{
    leftstrA/.store in = \leftstrA,
    leftstrA = black,
    leftstrB/.store in  = \leftstrB,
    leftstrB  = black,
    rightstr/.store in  = \rightstr,
    rightstr  = black,
    braidtso/.pic = {
    \centerarc[\leftstrA](0.5,0)(135:180:0.5);
    \draw[\leftstrA] (0.146,0.353) -- (0.354,0.646);
    \centerarc[\leftstrA](0,1)(315:360:0.5);
    \centerarc[\leftstrB](1,0)(135:180:0.5);
    \draw[\leftstrB] (0.646,0.353) -- (0.854,0.646);
    \centerarc[\leftstrB](0.5,1)(315:360:0.5);
    \centerarc[\rightstr](0.5,0)(0:45:0.5);
    \draw[\rightstr] (0.6,0.45) -- (0.4,0.5);
    \centerarc[\rightstr](0.5,1)(180:225:0.5);
    }; 
}

\tikzset{
  module/.store in = \modulecolour,
  module = blue,
  object/.store in  = \objectcolour,
  object  = black,
  modbr/.pic = {
    \draw[\modulecolour] (0.5,0) -- (0.5,0.4);
    \draw[\modulecolour] (0.5,0.6) -- (0.5,1.5);
    \centerarc[\objectcolour](0.5,0)(90:180:0.5);
    \centerarc[\objectcolour](0.5,0.75)(270:430:0.25);
    \centerarc[\objectcolour](0.5,1.5)(180:260:0.5);
    \draw[gray, dashed] (-0.1,0) -- (0.85,0) -- (0.85,1.5) -- (-0.1,1.5) -- (-0.1,0);
  }
}

\tikzset{
  moduleA/.store in = \moduleAcolour,
  moduleA = blue,
  moduleB/.store in = \moduleBcolour,
  moduleB = black,
  object/.store in  = \objectcolour,
  object  = black,
  modbrots/.pic = {
    \draw[\moduleAcolour] (0.5,0) -- (0.5,0.4);
    \draw[\moduleAcolour] (0.5,0.6) -- (0.5,1.5);
    \draw[\moduleBcolour] (0,0) -- (0,0.4);
    \draw[\moduleBcolour] (0,0.6) -- (0,1.5);
    \draw[\objectcolour] (0,0.5) -- (0.5,0.5);
    \centerarc[\objectcolour](0,0)(90:180:0.5);
    \centerarc[\objectcolour](0.5,0.75)(270:430:0.25);
    \draw[\objectcolour] (0.4,1) -- (0.1,1);
    \centerarc[\objectcolour](0.0,1.5)(180:260:0.5);
    \draw[gray, dashed] (-0.6,0) -- (0.85,0) -- (0.85,1.5) -- (-0.6,1.5) -- (-0.6,0);
  }
}

\tikzset{
  module/.store in = \modulecolour,
  module = blue,
  objectA/.store in  = \objectAcolour,
  objectA  = black,
  objectB/.store in  = \objectBcolour,
  objectB  = black,
  modbrtso/.pic = {
    \draw[\modulecolour] (0.5,0) -- (0.5,0.2);
    \draw[\modulecolour] (0.5,0.3) -- (0.5,0.45);
    \draw[\modulecolour] (0.5,0.6) -- (0.5,1.5);
    \centerarc[\objectAcolour](0,0)(90:180:0.5);
    \draw[\objectAcolour] (0,0.5) -- (0.5,0.5);
    \centerarc[\objectAcolour](0.5,0.75)(270:430:0.25);
    \draw[\objectAcolour] (0.4,1) -- (0,1);
    \centerarc[\objectAcolour](0,1.5)(180:270:0.5);
    \centerarc[\objectBcolour](0.25,0)(90:180:0.25);
    \draw[\objectBcolour] (0.25,0.25) -- (0.5,0.25);
    \centerarc[\objectBcolour](0.5,0.75)(270:440:0.5);
    \draw[\objectBcolour] (0.25,1.25) -- (0.4,1.25);
    \centerarc[\objectBcolour](0.25,1.5)(180:270:0.25);
    \draw[gray, dashed] (-0.6,0) -- (1.1,0) -- (1.1,1.5) -- (-0.6,1.5) -- (-0.6,0);
  }
}

\tikzset{
    module/.store in = \modulecolour,
    module = blue,
    algebra/.store in = \algebracolour,
    algebra = green,
    objectcolour/.store in = \objectcolour,
    objectcolour = black,
    HBAmodbr/.pic = {
    \pic[rightstr = \modulecolour, leftstr = \objectcolour] at (0,0) {braid}; \pic[rightstr = \algebracolour, leftstr = \objectcolour] at (0.5,0.5) {braid}; \filldraw[gray] (0.75,0.75) circle[radius=0.1]; \draw[\algebracolour] (1,0.3) circle[radius=0.1]; \draw[\algebracolour] (1,0.4) -- (1,0.5); \centerarc[\algebracolour](0,1)(0:90:0.5); \draw[\modulecolour] (0,0.5) -- (0,1.5); \pic[leftstr = \modulecolour, rightstr = \objectcolour] at (0,1.5) {braid}; \centerarc[\objectcolour](0.75,1)(0:90:0.25); \centerarc[\objectcolour](0.75,1.5)(180:270:0.25);
    }
}

\tikzset{
  thick blue/.style = {line width=0.8pt, draw=blue, double = cyan, double distance=1.6pt}
}
\tikzset{
  thick green/.style = {line width=0.8pt, draw=olive, double=green, double distance=1.6pt}
}
\tikzset{
  thick black/.style = {line width=0.8pt, draw=black, double=lightgray, double distance=1.6pt}
}

\title{Braided module categories from small braidings}
\author{Tyler Franke$^1$\footnote{\texttt{tyler.franke@uni-hamburg.de}}, Ingo Runkel$^1$\footnote{\texttt{ingo.runkel@uni-hamburg.de}}, Volker Schomerus$^2$\footnote{\texttt{volker.schomerus@desy.de}}\\[.5cm]
$^1$\normalsize\slshape Fachbereich Mathematik, Universit\"{a}t Hamburg,\\
\normalsize\slshape Bundesstra{\ss}e 55, 20146 Hamburg, Germany \\ \\
$^2$\normalsize\slshape Deutsches Elektronen-Synchrotron DESY,\\
\normalsize\slshape Notkestra{\ss}e 85, 22607 Hamburg, Germany
}

\date{}

\allowdisplaybreaks

\begin{document}

\maketitle

\begin{abstract}
We introduce small braidings on inclusions of monoidal categories $\cat{C}\hookrightarrow\cat{D}$, which allow one to braid objects in $\cat{C}$ with objects in $\cat{D}$. 
A small braiding makes $\cat{C}$ a braided monoidal category, and turns $\cat{D}$, as well as certain subcategories of $\cat{D}$, into braided module categories over $\cat{C}$. 
We prove versions of previously known reconstruction results for braided module categories which do not rely on finiteness assumptions. 
We construct a number of new examples coming from quantum groups, as well as the associated quantum moment maps.
\end{abstract}

\newpage

\tableofcontents

\section{Introduction}

Braided monoidal categories can be thought of as monoidal categories which provide representations of the braid groups of Coxeter type $A$ \cite{JS93}. Adopting this intuition, \emph{braided module categories}, first introduced in \cite{En04}, extend these representations to representations of braid groups of type $B$ \cite{Br12}.

Braided module categories arise naturally in the Morita theory of braided monoidal categories. They are objects in the endomorphism category of the identity morphism in the 4-category of braided tensor categories
\cite{BJS18}, see \cite[Sec.\,3.3]{Jo23} for a proof of this, or objects in the centre of the category of module categories over a braided monoidal category \cite[Thm.\,4.11]{DN20}. Braided module categories have provided a categorical interpretation of quantum $K$-matrices \cite{Ko17,LWY23}, and have found applications in quantum topology \cite{BBJ16} and category theory \cite{JR21}.

\subsection{Summary}

In this paper we present a simple new construction for braided module categories. Namely, 
let $\cat{C}$ and $\cat{D}$ be monoidal categories and let there be an inclusion of monoidal categories $\cat{C}\hookrightarrow\cat{D}$. We define a \emph{small braiding} (Definition~\ref{def:SmallBraiding}) which allows one to braid objects in $\cat{C}$ with objects in $\cat{D}$. A small braiding equips $\cat{C}$ with a braiding and any full subcategory of $\cat{D}$ which is closed under tensoring with objects in $\cat{C}$ with the structure of a braided module category. An advantage of this construction is that it does not rely on finiteness or linearity assumptions for the categories in question. As far as the authors are aware, most current systematic approaches for constructing braided module categories result in a category which is the category of finite-dimensional modules over an algebra \cite{Ko17}, or locally finite-dimensional modules \cite{BBJ16} (see also \cite[Sec.\,1.3.2]{BBJ15}). Our approach also falls outside of the construction via reflective centres of \cite{LWY23}. 

Given an algebra $A\in\cat{D}$, the category $\cat{D}_A$ of (right) $A$-modules in $\cat{D}$ is a (left) $\cat{D}$-module category. Furthermore, since $\cat{C}$ is a monoidal subcategory, $\cat{D}_A$ is also a $\cat{C}$-module category. For an algebra $A\in\cat{C}$, module braidings for $\cat{C}_A$ are known to be in bijection with
half-braided structures on $A$
\cite[Sec.\,2.5]{JR21}, and with quantum moment maps when $\cat{C}$ is ribbon \cite[Cor.\,4.11]{BBJ16} 
(see Definitions~\ref{def:half-braided-algebra} and~\ref{def:QMM}). 
If $A\in\cat{D}$ is not in the image of the inclusion, we can still show one direction of the bijection:
a half-braided structure on
$A$ results in a module braiding for $\cat{D}_A$. It remains true that if $\cat{C}$ is ribbon, half-braided structures for $A$ are in bijection with quantum moment maps.

\medskip

The definition of a small braiding is motivated by observing that the braiding for quantum groups $U_q\mathfrak{g}$ with quantum parameter $q$ not a root of unity is well-defined on a tensor product of type 1 modules if at least one of the modules is finite-dimensional. We extend this result to braid generic weight modules with finite-dimensional type 1 modules by showing that there is a monoidal inclusion of the category of finite-dimensional type 1 modules for $U_q\mathfrak{g}$ into the category of weight modules for the unrolled quantum group $U_q^H\mathfrak{g}$. From this we obtain the following small braiding:

\begin{theorem}[Theorem~\ref{thm:SmallBraidingUqg}]
    Let $\cat{C}$ be the category of finite-dimensional $U_q\mathfrak{g}$-modules of type 1 and $\cat{D}^H$ the category of $U_q^H\mathfrak{g}$-weight modules. For $q\in\C^\times$ not a root of unity, there is a small braiding on the monoidal inclusion $\cat{C}\hookrightarrow\cat{D}^H$ given by
    \begin{subequations}
        \begin{equation}
            \lc_{U,X}(u\otimes x) = s_{U,X}\Bigl( q^{\sum_{i,j=1}^nd_i(A^{-1})_{ij}H_i\otimes H_j}R(u\otimes x) \Bigr)
        \end{equation}
        \begin{equation}
            \rc_{X,U}(x\otimes u) = s_{X,U}\Bigl(q^{\sum_{i,j=1}^nd_i(A^{-1})_{ij}H_i\otimes H_j}R(x\otimes u)\Bigr) \;,
        \end{equation}
    \end{subequations}
    for all $U\in\cat{C}$, $X\in\cat{D}^H$, $u\in U$ and $x\in X$, where $R$ is the $R$-matrix for $U_q\mathfrak{g}$ and $s$ is the swap operation on the tensor product of vector spaces.
\end{theorem}

Using the above theorem and basic facts about the representation theory of $U_q^H\mathfrak{g}$ for quantum parameter $q$ not a root of unity we describe several concrete examples coming from this small braiding:
\begin{itemize}
    \item (Example~\ref{ex:O}) $\cat{O}^H_q$, the quantum category $\cat{O}$ of \cite[Def.\,3.1]{AM11} for the unrolled quantum group.
    \item (Example~\ref{ex:WeylOrbit}) Given $\Lambda\in\mathfrak{h}^\ast/\sim$, where $\lambda\sim\mu$ if $\lambda-\mu$ is in the weight lattice, $\cat{D}^H_\Lambda$, the full subcategory of $\cat{D}^H$ of modules with weight support in $\Lambda$.
    \item (Example~\ref{ex:Verma}) For $\Lambda\in\mathfrak{h}^\ast/\sim$ $q$-generic, $\cat{M}^H_\Lambda$ the full subcategory of $\cat{D}^H_\Lambda$ with objects isomorphic to a finite direct sum of Verma modules with highest weights in $\Lambda$.
\end{itemize}

In the Appendix we review how internal-end objects $\underline{\cat{D}}(V,V) =: \underline{\cat{D}}(V)$ in a monoidal category $\cat{D}$, if they exist, come with the canonical structure of an algebra, and how the internal-hom functor $\underline{\cat{D}}(V,\blank):\cat{D} \to \cat{D}_{\underline{\cat{D}}(V)}$ comes with a canonical natural transformation $\Phi_{U,W}^V: U\otimes\underline{\cat{D}}(V,W)\to \underline{\cat{D}}(V,U\otimes W)$. In Lemma~\ref{lem:PhisoIffRigid} we show that $\Phi_{U,\bblank}^\blank$ is an isomorphism if and only if $U$ has a left dual. Hence $\underline{\cat{D}}(V,\blank)$ has the canonical structure of a module functor when we consider $\cat{D}$ and $\cat{D}_{\underline{\cat{D}}(V,V)}$ as module categories over a left-rigid monoidal subcategory of $\cat{D}$. In the $U_q\mathfrak{g}$ setting described above, $\cat{C}$ is a full subcategory of left-rigid objects in $\cat{D}$. It is then natural to ask if internal-homs exist in $\cat{D}$, and whether we can equip them with quantum moment maps such that $\underline{\cat{D}}(V,\blank)$ becomes a braided module functor. 

We are able to construct quantum moment maps from the $\hbar$-adic Hopf algebra $U_\hbar\mathfrak{g}$, where $q = e^\hbar$ 
and $\hbar$ is now a formal parameter instead of a complex number. Using the fact that the category of topologically free modules $\cat{D}$ is closed monoidal, we construct half-braided structures on $\underline{\cat{D}}(V)$
using a small braiding between the category $\cat{C}$ of finite rank topologically free $U_\hbar \mathfrak{g}$-modules and $\cat{D}$. Since $\cat{C}$ is ribbon, we arrive at the following theorem.
\begin{theorem}[Theorem~\ref{thm:QMMHmod}]
   Let $\cat{D}$ be the category of topologically free $U_\hbar\mathfrak{g}$-modules and $\cat{C}$ the category of finite rank topologically free $U_\hbar\mathfrak{g}$-modules. For all $V\in\cat{D}$ there is a canonical quantum moment map for the internal-endomorphism algebra $\underline{\cat{D}}(V)$ of $V$.
\end{theorem}
Because the representation theory of topologically free $U_\hbar\mathfrak{g}$-weight modules is similar to that of $U_q^H\mathfrak{g}$-weight modules for $q$ not a root of unity, we are able to work with the above examples in the $\hbar$-adic setting. In Proposition~\ref{prop:VermaQMM} we observe that the braided module category $\cat{M}_\Lambda$, which is the $\hbar$-adic analogue of $\cat{M}_\Lambda^H$ consisting now of finite direct sums of $\hbar$-adic Verma modules $V_\hbar(\lambda)$, satisfies sufficient conditions (Proposition~\ref{prop:FFMCEmbedding}) such that the functor $\underline{\cat{D}}(V_\hbar(\lambda),\blank):\cat{M}_\Lambda \to \cat{D}_{\underline{\cat{D}}(V_\hbar(\lambda))}$ is a fully faithful (but not essentially surjective) braided $\cat{C}$-module functor.

\subsection{Outlook}\label{sec:outlook}

Our main motivation comes from an application in defect skein theory. Namely, organising defect data of skein modules according to the Morita theory of braided tensor categories \cite{BJS18}, one is led to label line defects with braided module categories and junctions between line defects with braided module functors. 
In upcoming work, we will study examples of defect skein modules with defect data coming from the examples discussed here \cite{paper2-skeins}.
See \cite{BBJ16,BJ25,G26,KR26,BGJV26,CJR26} for previous work of a similar flavour on defect skein modules.

In Remark~\ref{rem:Not an equivalence} we show that the inclusion of braided module categories $\cat{M}_\Lambda\hookrightarrow\cat{D}_{\underline{\cat{D}}(V_\hbar(\lambda))}$ is not an equivalence. Tools for monadic reconstruction of module categories were developed in \cite{SZ24} precisely in order to deal with non-finite and non-rigid settings. Any $V_\hbar(\lambda)\in\cat{M}_\Lambda$ satisfies similar properties to the $\cat{C}$-projective $\cat{C}$-generator requirement for the reconstruction theorem \cite[Thm.\,5.1]{SZ24}. It is then natural to ask whether there is a monadic version of a quantum moment map, allowing one to reconstruct $\cat{M}_\Lambda$.

The quantum group setting is amenable to constructing many more examples of braided module categories from small braidings. One could study the tensor product of finite-dimensional modules with generalised Verma modules in \cite[Prop.\,3.7]{Mu04} to construct braided module categories similar to $\cat{M}_\Lambda$. These examples should also have counterparts when the quantum parameter is specialised to a complex number. Examples beyond highest weight modules could be constructed by studying the tensor product of the torsion-free modules of \cite{FHW14} with finite-dimensional modules.

\subsection{Physics context and motivation}\label{subsec:physics}
 
Many of the structures in this paper have a long history in mathematical physics, where they often appeared well before their categorical formulation.
The purpose of this subsection is twofold. First, it collects the relevant physics literature in one place, so that the exposition in the main text does not have to be interrupted by these references. Second, it is meant to help readers with a background in physics to follow the constructions in the main text, by relating the categorical notions we use to their counterparts in
physics.
 
Ribbon categories such as $\mathcal{C}$, and the skein theory built from them, are familiar from Chern--Simons theory and rational conformal field theory, where they describe Wilson lines and the braiding of chiral vertex
operators~\cite{Wi89,MS89}. 
The additional structure of a braided module category
has an equally long history in integrable systems with boundaries. Its defining relation in the braid group of type~$B$, $\tau\sigma_1\tau\sigma_1 = \sigma_1\tau\sigma_1\tau$, 
is the reflection equation, which was introduced by Cherednik for the factorised scattering of particles on a half-line~\cite{Ch84} and by Sklyanin for integrable spin chains with boundaries~\cite{Sk88}. 
The module braiding plays the role of the boundary $K$-matrix, and the spectral-parameter-independent reflection equation relevant here was studied, together with its relation to quantum groups, in~\cite{KS92}. In this language, the axioms \eqref{eq:BrM1} and \eqref{eq:BrM2} satisfied by the module braiding express the dressing of a boundary by an additional site~\cite{Sk88} and the fusion of
$K$-matrices~\cite{MN92}.
 
The coend $\mathcal{L}$, quantum moment maps (cf.\ Definition~\ref{def:QMM}) and internal endomorphism algebras 
(cf.\ Theorem~\ref{thm:QMMHmod}) appeared in the combinatorial quantisation of Chern--Simons theory~\cite{AGS94a,AGS94b,BR94}. There, $\mathcal{L}$ is the algebra generated by the
holonomy around a single puncture. Its defining relations are those of the
reflection equation algebra of~\cite{KS92}, and this loop algebra is closely
related to the quantum group itself~\cite{Al93}, via the map $\varphi$ of
Remark~\ref{rem:Mudrov}.  A quantum moment map $\mu\colon\mathcal{L}\to A$ realises 
the holonomy in an algebra~$A$, and the moment map condition \eqref{eq:QMM}
expresses that conjugation by the holonomy implements the quantum gauge symmetry. For a puncture carrying a representation~$V$, the holonomy acts on $V$ through $R_{21}R$~\cite{AS95}, which is formula
\eqref{eq:QMMexplicit},
and quantum Hamiltonian reduction along such maps underlies the moduli algebras of~\cite{AGS94b}. Classically, these structures quantise the group-valued moment maps of~\cite{AMM97}. 
The images $\mu(\mathcal{L})\subset \underline{\cat{D}}(V)$, 
i.e.\ quantised conjugacy classes, arose independently in string theory as world-volume algebras of D-branes in WZW models, which are localised along conjugacy classes and become fuzzy spheres at large level~\cite{ARS99}. 
More generally, internal endomorphism algebras organise boundary fields in the topological field theory construction of conformal field theory correlators~\cite{FRS02}.
 
Finally, our main motivation, the use of braided module categories $\cat{M}$
as labels of line defects in skein theory (Section~\ref{sec:outlook}), has a
natural counterpart in gauge theory. There, defects of codimension two are
described by the monodromy of the gauge field around them. In four dimensions,
Gukov and Witten introduced surface operators through a singularity of the
gauge field that prescribes this monodromy, with continuous parameters that fix
its conjugacy class~\cite{GW06}. In three-dimensional Chern--Simons theory,
Wilson lines can be described in the same way, but their monodromy is quantised
by the integral highest weight of the representation~\cite{Wi11}. In skein
theory, the monodromy around a defect line carrying $M\in\cat{M}$ is encoded
in the module braiding $e_{U,M}$, i.e.\ in the effect of a Wilson line labelled
by $U\in\mathcal{C}$ winding once around it. For $\cat{M}=\mathcal{C}$ this
is the double braiding (Remark~\ref{rem:bmc-definitions-compare}), and the
defect lines are ordinary Wilson lines. Verma modules of generic weight, on the
other hand, label defect lines whose monodromy depends on the continuous
parameter~$\lambda$, as for the surface operators of Gukov and Witten. The same
distinction between discrete and continuous labels appears in Liouville theory
with $q=e^{i\pi b^2}$: degenerate fields act on generic primaries $V_\alpha$ as
finite-dimensional modules act on the Verma modules of Example~\ref{ex:Verma},
with $\lambda+1=(Q-2\alpha)/b$, and their monodromy around $V_\alpha$ is given
by the module braiding. That the fusion of a degenerate field with $V_\alpha$
has only finitely many channels is the key input in Teschner's derivation of
the three-point function~\cite{Te95}.
For $\mathfrak{g}=\mathfrak{sl}_2$, the category $\cat{M}^H_\Lambda$ also describes the
$U_q\mathfrak{sl}_2$-invariant open XXZ chain~\cite{PS90} with a Verma module
at the boundary~\cite{CGS22}, whose commutant is the blob algebra~\cite{MS93}.
In the morphism spaces of $\cat{M}_\Lambda^H$ one can find the morphism spaces 
of the blob category of \cite{ILZ18,LNV20} and hence also the blob algebras of \cite{MS93}.

\subsection{Organisation}

In Section~\ref{sec:BrModCat} we review the basics of braided module categories and braided module functors. Then in Section~\ref{ssec:HBAs and QMMs} we review the construction of braided module categories via half-braided algebras and quantum moment maps. We then introduce small braidings in Section~\ref{ssec:SmallBraidings}, and establish generalised versions of half-braided algebras and quantum moment maps. In Section~\ref{sec:BrModCatsforQuantumGroups} we study braided module categories coming from quantum groups. We begin by discussing the small braiding between finite-dimensional type 1 modules for $U_q\mathfrak{g}$ and weight modules for $U_q^H\mathfrak{g}$ in Section~\ref{ssec:Small braidings for Uqg}. We then introduce concrete examples in this setting in Section~\ref{ssec:Examples for generic q}. We conclude in Section~\ref{ssec:QMMs from internal-homs} by switching settings to the $\hbar$-adic quantum group, where we construct half-braided algebras and quantum moment maps using internal-homs. In Appendix~\ref{app:closed} we prove some results about algebras in closed categories which we use in Section~\ref{ssec:QMMs from internal-homs}.

\subsection*{Acknowledgements}

We thank Kevin Walker, Patrick Kinnear, Christopher Raymond, Mateusz Stroiński and Tony Zorman for discussions. All authors acknowledge support from the Deutsche Forschungsgemeinschaft (DFG, German Research Foundation) – SFB 1624 – ``Higher structures, moduli spaces and integrability'' – 506632645; IR and VS are in addition supported under Germany's Excellence Strategy - EXC 2121 ``Quantum Universe'' - 390833306.

\subsection*{AI usage disclosure}

During the preparation of this manuscript we used Claude and ChatGPT to search for relevant literature, as well as to proof-read a draft during the final editing. None of the text or diagrams in this manuscript are generated by AI. We hold full responsibility for the validity of the results in this manuscript.

\newpage

\section{Categorical preliminaries}

In this section we review the notion of a braided module category over a braided monoidal category and its corresponding graphical calculus. We then review the results of \cite{BBJ16} and \cite{JR21} which construct braided module categories from so-called quantum moment maps and half-braided algebras respectively. Finally, we define small braidings, and develop the analogues of quantum moment maps and half-braided algebras in this new setting.

Given a category $\cat{C}$ and $U,V \in \cat{C}$, we denote the space of morphisms $U\to V$ as $\cat{C}(U,V)$. Our conventions for a braided monoidal category $(\cat{C},\otimes,\1,\alpha,l,r,\sigma)$ follow those of \cite[Chs.~XI-XIII]{Ka95}. Namely, the direction of the associator is
\begin{equation}
    \alpha: (\blank\otimes\bblank)\otimes\bblank \Rightarrow \blank\otimes(\bblank\otimes\bblank)
    \ .
\end{equation}
Our convention and notation for a rigid monoidal category is the same as that of \cite[Sec.\,2.10]{EGNO15}. We define a ribbon category as in \cite[Def.\,8.10.1]{EGNO15}. We abbreviate the identity morphism on objects $\Id_V$ by simply writing the object $V$ when there is no ambiguity in doing so.

We frequently employ a diagrammatic calculus for braided and rigid monoidal categories. Diagrams are to be read from the bottom up, with the structural morphisms being depicted as
\begin{equation}
    \sigma_{U,V} = \begin{tikzpicture}[baseline={(0,0.25)}]
        \pic at (0,0) {braid}; \draw[-Latex] (0,0.5) -- (0,0.75); \draw[-Latex] (0.5,0.5) -- (0.5,0.75); \node at (0,-0.2) {$U$}; \node at (0.5,-0.2) {$V$}; \node at (0,0.95) {$V$}; \node at (0.5,0.95) {$U$};
    \end{tikzpicture}\;, \quad \ev_V = \begin{tikzpicture}[baseline={(0,0.1)}] \centerarc[](0,0)(0:180:0.5); \draw[-Latex] (0,0.5) -- (-0.2,0.5); \node at (-0.5,-0.2) {$V^\ast$}; \node at (0.5,-0.2) {$V$}; \end{tikzpicture} \;,\quad \coev_V = \begin{tikzpicture}[baseline={(0,-0.2)}] \centerarc[](0,0)(180:360:0.5); \draw[-Latex] (0,-0.5) -- (-0.2,-0.5); \node at (-0.5,0.2) {$V$}; \node at (0.5,0.2) {$V^\ast$}; \end{tikzpicture} \;,
\end{equation}
and if our category has pivotal structure $\delta_V:V\to V^{\ast\ast}$,
\begin{equation}
    \ev_{V^\ast}\circ(\delta_V\otimes V^\ast) = \begin{tikzpicture}[baseline={(0,0.1)}] \centerarc[](0,0)(0:180:0.5); \draw[-Latex] (0,0.5) -- (0.2,0.5); \node at (-0.5,-0.2) {$V$}; \node at (0.5,-0.2) {$V^\ast$}; \end{tikzpicture} \;,\quad (V^\ast\otimes\delta_V^{-1})\circ\coev_{V^\ast} = \begin{tikzpicture}[baseline={(0,-0.2)}] \centerarc[](0,0)(180:360:0.5); \draw[-Latex] (0,-0.5) -- (0.2,-0.5); \node at (-0.5,0.2) {$V^\ast$}; \node at (0.5,0.2) {$V$}; \end{tikzpicture} \;.
\end{equation}

\subsection{Braided module categories}\label{sec:BrModCat}

Let $\cat{C}$ be a braided monoidal category.
\begin{definition}[{\cite[Def.\,7.1.2]{EGNO15}}]
    A \emph{(left) $\cat{C}$-module category} is a category $\cat{M}$ along with an \emph{action bifunctor} $\blank\triangleright\bblank:\cat{C}\times\cat{M}\to \cat{M}$ and a \emph{modulator}, which is a natural isomorphism
    \begin{equation}
        \alpha^\cat{M} : (\blank\otimes \bblank) \triangleright \bbblank \Rightarrow \blank\triangleright(\bblank\triangleright \bbblank) \;,
    \end{equation}
    as well as a unitor isomorphism $l^\cat{M} : \1\triangleright \blank \Rightarrow \Id_\cat{M}$. 
    The modulator is required to satisfy the mixed pentagon identity, and the unitor the mixed triangle identity, both of which can be found in \cite[Sec.\,7.1]{EGNO15}.
\end{definition}
\begin{example}[{\cite[Sec.\,7.8]{EGNO15}}]
Let $(A,m,u)$ be an algebra in $\cat{C}$. The category $\cat{C}_A$ of right $A$-modules in $\cat{C}$ is a left $\cat{C}$-module category.
\end{example}
\begin{definition}[{\cite[Def.\,7.2.1]{EGNO15}}]
    A \emph{$\cat{C}$-module functor $F:(\cat{M},\triangleright)\to(\cat{N},\triangleright^\prime)$} is a functor $F:\cat{M}\to\cat{N}$ together with a natural isomorphism $F_{\blank,\bblank}:\blank\triangleright^\prime F(\bblank) \Rightarrow F(\blank\triangleright \bblank)$ such that the following diagrams commute for all objects $U,V$ in $\cat{C}$ and $M$ in $\cat{M}$:
    \begin{subequations}
    \begin{equation}
        \begin{tikzcd}
            & (U\otimes V)\triangleright^\prime F(M) & \\ U\triangleright^\prime(V\triangleright^\prime
            F(M)) & & F((U\otimes V)\triangleright M) \\ U\triangleright^\prime F(V \triangleright M) & & F(U\triangleright(V\triangleright M))
            \arrow["{F_{U\otimes V,M}}", from=1-2 , to=2-3 ]\arrow["{F(\alpha^\cat{M})}", from=2-3 , to=3-3 ]
            \arrow["{\alpha^\cat{N}}"', from=1-2 , to=2-1 ]\arrow["{U\triangleright^\prime F_{V,M}}"', from=2-1 , to=3-1 ]\arrow["{F_{U,V\triangleright M}}", from=3-1 , to=3-3 ]
        \end{tikzcd} \;,
    \end{equation}
    \begin{equation}
        \begin{tikzcd}
            \1\triangleright^\prime F(M) && F(\1\triangleright M) \\
            &F(M)&
            \arrow["{F_{\1,M}}", from=1-1, to=1-3]\arrow["{F(l^\cat{M}_M)}", from=1-3, to=2-2]
            \arrow["{l^\cat{N}_{F(M)}}"', from=1-1, to=2-2]
        \end{tikzcd}
    \end{equation}
    \end{subequations}
\end{definition}
\begin{example}[{\cite[Sec.\,7.8]{EGNO15}}]
    Assume $\cat{C}$ is 
    in addition abelian and that the tensor product is right exact. Let $(A,m,u)$ and $(B,m^\prime,u^\prime)$ be algebras in $\cat{C}$ and $M$ be an $(A,B)$-bimodule in $\cat{C}$. The functor $\blank\otimes_A M:\cat{C}_A\to\cat{C}_B$ is a $\cat{C}$-module functor.
\end{example}
\begin{definition}
    A \emph{(left) braided $\cat{C}$-module category} is a $\cat{C}$-module category $(\cat{M},\triangleright$, $\alpha^\cat{M},l^\cat{M})$ along with a \emph{module braiding}, which is a natural automorphism of the action functor
    \begin{equation}
        e_{\blank,\bblank}: \blank\triangleright \bblank \rightarrow \blank\triangleright \bblank \;,
    \end{equation}
    such that the following diagrams commute for all $M\in\cat{M}$ and $U,V\in\cat{C}$:
    \begin{subequations}
            \begin{equation}\label{eq:BrM1}
                \begin{tikzcd}
                    (U\otimes V)\triangleright M
                        & (V\otimes U)\triangleright M \\
                    U\triangleright(V\triangleright M)
                        & V\triangleright(U\triangleright M) \\
                    U\triangleright(V\triangleright M)
                        & V\triangleright(U\triangleright M) \\
                    (U\otimes V)\triangleright M
                        & (V\otimes U)\triangleright M
                    \arrow["{\sigma_{U,V}\triangleright M}",
                        from=1-1, to=1-2]
                    \arrow["{\alpha^\cat{M}}"',
                        from=1-1, to=2-1]
                    \arrow["{e_{U,V\triangleright M}}"',
                        from=2-1, to=3-1]
                    \arrow["{(\alpha^\cat{M})^{-1}}"',
                        from=3-1, to=4-1]
                    \arrow["{\alpha^\cat{M}}",
                        from=1-2, to=2-2]
                    \arrow["{V\triangleright e_{U,M}}",
                        from=2-2, to=3-2]
                    \arrow["{(\alpha^\cat{M})^{-1}}",
                        from=3-2, to=4-2]
                    \arrow["{\sigma_{V,U}\triangleright M}"',
                        from=4-2, to=4-1]
                \end{tikzcd}
            \end{equation}
            \begin{equation}\label{eq:BrM2}
                \begin{tikzcd}
                    (U\otimes V)\triangleright M && (U\otimes V)\triangleright M \\ (V\otimes U)\triangleright M && U\triangleright(V\triangleright M) \\ V\triangleright(U\triangleright M) && U\triangleright(V\triangleright M) \\ V\triangleright(U\triangleright M) && (U\otimes V)\triangleright M \\ & (V\otimes U)\triangleright M &
                    \arrow["{e_{U\otimes V,M}}", from=1-1, to=1-3]
                    \arrow["{\sigma^{-1}_{V,U}\triangleright M}"', from=1-1, to=2-1] \arrow["\alpha^\cat{M}"', from=2-1, to=3-1] \arrow["{V\triangleright e_{U,M}}"', from=3-1, to=4-1] \arrow["{(\alpha^\cat{M})^{-1}}", from=4-1, to=5-2] \arrow["{\sigma_{V,U}\triangleright M}", from=5-2, to=4-3] \arrow["{\alpha^\cat{M}}"', from=4-3, to=3-3] \arrow["{U\triangleright e_{V,M}}"', from=3-3, to=2-3] \arrow["{(\alpha^\cat{M})^{-1}}"', from=2-3, to=1-3]
                \end{tikzcd} 
            \end{equation}
    \end{subequations}
\end{definition}

We note that naturality of the associativity and unitor isomorphisms and condition \eqref{eq:BrM2} for $U=V=\1$ imply $e_{\1,M} = e_{\1,M}\circ e_{\1,M}$. Since $e_{\1,M}$ is an isomorphism, this gives $e_{\1,M} = \Id_{\1 \triangleright M}$.

\begin{remark}\label{rem:bmc-definitions-compare}
    Braided module categories were first introduced in \cite{En04}. Several different conventions have emerged regarding their definition, the complete list (to the authors' knowledge) being \cite[Sec.\,4.3]{En04}, \cite[Sec.\,5.1]{Br12}, \cite[Def.\,3.4]{BBJ16}, \cite[Def.\,2.1]{Ko17}, \cite[Def.\,4.1]{DN20}, \cite[Sec.\,2.2]{JR21} and \cite[Def.\,3.5]{LWY23}. We choose the same convention as \cite[Def.\,3.5]{LWY23}, allowing for $\cat{C}$ to be a braided module category over itself with module braiding $e_{U,V} = \sigma_{V,U}\circ\sigma_{U,V}$. Every definition listed is equivalent when the appropriate assumptions are made on $\cat{C}$, see \cite[Rem.\,3.15]{LWY23} for a detailed comparison. 
\end{remark}
\begin{remark}
    It is a standard result that any monoidal category is monoidally equivalent to a strict monoidal category (see for example \cite[Thm.\,2.8.5]{EGNO15}). A similar result applies for module categories: A $\cat{C}$-module category is equivalent to a module category over the strictification of $\cat{C}$ such that the mixed associator and unitor are identities \cite[Sec.\,5.1]{Br12}. We call such a module category a \emph{strict module category}. 
\end{remark}
\begin{definition}
    A $\cat{C}$-module functor $F:(\cat{M},\triangleright)\to(\cat{N},\triangleright^\prime)$ of braided module categories $(\cat{M},e)$ and $(\cat{N},e^\prime)$ is \emph{braided} if the diagram
    \begin{equation}\label{eq:BrModFunct}
        \begin{tikzcd}
            V\triangleright^\prime F(M) & V\triangleright^\prime F(M) \\ F(V\triangleright M) & F(V\triangleright M)
            \arrow["{e^\prime_{V,F(M)}}", from=1-1, to=1-2] \arrow["{F_{V,M}}", from=1-2, to=2-2]
            \arrow["{F_{V,M}}"', from=1-1, to=2-1] \arrow["{F(e_{V,M})}", from=2-1, to=2-2]
        \end{tikzcd}
    \end{equation}
    commutes for all
    $V \in \cat{C}$, $M \in \cat{M}$.
\end{definition}

Let $\cat{M}$ be a braided $\cat{C}$-module category. We use the following graphical notation for the module braiding.
\begin{equation}
    e_{V,M} = \begin{tikzpicture}[baseline={(0,0.8)}] \pic at (0,0) {modbr}; \draw (0,-0.3) -- (0,0); \draw[-Latex] (0,1.5) -- (0,2); \draw[blue] (0.5,-0.3) -- (0.5,0); \draw[blue,-Latex] (0.5,1.5) -- (0.5,2); \node at (0,-0.5) {$V$}; \node[blue] at (0.5,-0.5) {$M$}; \node at (0,2.2) {$V$}; \node[blue] at (0.5,2.2) {$M$}; \end{tikzpicture}
\end{equation}
The graphical notation for a braided module category implicitly uses the equivalence to a strict braided module category over a strict braided monoidal category. The conditions \eqref{eq:BrM1} and \eqref{eq:BrM2} are quite natural in the graphical notation,
\begin{subequations}
    \begin{equation}\label{eq:BrMpicturesA}
          \begin{tikzpicture}[baseline={(0,0.8)}]
             \pic at (0,0) {modbrots}; \draw (-0.5,-0.5) -- (-0.5,0); \node at (-0.5,-0.7) {$U$}; \draw (0,-0.5) -- (0,0); \node at (0,-0.7) {$V$}; \draw[blue] (0.5,-0.5) -- (0.5,0); \node[blue] at (0.5,-0.7) {$M$}; \draw[-Latex] (0,1.5) -- (0,2); \node at (-0.5,2.2) {$U$}; \draw[-Latex] (-0.5,1.5) -- (-0.5,2); \node at (0,2.2) {$V$}; \draw[-Latex, blue] (0.5,1.5) -- (0.5,2); \node[blue] at (0.5,2.2) {$M$};
        \end{tikzpicture} 
        ~= 
        \begin{tikzpicture}[baseline={(0,0.8)}]
            \pic at (0,0) {modbr}; \pic at (-0.5,1.5) {braid}; \pic at (-0.5,-0.5) {braid}; \draw (-0.5,0) -- (-0.5,1.5); \draw[-Latex] (-0.5,2) -- (-0.5,2.25); \draw[-Latex] (0,2) -- (0,2.25); \draw[blue] (0.5,-0.5) -- (0.5,0); \draw[blue,-Latex] (0.5,1.5) -- (0.5,2.25); \node at (-0.5,-0.7) {$U$}; \node at (0,-0.7) {$V$}; \node[blue] at (0.5,-0.7) {$M$}; \node at (-0.5,2.45) {$U$}; \node at (0,2.45) {$V$}; \node[blue] at (0.5,2.45) {$M$};
        \end{tikzpicture}
    \end{equation}
    \begin{equation}\label{eq:BrMpicturesB}
        \begin{tikzpicture}[baseline={(0,0.8)}]
             \pic at (0,0) {modbrtso}; \draw[-Latex] (-0.5,1.5) -- (-0.5,2); \draw[-Latex] (0,1.5) -- (0,2); \draw[blue] (0.5,-0.5) -- (0.5,0); \draw[blue,-Latex] (0.5,1.5) -- (0.5,2); \draw (0,-0.5) -- (0,0); \draw (-0.5,-0.5) -- (-0.5,0); \node at (-0.5,-0.7) {$U$}; \node at (0,-0.7) {$V$}; \node[blue] at (0.5,-0.7) {$M$}; \node at (-0.5,2.2) {$U$}; \node at (0,2.2) {$V$}; \node[blue] at (0.5,2.2) {$M$};
        \end{tikzpicture}
        =
        \begin{tikzpicture}[baseline={(0,0.8)}]
            \pic at (-0.5,-0.5) {invbraid}; \pic at (0,0) {modbr};  \pic at (-0.5,1.5) {braid}; \pic at (0,2) {modbr}; \node at (-0.5,-0.7) {$U$}; \draw (-0.5,0) -- (-0.5,1.5); \draw[-Latex] (-0.5,2) -- (-0.5,4); \node at (-0.5,4.2) {$U$}; \draw[-Latex] (0,3.5) -- (0,4); \node at (0,4.2) {$V$}; \node at (0,-0.7) {$V$}; \draw[blue] (0.5,-0.5) -- (0.5,0); \draw[blue] (0.5,1.5) -- (0.5,2); \draw[blue,-Latex] (0.5,3.5) -- (0.5,4); \node[blue] at (0.5,-0.7) {$M$}; \node[blue] at (0.5,4.2) {$M$};
        \end{tikzpicture} \;.
    \end{equation}
\end{subequations}

\subsection{Half-braided algebras and quantum moment maps}\label{ssec:HBAs and QMMs}

Let $\cat{C}$ be a braided monoidal category. 

In \cite{JR21} the theory of half-braided algebras, including the results of this section, was developed under the assumption that $\cat{C}$ is a fusion category. However, neither semisimplicity, rigidity, nor finiteness is required in the proofs of the statements we use. We include the proofs below as we refer to them in the sequel.

\begin{proposition}[{\cite[Lemma~2.38]{JR21}}]\label{Prop:BrModCatiffHalfBr}
    Let $(A,m,u)$ be an algebra in $\cat{C}$. There is a bijective correspondence between module braidings for $\cat{C}_A$ and half-braidings 
$\gamma^A:\blank\otimes A \xRightarrow{\sim} A\otimes\blank$, depicted
\begin{equation}
    \gamma^A_V = 
    \begin{tikzpicture}[baseline={(0,0.25)}]
        \pic[rightstr=green] at (0,0) {braid}; \filldraw[gray] (0.25,0.25) circle[radius=0.1]; \draw[green,-Latex] (0,0.5) -- (0,0.75); \draw[-Latex] (0.5,0.5) -- (0.5,0.75); \node at (0,-0.2) {$V$}; \node[green] at (0.5,-0.2) {$A$}; \node[green] at (0,0.95) {$A$}; \node at (0.5,0.95) {$V$};
    \end{tikzpicture} \;, \quad (\gamma^A_V)^{-1} = 
    \begin{tikzpicture}[baseline={(0,0.25)}]
        \pic[leftstr=green] at (0,0) {invbraid}; \filldraw[gray] (0.25,0.25) circle[radius=0.1]; \draw[-Latex] (0,0.5) -- (0,0.75); \draw[green,-Latex] (0.5,0.5) -- (0.5,0.75); \node[green] at (0,-0.2) {$A$}; \node at (0.5,-0.2) {$V$}; \node at (0,0.95) {$V$}; \node[green] at (0.5,0.95) {$A$};
    \end{tikzpicture} \;,
\end{equation}
satisfying the hexagon \cite[Eq.\,(7.41)]{EGNO15}, such that the diagram
    \begin{equation}\label{eq:HBAdiagram}
        \begin{tikzcd}[column sep=1.2]
            & (A\otimes V)\otimes A && A\otimes(A\otimes V) & \\ (V\otimes A)\otimes A && V\otimes A && A\otimes V \\ & (A\otimes V)\otimes A && A\otimes(A\otimes V) &
            \arrow["\gamma^A_V\otimes A", from=2-1, to=1-2]\arrow["(A\otimes\sigma_{A,V}^{-1})\circ\alpha", from=1-2, to=1-4]\arrow["(m\otimes V)\circ\alpha^{-1}", from=1-4, to=2-5]
            \arrow["(V\otimes m)\circ\alpha", from=2-1, to=2-3]\arrow["\gamma^A_V", from=2-3, to=2-5]
            \arrow["\sigma_{V,A}\otimes A"', from=2-1, to=3-2]\arrow["(A\otimes\gamma^A_V)\circ\alpha"', from=3-2, to=3-4]\arrow["(m\otimes V)\circ\alpha^{-1}"', from=3-4, to=2-5]
        \end{tikzcd}
    \end{equation}
    commutes for all objects $V$ in $\cat{C}$. The condition \eqref{eq:HBAdiagram} is depicted
    \begin{equation}\label{eq:HBA}
        \begin{tikzpicture}[baseline={(0,0.8)}]
            \centerarc[green](0,0.3)(0:180:0.5); \draw[green] (-0.5,0) -- (-0.5,0.3); \draw[green] (0.5,0) -- (0.5,0.3);
            \draw[green] (0,0.8) -- (0,1); \pic[rightstr = green] at (-0.5,1) {braid};             \filldraw[gray] (-0.25,1.25) circle[radius=0.1]; \draw[green, -Latex] (-0.5,1.5) -- (-0.5,1.8); \draw[-Latex] (0,1.5) -- (0,1.8);
            \centerarc[](-0.75,1)(270:360:0.25); \centerarc[](-0.75,0.5)(90:180:0.25); \draw (-1,0) -- (-1,0.5);
            \node at (-1,-0.2) {$V$}; \node[green] at (-0.5,-0.2) {$A$}; \node[green] at (0.5,-0.2) {$A$};
            \node at (0,2) {$V$}; \node[green] at (-0.5,2) {$A$};
        \end{tikzpicture}
        =
        \begin{tikzpicture}[baseline={(0,0.8)}]
            \pic[rightstr = green] at (-0.5,0) {braid}; \filldraw[gray] (-0.25,0.25) circle[radius=0.1]; \pic[rightstr = green] at (0,0.5) {invbraid};
            \node at (-0.5,-0.2) {$V$}; \node[green] at (0,-0.2) {$A$}; \node[green] at (0.5,-0.2) {$A$};
            \draw[-Latex] (0.5,1) -- (0.5,1.8); \node at (0.5,2) {$V$}; \draw[green] (-0.5,0.5) -- (-0.5,1); \draw[green] (0.5,0) -- (0.5,0.5); \centerarc[green](-0.25,1)(0:180:0.25); \draw[green, -Latex] (-0.25,1.25) -- (-0.25,1.8); \node[green] at (-0.25,2) {$A$};
        \end{tikzpicture}
        =
        \begin{tikzpicture}[baseline={(0,0.8)}]
            \pic[rightstr = green] at (-0.5,0) {braid}; \pic[rightstr = green] at (0,0.5) {braid}; \filldraw[gray] (0.25,0.75) circle[radius=0.1];
            \node at (-0.5,-0.2) {$V$}; \node[green] at (0,-0.2) {$A$}; \node[green] at (0.5,-0.2) {$A$};
            \draw[-Latex] (0.5,1) -- (0.5,1.8); \node at (0.5,2) {$V$}; \draw[green] (-0.5,0.5) -- (-0.5,1); \draw[green] (0.5,0) -- (0.5,0.5); \centerarc[green](-0.25,1)(0:180:0.25); \draw[green, -Latex] (-0.25,1.25) -- (-0.25,1.8); \node[green] at (-0.25,2) {$A$};
        \end{tikzpicture} \;.
    \end{equation}
    For all $U,V\in\cat{C}$, the hexagon \cite[Eq.\,(7.41)]{EGNO15} is depicted
    \begin{equation}\label{eq:half-hexagon}
        \begin{tikzpicture}[baseline={(0,0.8)}]
            \draw (0,0) -- (0,0.2);
            \draw (0.5,0) -- (0.5,0.2);
            \draw[green] (1,0) -- (1,0.2);
            \pic[rightstr=green] at (0,0.2) {braidtso};
            \draw[green,-Latex] (0,1.2) -- (0,1.7);
            \draw[-Latex] (0.5,1.2) -- (0.5,1.7);
            \draw[-Latex] (1,1.2) -- (1,1.7);
            \node[green] at (1,-0.2) {$A$}; \node[green] at (0,1.9) {$A$};
            \node at (0,-0.2) {$U$}; \node at (0.5,1.9) {$U$};
            \node at (0.5,-0.2) {$V$}; \node at (1,1.9) {$V$};
            \filldraw[gray] (0.5,0.7) circle[radius=0.3];
        \end{tikzpicture}
        =
        \begin{tikzpicture}[baseline={(0,0.8)}]
            \draw (0,0) -- (0,0.7);
            \draw (0.5,0) -- (0.5,0.2);
            \draw[green] (1,0) -- (1,0.2);
            \pic[rightstr = green] at (0.5,0.2) {braid}; \filldraw[gray] (0.75,0.45) circle[radius=0.1];
            \pic[rightstr = green] at (0,0.7) {braid}; \filldraw[gray] (0.25,0.95) circle[radius=0.1];
            \draw[green,-Latex] (0,1.2) -- (0,1.7);
            \draw[-Latex] (0.5,1.2) -- (0.5,1.7);
            \draw[-Latex] (1,0.7) -- (1,1.7);
            \node[green] at (1,-0.2) {$A$}; \node[green] at (0,1.9) {$A$};
            \node at (0,-0.2) {$U$}; \node at (0.5,1.9) {$U$};
            \node at (0.5,-0.2) {$V$}; \node at (1,1.9) {$V$};
        \end{tikzpicture} \;.
    \end{equation}
\end{proposition}
\begin{proof}
    The proof is essentially identical to that of \cite[Lemma~2.38]{JR21}, although our convention differs from theirs. We reproduce it string-diagrammatically. We show via a direct computation that a module braiding can be obtained from a half-braiding satisfying \eqref{eq:HBA} via
    \begin{equation}\label{eq:MBrfromHBr}
        e_{V,M} := \begin{tikzpicture}[baseline={(0,0.8)}]
            \pic at (0,0) {HBAmodbr}; 
            \draw (0,-0.3) -- (0,0); \node at (0,-0.5) {$V$}; \draw[-Latex] (0,2) -- (0,2.5); \node at (0,2.7) {$V$};
            \draw[blue] (0.5,-0.3) -- (0.5,0); \node[blue] at (0.5,-0.5) {$M$}; \draw[blue,-Latex] (0.5,2) -- (0.5,2.5); \node[blue] at (0.5,2.7) {$M$};
        \end{tikzpicture} \;,
    \end{equation}
    and a half-braiding satisfying \eqref{eq:HBA} from a module braiding via
    \begin{equation}\label{eq:HbrfromMBr}
        \gamma^A_V := \begin{tikzpicture}[baseline={(0,0.8)}]
            \pic[module = green] at (0,0) {modbr}; \pic[rightstr = green] at (0,1.5) {invbraid}; \draw (0,-0.3) -- (0,0); \node at (0,-0.5) {$V$}; \draw[green] (0.5,-0.3) -- (0.5,0); \node[green] at (0.5,-0.5) {$A$}; \draw[green,-Latex] (0,2) -- (0,2.5); \node[green] at (0,2.7) {$A$}; \draw[-Latex] (0.5,2) -- (0.5,2.5); \node at (0.5,2.7) {$V$};
        \end{tikzpicture}\;.
    \end{equation} 
Showing that the above assignments are indeed mutually inverse is a straightforward computation using the properties of half-braidings and module braidings. 

    We now show that \eqref{eq:MBrfromHBr} defines a module braiding.
    Let $U,V$ be objects in $\cat{C}$ and $(M,\rho)$ an object in $\cat{C}_A$. Note that $V\triangleright M = (V\otimes M,V\otimes \rho)$ 
    in the strict setting. A straightforward computation using \eqref{eq:HBA} shows that the inverse for $e_{\blank,\bblank}$ defined in \eqref{eq:MBrfromHBr} is
    \begin{equation}
            e_{V,M}^{-1}  = \begin{tikzpicture}[baseline={(0,0.8)}]
            \node at (0,-0.2) {$V$}; \node[blue] at (0.5,-0.2) {$M$};
            \draw (0,0) -- (0,0.3); \draw[blue] (0.5,0) -- (0.5,0.3);
            \pic[rightstr=blue] at (0,0.3) {invbraid};
            \centerarc[](0.75,0.8)(90:180:0.25);
            \centerarc[](0.75,1.3)(270:360:0.25);
            \pic[leftstr=green] at (0.5,1.3) {invbraid}; \filldraw[gray] (0.75,1.55) circle[radius=0.1];
            \draw[green] (0.5,1.2) circle[radius=0.1];
            \draw (0.5,1.8) -- (0.5,2.3);
            \draw[blue] (0,0.8) -- (0,2.3);
            \pic[leftstr=blue] at (0,2.3) {invbraid};
            \draw[-Latex] (0,2.8) -- (0,4);
            \centerarc[green](0.5,2.8)(0:90:0.5);
            \draw[blue,-Latex] (0.5,2.8) -- (0.5,4);
            \draw[green] (1,1.8) -- (1,2.8);
            \node at (0,4.2) {$V$};
            \node[blue] at (0.5,4.2) {$M$};
        \end{tikzpicture} \;.
    \end{equation}
    Now we show that it satisfies the condition \eqref{eq:BrMpicturesA}:
    \begin{equation}
        \begin{tikzpicture}[baseline={(0,0.6)}]
             \pic at (0,0) {modbrots}; \draw (-0.5,-0.5) -- (-0.5,0); \node at (-0.5,-0.7) {$U$}; \draw (0,-0.5) -- (0,0); \node at (0,-0.7) {$V$}; \draw[blue] (0.5,-0.5) -- (0.5,0); \node[blue] at (0.5,-0.7) {$M$}; \draw[-Latex] (0,1.5) -- (0,2); \node at (-0.5,2.2) {$U$}; \draw[-Latex] (-0.5,1.5) -- (-0.5,2); \node at (0,2.2) {$V$}; \draw[-Latex, blue] (0.5,1.5) -- (0.5,2); \node[blue] at (0.5,2.2) {$M$};
        \end{tikzpicture}
        \overset{\eqref{eq:MBrfromHBr}}{=}
        \begin{tikzpicture}[baseline={(0,0.8)}]
            \pic[rightstrB = blue] at (0,0) {braidots}; 
            \pic[rightstr = green] at (1,1) {braid}; \filldraw[gray] (1.25,1.25) circle[radius=0.1]; \draw[green] (1.5,0.8) circle[radius=0.1]; \draw[green] (1.5,0.9) -- (1.5,1); \centerarc[green](0.5,1.5)(0:90:0.5); \draw[blue] (0.5,1) -- (0.5,2); \centerarc[](1.25,1.5)(0:90:0.25); \centerarc[](1.25,2)(180:270:0.25);
            \pic[leftstrB = blue] at (0,2) {braidtso};
            \node at (0,-0.5) {$U$}; \draw (0,-0.3) -- (0,0); \draw (0,1) -- (0,2); \draw[-Latex] (0,3) -- (0,3.5); \node at (0,3.7) {$U$};
            \node at (0.5,-0.5) {$V$}; \draw (0.5,-0.3) -- (0.5,0); \draw[-Latex] (0.5,3) -- (0.5,3.5); \node at (0.5,3.7) {$V$};
            \node[blue] at (1,-0.5) {$M$}; \draw[blue] (1,-0.3) -- (1,0); \draw[blue,-Latex] (1,3) -- (1,3.5); \node[blue] at (1,3.7) {$M$};
        \end{tikzpicture}
        =
        \begin{tikzpicture}[baseline={(0,0.8)}]
            \pic at (0,0) {braid}; \pic at (0.5,0.5) {HBAmodbr}; \pic at (0,2.5) {braid}; \node[blue] at (1,-0.5) {$M$}; \draw[blue] (1,-0.3) -- (1,0.5); \draw[blue,-Latex] (1,2.5) -- (1,3.5); \node[blue] at (1,3.7) {$M$};
            \node at (0.5,-0.5) {$V$}; \draw (0.5,-0.3) -- (0.5,0); \node at (0.5,3.7) {$V$};
            \node at (0,-0.5) {$U$}; \draw (0,-0.3) -- (0,0); \draw (0,0.5) -- (0,2.5); \draw[-Latex] (0.5,3) -- (0.5,3.5); \draw[-Latex] (0,3) -- (0,3.5); \node at (0,3.7) {$U$};
        \end{tikzpicture}
        \overset{\eqref{eq:MBrfromHBr}}{=}
        \begin{tikzpicture}[baseline={(0,0.6)}]
            \pic at (0,0) {modbr}; \pic at (-0.5,1.5) {braid}; \pic at (-0.5,-0.5) {braid}; \draw (-0.5,0) -- (-0.5,1.5); \draw[-Latex] (-0.5,2) -- (-0.5,2.25); \draw[-Latex] (0,2) -- (0,2.25); \draw[blue] (0.5,-0.5) -- (0.5,0); \draw[blue,-Latex] (0.5,1.5) -- (0.5,2.25); \node at (-0.5,-0.7) {$U$}; \node at (0,-0.7) {$V$}; \node[blue] at (0.5,-0.7) {$M$}; \node at (-0.5,2.45) {$U$}; \node at (0,2.45) {$V$}; \node[blue] at (0.5,2.45) {$M$};
        \end{tikzpicture}
    \end{equation}
    Now for the condition \eqref{eq:BrMpicturesB}:
    \begin{align}
        \begin{tikzpicture}[baseline={(0,0.6)}]
             \pic at (0,0) {modbrtso}; \draw[-Latex] (-0.5,1.5) -- (-0.5,2); \draw[-Latex] (0,1.5) -- (0,2); \draw[blue] (0.5,-0.5) -- (0.5,0); \draw[blue,-Latex] (0.5,1.5) -- (0.5,2); \draw (0,-0.5) -- (0,0); \draw (-0.5,-0.5) -- (-0.5,0); \node at (-0.5,-0.7) {$U$}; \node at (0,-0.7) {$V$}; \node[blue] at (0.5,-0.7) {$M$}; \node at (-0.5,2.2) {$U$}; \node at (0,2.2) {$V$}; \node[blue] at (0.5,2.2) {$M$};
        \end{tikzpicture}
        &\overset{\eqref{eq:MBrfromHBr}}{=}
        \begin{tikzpicture}[baseline={(0,1.1)}]
            \pic[rightstr = blue] at (0,0) {braidtso}; \pic[rightstr = green] at (0.5,1) {braidtso}; \filldraw[gray] (1,1.5) circle[radius=0.3]; \centerarc[](0.75,2)(0:90:0.25); \centerarc[](0.75,2.5)(180:270:0.25); \centerarc[](1.25,2)(0:90:0.25); \centerarc[](1.25,2.5)(180:270:0.25); \draw[green] (1.5,0.9) circle[radius=0.1];\draw[blue] (0,1) -- (0,2.5); \centerarc[green](0,2)(0:90:0.5); \pic[leftstr=blue] at (0,2.5) {braidots};
            \node at (0,-0.2) {$U$}; \node at (0.5,-0.2) {$V$}; \node[blue] at (1,-0.2) {$M$};
            \draw[-Latex] (0,3.5) -- (0,4); \node at (0,4.2) {$U$}; \draw[-Latex] (0.5,3.5) -- (0.5,4); \node at (0.5,4.2) {$V$}; \draw[blue,-Latex] (1,3.5) -- (1,4); \node[blue] at (1,4.2) {$M$};
        \end{tikzpicture}
    \overset{\eqref{eq:half-hexagon}}{=}
        \begin{tikzpicture}[baseline={(0,1.1)}]
            \pic[rightstr = blue] at (0.5,0) {braid}; \pic[rightstr = green] at (1,0.5) {braid}; \filldraw[gray] (1.25,0.75) circle[radius=0.1]; \centerarc[green](0.75,1)(0:180:0.25); \draw[green] (0.5,0.9) circle[radius = 0.1]; \pic[rightstr = green] at (0.25,1.5) {braid}; \filldraw[gray] (0.5,1.75) circle[radius=0.1]; \draw[green] (0.75,1.25) -- (0.75,1.5); \centerarc[green](-0.25,2)(0:90:0.5); \draw[blue] (-0.25,1.5) -- (-0.25,2.5); \pic[leftstr = blue] at (-0.25,2.5) {braid}; \pic[leftstr = blue] at (0.25,3) {braid}; \draw (1.5,1) -- (0.75,3); \centerarc[](0.5,2.5)(180:270:0.25); \centerarc[](0.5,2)(0:90:0.25);
            \draw[blue] (0.5,0.5) -- (0.25,1); \pic[rightstr = blue] at (-0.25,1) {braid};
            \draw[green] (1.5,0.4) circle[radius=0.1];
            \node at (-0.25,-0.2) {$U$}; \draw (-0.25,0) -- (-0.25,1);
            \node at (0.5,-0.2) {$V$};
            \node[blue] at (1,-0.2) {$M$};
            \draw[blue,-Latex] (0.75,3.5) -- (0.75,4); \node[blue] at (0.75,4.2) {$M$};
            \draw[-Latex] (0.25,3.5) -- (0.25,4); \node at (0.25,4.2) {$V$};
            \draw[-Latex] (-0.25,3) -- (-0.25,4); \node at (-0.25,4.2) {$U$};
        \end{tikzpicture}
        \overset{\eqref{eq:HBA}}{=}
        \begin{tikzpicture}[baseline={(0,-0.4)}]
            \pic[rightstr = green] at (-0.5,0) {braid}; \filldraw[gray] (-0.25,0.25) circle[radius=0.1]; \pic[rightstr = green] at (0,0.5) {invbraid}; \draw[green] (-0.5,0.5) -- (-0.5,1); \draw[green] (0.5,0) -- (0.5,0.5); \centerarc[green](-0.25,1)(0:180:0.25); \draw[green] (0,-0.1) circle[radius=0.1];
            \pic[rightstr = green] at (0.5,-0.5) {braid}; \filldraw[gray] (0.75,-0.25) circle[radius=0.1]; \draw[green] (1,-0.6) circle[radius=0.1];
            \pic[rightstr=blue] at (0,-1.5) {braid}; \centerarc[blue](-0.25,-1)(0:90:0.25);  \centerarc[blue](-0.25,-0.5)(180:270:0.25);
            \pic[rightstr=blue] at (-1,-0.5) {braid};
            \centerarc[green](-0.75,1.25)(0:90:0.5); \draw[green] (-0.75,1.75) -- (-1,1.75); \draw[blue] (-1,0) -- (-1,2.25);
            \pic[leftstr=blue] at (-1,2.25) {braid};\draw[smooth] plot coordinates {(0.5,1) (0.4,1.3) (-0.4,1.95) (-0.5,2.25)};
            \pic[leftstr=blue] at (-0.5,2.75) {braid}; \draw[smooth] plot coordinates {(1,1.5) (0.9,1.8) (0.1,2.45) (0,2.75)}; \draw (1,1.5) -- (1,0);
            \draw (0.5,-1) -- (0.5,-0.5);
            \node at (-1,-1.7) {$U$}; \draw (-1,-1.5) -- (-1,-0.5);
            \node at (0,-1.7) {$V$}; \node[blue] at (0.5,-1.7) {$M$};
            \draw[-Latex] (-1,2.75) -- (-1,3.75); \node at (-1,3.95) {$U$};
            \draw[-Latex] (-0.5,3.25) -- (-0.5,3.75); \node at (-0.5,3.95) {$V$};
            \draw[blue,-Latex] (0,3.25) -- (0,3.75); \node[blue] at (0,3.95) {$M$};
        \end{tikzpicture}
        =
        \begin{tikzpicture}[baseline={(0,-0.4)}]
            \pic[rightstr = green] at (-0.5,0) {braid}; \filldraw[gray] (-0.25,0.25) circle[radius=0.1]; \pic[rightstr = green] at (0,0.5) {invbraid}; \draw[green] (-0.5,0.5) -- (-0.5,1); \centerarc[green](-1,1)(0:90:0.5); \centerarc[green](-1,1.5)(0:90:0.5); \draw[smooth, green] plot coordinates {(0,1) (-0.1,1.1) (-0.4,1.4) (-0.5,1.5)}; \draw[green] (0.5,0) -- (0.5,0.5); \draw[green] (0,-0.1) circle[radius=0.1];
            \pic[rightstr = green] at (0.5,-0.5) {braid}; \filldraw[gray] (0.75,-0.25) circle[radius=0.1]; \draw[green] (1,-0.6) circle[radius=0.1];
            \pic[rightstr=blue] at (0,-1.5) {braid}; \centerarc[blue](-0.25,-1)(0:90:0.25);  \centerarc[blue](-0.25,-0.5)(180:270:0.25);
            \pic[rightstr=blue] at (-1,-0.5) {braid};
            \draw[blue] (-1,0) -- (-1,2.25);
            \pic[leftstr=blue] at (-1,2.25) {braid};\draw[smooth] plot coordinates {(0.5,1) (0.4,1.3) (-0.4,1.95) (-0.5,2.25)};
            \pic[leftstr=blue] at (-0.5,2.75) {braid}; \draw[smooth] plot coordinates {(1,1.5) (0.9,1.8) (0.1,2.45) (0,2.75)}; \draw (1,1.5) -- (1,0);
            \draw (0.5,-1) -- (0.5,-0.5);
            \node at (-1,-1.7) {$U$}; \draw (-1,-1.5) -- (-1,-0.5);
            \node at (0,-1.7) {$V$}; \node[blue] at (0.5,-1.7) {$M$};
            \draw[-Latex] (-1,2.75) -- (-1,3.75); \node at (-1,3.95) {$U$};
            \draw[-Latex] (-0.5,3.25) -- (-0.5,3.75); \node at (-0.5,3.95) {$V$};
            \draw[blue,-Latex] (0,3.25) -- (0,3.75); \node[blue] at (0,3.95) {$M$};
        \end{tikzpicture}\nonumber\\
        &=
        \begin{tikzpicture}[baseline={(0,0.6)}]
            \pic at (-0.5,-0.5) {invbraid}; \pic at (0,0) {HBAmodbr};  \pic at (-0.5,2) {braid}; \pic at (0,2.5) {HBAmodbr}; \node at (-0.5,-0.7) {$U$}; \draw (-0.5,0) -- (-0.5,2); \draw[-Latex] (-0.5,2.5) -- (-0.5,5); \node at (-0.5,5.2) {$U$}; \draw[-Latex] (0,4.5) -- (0,5); \node at (0,5.2) {$V$}; \node at (0,-0.7) {$V$}; \draw[blue] (0.5,-0.5) -- (0.5,0); \draw[blue] (0.5,2) -- (0.5,2.5); \draw[blue,-Latex] (0.5,4.5) -- (0.5,5); \node[blue] at (0.5,-0.7) {$M$}; \node[blue] at (0.5,5.2) {$M$};
        \end{tikzpicture} 
        \overset{\eqref{eq:MBrfromHBr}}{=}
        \begin{tikzpicture}[baseline={(0,0.8)}]
            \pic at (-0.5,-0.5) {invbraid}; \pic at (0,0) {modbr};  \pic at (-0.5,1.5) {braid}; \pic at (0,2) {modbr}; \node at (-0.5,-0.7) {$U$}; \draw (-0.5,0) -- (-0.5,1.5); \draw[-Latex] (-0.5,2) -- (-0.5,4); \node at (-0.5,4.2) {$U$}; \draw[-Latex] (0,3.5) -- (0,4); \node at (0,4.2) {$V$}; \node at (0,-0.7) {$V$}; \draw[blue] (0.5,-0.5) -- (0.5,0); \draw[blue] (0.5,1.5) -- (0.5,2); \draw[blue,-Latex] (0.5,3.5) -- (0.5,4); \node[blue] at (0.5,-0.7) {$M$}; \node[blue] at (0.5,4.2) {$M$};
        \end{tikzpicture} \;.
    \end{align}    
It remains to check that each $e_{V,M}$ is a morphism in $\cat{C}_A$ (i.e.\ an intertwiner of $A$-modules):
    \begin{align}
        \begin{tikzpicture}[baseline={(0,0.5)}]
            \pic at (0,0) {modbr}; \centerarc[green](0.5,1.25)(0:90:0.5); \draw (0,-0.3) -- (0,0); \node at (0,-0.5) {$V$}; \draw[blue] (0.5,-0.3) -- (0.5,0); \node[blue] at (0.5,-0.5) {$M$}; \draw[green] (1,-0.3) -- (1,1.25); \node[green] at (1,-0.5) {$A$}; \draw[-Latex] (0,1.5) -- (0,2.25); \draw[blue,-Latex] (0.5,1.5) -- (0.5,2.25);
            \node at (0,2.45) {$V$}; \node[blue] at (0.5,2.45) {$M$};
        \end{tikzpicture}
        &\overset{\eqref{eq:MBrfromHBr}}{=}
        \begin{tikzpicture}[baseline={(0,0.8)}]
            \pic at (0,0) {HBAmodbr}; \node at (0,-0.2) {$V$}; \node[blue] at (0.5,-0.2) {$M$}; \draw[-Latex] (0,2) -- (0,2.5); \node at (0,2.7) {$V$}; \centerarc[green](0.5,1)(0:90:1); \draw[green] (1.5,0) -- (1.5,1); \node[green] at (1.5,-0.2) {$A$}; \draw[blue,-Latex] (0.5,2) -- (0.5,2.5);  \node at (0.5,2.7) {$M$};
        \end{tikzpicture}
        =
        \begin{tikzpicture}[baseline={(0,0.8)}]
            \pic[rightstr = blue] at (0,0) {braid}; \pic[rightstr = green] at (0.5,0.5) {braid}; \filldraw[gray] (0.75,0.75) circle[radius=0.1]; \draw[green] (1,0.3) circle[radius=0.1]; \draw[green] (1,0.4) -- (1,0.5); \pic[leftstr = green] at (0.5,1) {braid}; \centerarc[green](1.25,1.5)(0:180:0.25); \draw[blue] (0,0.5) -- (0,1.5); \pic[leftstr = blue] at (0,1.5) {braid};
            \node at (0,-0.2) {$V$}; \node[blue] at (0.5,-0.2) {$M$}; \draw[-Latex] (0,2) -- (0,2.5); \node at (0,2.7) {$V$}; \centerarc[green](1,1.75)(0:90:0.25); \draw[green](1,2) -- (0.5,2); \draw[green] (1.5,0) -- (1.5,1.5); \node[green] at (1.5,-0.2) {$A$}; \draw[blue,-Latex] (0.5,2) -- (0.5,2.5);  \node at (0.5,2.7) {$M$};
        \end{tikzpicture}
        \overset{\eqref{eq:HBA}}{=}
        \begin{tikzpicture}[baseline={(0,0.8)}]
            \pic[rightstr = blue] at (0,0) {braid}; 
            \centerarc[green](1,0.5)(0:180:0.25); \draw[green] (0.75,0.4) circle[radius=0.1]; \draw[green] (1.25,0) -- (1.25,0.5); \node[green] at (1.25,-0.2) {$A$}; \pic[rightstr = green] at (0.5,0.75) {braid}; \filldraw[gray] (0.75,1) circle[radius=0.1]; \draw (0.5,0.5) -- (0.5,0.75); \centerarc[green](0,1.25)(0:90:0.5); \draw[blue] (0,0.5) -- (0,1.75); \pic[leftstr=blue] at (0,1.75) {braid}; \centerarc[](0.75,1.75)(180:270:0.25); \centerarc[](0.75,1.25)(0:90:0.25);
            \node at (0,-0.2) {$V$}; \node[blue] at (0.5,-0.2) {$M$}; \draw[-Latex] (0,2.25) -- (0,2.5); \node at (0,2.7) {$V$}; \draw[blue,-Latex] (0.5,2.25) -- (0.5,2.5);  \node at (0.5,2.7) {$M$};
        \end{tikzpicture}
        =
        \begin{tikzpicture}[baseline={(0,0.8)}]
            \pic[rightstr = blue] at (0,0) {braid}; 
            \centerarc[green](1,0.5)(0:180:0.25); \draw[green] (1.25,0.4) circle[radius=0.1]; \draw[green] (0.75,0) -- (0.75,0.5); \node[green] at (0.75,-0.2) {$A$}; \pic[rightstr = green] at (0.5,0.75) {braid}; \filldraw[gray] (0.75,1) circle[radius=0.1]; \draw (0.5,0.5) -- (0.5,0.75); \centerarc[green](0,1.25)(0:90:0.5); \draw[blue] (0,0.5) -- (0,1.75); \pic[leftstr=blue] at (0,1.75) {braid}; \centerarc[](0.75,1.75)(180:270:0.25); \centerarc[](0.75,1.25)(0:90:0.25);
            \node at (0,-0.2) {$V$}; \node[blue] at (0.5,-0.2) {$M$}; \draw[-Latex] (0,2.25) -- (0,2.5); \node at (0,2.7) {$V$}; \draw[blue,-Latex] (0.5,2.25) -- (0.5,2.5);  \node at (0.5,2.7) {$M$};
        \end{tikzpicture} \nonumber\\
        &\overset{\eqref{eq:HBA}}{=}
        \begin{tikzpicture}[baseline={(0,0.8)}]
            \pic[rightstr=blue] at (0,0) {braid}; \pic[rightstr=green] at (0.5,0.5) {braid}; \pic[rightstr=green] at (1,1) {braid}; \draw[green] (1,0) -- (1,0.5); \node[green] at (1,-0.2) {$A$}; \filldraw[gray] (1.25,1.25) circle[radius=0.1]; \draw[green] (1.5,0.9) circle[radius=0.1]; \centerarc[green](0.75,1.5)(0:180:0.25); \draw[green] (0.5,1) -- (0.5,1.5); \centerarc[green](0.25,1.75)(0:90:0.5); \draw[green] (0.25,2.25) -- (0,2.25); \draw[blue] (0,0.5) -- (0,2.5); \pic[leftstr=blue] at (0,2.5) {braid};
            \draw[smooth] plot coordinates {(1.5,1.5) (1.4,1.7) (0.6,2.3) (0.5,2.5)};
             \node at (0,-0.2) {$V$}; \node[blue] at (0.5,-0.2) {$M$}; \draw[-Latex] (0,3) -- (0,3.25); \node at (0,3.45) {$V$}; \draw[blue,-Latex] (0.5,3) -- (0.5,3.25);  \node at (0.5,3.45) {$M$};
        \end{tikzpicture}
        =
        \begin{tikzpicture}[baseline={(0,0.8)}]
            \pic at (0,0.5) {HBAmodbr}; \draw (0,0) -- (0,0.5); \draw[blue] (0.5,0) -- (0.5,0.5); \centerarc[green](0.5,0)(0:90:0.5); \node at (0,-0.2) {$V$}; \node[blue] at (0.5,-0.2) {$M$}; \node[green] at (1,-0.2) {$A$}; \draw[-Latex] (0,2.5) -- (0,2.8); \node at (0,3) {$V$}; \draw[blue,-Latex] (0.5,2.5) -- (0.5,2.8); \node[blue] at (0.5,3) {$M$};
        \end{tikzpicture}
        \overset{\eqref{eq:MBrfromHBr}}{=} 
        \begin{tikzpicture}[baseline={(0,0.8)}]
            \pic at (0,0.75) {modbr}; \centerarc[green](0.5,0)(0:90:0.5); \draw (0,0) -- (0,0.75); \node at (0,-0.2) {$V$}; \draw[blue] (0.5,0) -- (0.5,0.75); \node[blue] at (0.5,-0.2) {$M$}; \node[green] at (1,-0.2) {$A$}; \draw[-Latex] (0,2.25) -- (0,2.75); \node at (0,2.95) {$V$}; \draw[blue,-Latex] (0.5,2.25) -- (0.5,2.75); \node[blue] at (0.5,2.95) {$M$};
        \end{tikzpicture}
    \end{align}
    
Finally, let us show that $\gamma^A$ in \eqref{eq:HbrfromMBr} is a half-braiding satisfying \eqref{eq:HBA}. 
It is invertible by virtue of the fact that the module braidings and braiding in $\cat{C}$ are invertible. The hexagon \eqref{eq:half-hexagon} is also a straightforward string diagram computation using the hexagon for the braiding and the conditions \eqref{eq:BrMpicturesA} and \eqref{eq:BrMpicturesB}. We check condition \eqref{eq:HBA} directly. The first equality in \eqref{eq:HBA} follows from
    \begin{equation*} 
        \begin{tikzpicture}[baseline={(0,0.8)}]
            \centerarc[green](0,0.3)(0:180:0.5); \draw[green] (-0.5,0) -- (-0.5,0.3); \draw[green] (0.5,0) -- (0.5,0.3);
            \draw[green] (0,0.8) -- (0,1); \pic[rightstr = green] at (-0.5,1) {braid};             \filldraw[gray] (-0.25,1.25) circle[radius=0.1]; \draw[green, -Latex] (-0.5,1.5) -- (-0.5,1.8); \draw[-Latex] (0,1.5) -- (0,1.8);
            \centerarc[](-0.75,1)(270:360:0.25); \centerarc[](-0.75,0.5)(90:180:0.25); \draw (-1,0) -- (-1,0.5);
            \node at (-1,-0.2) {$V$}; \node[green] at (-0.5,-0.2) {$A$}; \node[green] at (0.5,-0.2) {$A$};
            \node at (0,2) {$V$}; \node[green] at (-0.5,2) {$A$};
        \end{tikzpicture}
        \overset{\eqref{eq:HbrfromMBr}}{=}
        \begin{tikzpicture}[baseline={(0,0.8)}]
            \centerarc[green](0,0.5)(0:180:0.5); \draw[green] (0,1) -- (0,1.25); \pic[module = green] at (-0.5,1.25) {modbr}; \pic[rightstr = green] at (-0.5,2.75) {invbraid}; \draw[green,-Latex] (-0.5,3.25) -- (-0.5,3.75); \node[green] at (-0.5,3.95) {$A$}; \draw[-Latex] (0,3.25) -- (0,3.75); \node at (0,3.95) {$V$}; \draw[smooth] plot coordinates {(-1,0) (-1,0.5) (-0.9,0.8) (-0.6,1.05) (-0.5,1.25)}; \node at (-1,-0.2) {$V$}; \draw[green] (-0.5,0) -- (-0.5,0.5); \draw[green] (0.5,0) -- (0.5,0.5); \node[green] at (-0.5,-0.2) {$A$}; \node[green] at (0.5,-0.2) {$A$};
        \end{tikzpicture}
        \overset{\substack{e_{V,A} \\ \text{intertwiner}}}{=}
        \begin{tikzpicture}[baseline={(0,0.8)}]
           \pic[module = green] at (0,0.25) {modbr}; \pic[rightstr = green] at (0,1.75) {invbraid}; \pic[rightstr = green] at (0.5,2.25) {invbraid}; \centerarc[green](0.25,2.75)(0:180:0.25); \draw[green,-Latex](0.25,3) -- (0.25,3.5); \node[green] at (0.25,3.7) {$A$}; \draw[green] (0,2.25) -- (0,2.75);
           \draw (0,0) -- (0,0.25); \node at (0,-0.2) {$V$}; \draw[green] (0.5,0) -- (0.5,0.25); \node[green] at (0.5,-0.2) {$A$}; \draw[green] (1,0) -- (1,2.25); \node[green] at (1,-0.2) {$A$}; \draw[-Latex] (1,2.75) -- (1,3.5); \node at (1,3.7) {$V$};
        \end{tikzpicture}
        \overset{\eqref{eq:HbrfromMBr}}{=}
        \begin{tikzpicture}[baseline={(0,0.8)}]
            \pic[rightstr = green] at (-0.5,0) {braid}; \filldraw[gray] (-0.25,0.25) circle[radius=0.1]; \pic[rightstr = green] at (0,0.5) {invbraid};
            \node at (-0.5,-0.2) {$V$}; \node[green] at (0,-0.2) {$A$}; \node[green] at (0.5,-0.2) {$A$};
            \draw[-Latex] (0.5,1) -- (0.5,1.8); \node at (0.5,2) {$V$}; \draw[green] (-0.5,0.5) -- (-0.5,1); \draw[green] (0.5,0) -- (0.5,0.5); \centerarc[green](-0.25,1)(0:180:0.25); \draw[green, -Latex] (-0.25,1.25) -- (-0.25,1.8); \node[green] at (-0.25,2) {$A$};
        \end{tikzpicture} \;.
    \end{equation*}
    The second equality in \eqref{eq:HBA} follows from
    \begin{equation*}
        \begin{tikzpicture}[baseline={(0,0.8)}]
            \centerarc[green](0,0.3)(0:180:0.5); \draw[green] (-0.5,0) -- (-0.5,0.3); \draw[green] (0.5,0) -- (0.5,0.3);
            \draw[green] (0,0.8) -- (0,1); \pic[rightstr = green] at (-0.5,1) {braid};             \filldraw[gray] (-0.25,1.25) circle[radius=0.1]; \draw[green, -Latex] (-0.5,1.5) -- (-0.5,1.8); \draw[-Latex] (0,1.5) -- (0,1.8);
            \centerarc[](-0.75,1)(270:360:0.25); \centerarc[](-0.75,0.5)(90:180:0.25); \draw (-1,0) -- (-1,0.5);
            \node at (-1,-0.2) {$V$}; \node[green] at (-0.5,-0.2) {$A$}; \node[green] at (0.5,-0.2) {$A$};
            \node at (0,2) {$V$}; \node[green] at (-0.5,2) {$A$};
        \end{tikzpicture}
        \overset{\eqref{eq:HbrfromMBr}}{=}
        \begin{tikzpicture}[baseline={(0,0.8)}]
            \centerarc[green](0,0.5)(0:180:0.5); \draw[green] (0,1) -- (0,1.25); \pic[module = green] at (-0.5,1.25) {modbr}; \pic[rightstr = green] at (-0.5,2.75) {invbraid}; \draw[green,-Latex] (-0.5,3.25) -- (-0.5,3.75); \node[green] at (-0.5,3.95) {$A$}; \draw[-Latex] (0,3.25) -- (0,3.75); \node at (0,3.95) {$V$}; \draw[smooth] plot coordinates {(-1,0) (-1,0.5) (-0.9,0.8) (-0.6,1.05) (-0.5,1.25)}; \node at (-1,-0.2) {$V$}; \draw[green] (-0.5,0) -- (-0.5,0.5); \draw[green] (0.5,0) -- (0.5,0.5); \node[green] at (-0.5,-0.2) {$A$}; \node[green] at (0.5,-0.2) {$A$};
        \end{tikzpicture}
        \overset{\substack{e_{V,\blank} 
        \\ 
        \text{natural}}}{=} 
        \begin{tikzpicture}[baseline={(0,0.8)}]
            \pic[moduleA = green, moduleB = green] at (0,0.25) {modbrots}; \pic[rightstrA = green, rightstrB = green] at (-0.5,1.75) {braiduts}; \centerarc[green](-0.25,2.75)(0:180:0.25); \draw[green, -Latex] (-0.25,3) -- (-0.25,3.75);\node[green] at (-0.25,3.95) {$A$}; \draw[-Latex] (0.5,2.75) -- (0.5,3.75);\node at (0.5,3.95) {$V$};
            \draw (-0.5,0) -- (-0.5,0.25); \node at (-0.5,-0.2) {$V$}; \draw[green] (0,0) -- (0,0.25); \draw[green] (0.5,0) -- (0.5,0.25); \node[green] at (0,-0.2) {$A$}; \node[green] at (0.5,-0.2) {$A$};
        \end{tikzpicture}
        \overset{\eqref{eq:BrM1}}{=}
        \begin{tikzpicture}[baseline={(0,0.8)}]
            \pic[rightstr = green] (0,0) {braid}; \pic[module = green] at (0.5,0.5) {modbr}; \pic[rightstr = green] at (0.5,2) {invbraid}; \draw[green] (0,0.5) -- (0,2.5); \centerarc[green](0.25,2.5)(0:180:0.25); \draw[green,-Latex] (0.25,2.75) -- (0.25,3.5); \node[green] at (0.25,3.7) {$A$}; \draw[-Latex] (1,2.5) -- (1,3.5); \node at (1,3.7) {$V$};
            \node at (0,-0.2) {$V$}; \node[green] at (0.5,-0.2) {$A$}; \node[green] at (1,-0.2) {$A$}; \draw[green] (1,0) -- (1,0.5);
        \end{tikzpicture} 
        \overset{\eqref{eq:HbrfromMBr}}{=}
        \begin{tikzpicture}[baseline={(0,0.8)}]
            \pic[rightstr = green] at (-0.5,0) {braid}; \pic[rightstr = green] at (0,0.5) {braid}; \filldraw[gray] (0.25,0.75) circle[radius=0.1];
            \node at (-0.5,-0.2) {$V$}; \node[green] at (0,-0.2) {$A$}; \node[green] at (0.5,-0.2) {$A$};
            \draw[-Latex] (0.5,1) -- (0.5,1.8); \node at (0.5,2) {$V$}; \draw[green] (-0.5,0.5) -- (-0.5,1); \draw[green] (0.5,0) -- (0.5,0.5); \centerarc[green](-0.25,1)(0:180:0.25); \draw[green, -Latex] (-0.25,1.25) -- (-0.25,1.8); \node[green] at (-0.25,2) {$A$};
        \end{tikzpicture} \;.
    \end{equation*}
\end{proof}

\begin{definition}[{\cite[Def.\,2.31]{JR21}}]\label{def:half-braided-algebra}
    An algebra $(A,m,u)$ in $\cat{C}$ equipped with a half-braiding $\gamma^A$ satisfying \eqref{eq:HBA} is called a \emph{half-braided algebra}, and $\gamma^A$ is called a \textit{half-braided structure} on $A$.
\end{definition}

\begin{remark}
Setting $A = \1$, the trivial algebra, we observe that a module braiding for the regular $\cat{C}$-module category is the same as a monoidal natural automorphism of the identity functor, recovering \cite[Prop.\,2.2]{HS26}.
\end{remark}

Assume $\cat{C}$ is ribbon (and hence in particular rigid).
For our next braided module category construction, we require the coend of the functor $(\blank)^\ast\otimes\bblank:\cat{C}^{\mathrm{op}}\times\cat{C}\to\cat{C}$, denoted
\begin{equation}
    \mathcal{L} = \int^{V\in\cat{C}} V^\ast\otimes V \;.
\end{equation}
See \cite[Ch.\,5]{KL01} for a detailed review of coends in braided monoidal categories. The coend need not exist in $\cat{C}$, and often does not (for example, the coend over the category of finite-dimensional modules over an infinite-dimensional ribbon Hopf algebra is an infinite-dimensional algebra). In the construction below we need to assume $\mathcal{L}$ exists as an object in $\cat{C}$, so we do so in what follows.
By definition, $\mathcal{L}$ comes equipped with a dinatural transformation $i = (i_V: V^\ast\otimes V \to \mathcal{L})_{V\in\cat{C}}$ satisfying the following universal property: Given $A\in\cat{C}$ and a dinatural transformation 
$b:(\blank)^\ast\otimes\bblank\ddot{\to}A$, there exists a unique morphism $\beta:\mathcal{L}\to A$ 
in $\cat{C}$ such that the diagram
\begin{equation}\label{eq:Coend-UP}
    \begin{tikzcd}
    	& {V^\ast\otimes V} & \\
    	\mathcal{L} && A
    	\arrow["{i_V}"', from=1-2, to=2-1]
    	\arrow["{b_V}", from=1-2, to=2-3]
    	\arrow["{\beta}"', dashed, from=2-1, to=2-3]
    \end{tikzcd} 
\end{equation}
commutes for all $V$ in $\cat{C}$. We depict the dinatural transformation $i$ graphically as
\begin{equation}
    i_V = \begin{tikzpicture}[baseline={(0,0.8)}]
        \node at (0,-0.2) {$V^\ast$}; \node at (0.5,-0.2) {$V$}; \draw[-Latex] (0,1) -- (0,0.3); \draw (0,1) -- (0,0); \draw[-Latex] (0.5,0) -- (0.5,0.8); \draw (0.5,1) -- (0.5,0); \draw[gray, ultra thick] (-0.2,1) -- (0.7,1); \draw[-Latex] (0.25,1) -- (0.25,2); \node at (0.25,2.2) {$\mathcal{L}$}; \node[gray] at (1,1) {$i_V$};
    \end{tikzpicture} \;.
\end{equation}
Because $\cat{C}$ is rigid, $\blank\otimes V$ and $V\otimes\blank$ admit right adjoints and hence preserve all colimits. In particular, they preserve $\mathcal{L}$ for all $V\in\cat{C}$. Thus the object $\mathcal{L}\otimes\mathcal{L}$ is exhibited as a coend over $\cat{C}\times\cat{C}$ with dinatural transformation $(i_U\otimes i_V)_{U,V\in\cat{C}}$. The morphisms $m:\mathcal{L}\otimes\mathcal{L}\to\mathcal{L}$ and $u:\1\to\mathcal{L}$ defined via
\begin{subequations}
    \begin{equation}
            m\circ(i_U\otimes i_V) = \begin{tikzpicture}[baseline={(0,0.8)}]
            \node at (0,-0.2) {$U^\ast$}; \node at (0.5,-0.2) {$U$}; \draw[-Latex] (0,1) -- (0,0.3); \draw (0,1) -- (0,0); \draw[-Latex] (0.5,0) -- (0.5,0.8); \draw (0.5,1) -- (0.5,0); \draw[gray] (-0.2,1) -- (0.7,1); \node[gray] at (1,1) {$i_U$};
            \node at (1.5,-0.2) {$V^\ast$}; \node at (2,-0.2) {$V$}; \draw[-Latex] (1.5,1) -- (1.5,0.3); \draw (1.5,1) -- (1.5,0); \draw[-Latex] (2,0) -- (2,0.8); \draw (2,1) -- (2,0); \draw[gray] (1.3,1) -- (2.2,1); \node[gray] at (2.5,1) {$i_V$};
            \centerarc[](1,1)(0:180:0.75); \draw[-Latex] (1,1.75) -- (1,2.5); \node at (1,2.7) {$\mathcal{L}$};
        \end{tikzpicture}
        :=
        \begin{tikzpicture}[baseline={(0,0.8)}]
            \node at (0,-0.2) {$U^\ast$}; \node at (0.5,-0.2) {$U$}; \node at (1,-0.2) {$V^\ast$}; \node at (1.5,-0.2) {$V$}; \pic at (0,0) {braidtso}; \draw[-Latex] (1.5,0) -- (1.5,1.4); \draw (1.5,0) -- (1.5,1.5); \draw[gray] (-0.2,1.5) -- (1.7,1.5); \draw[-Latex] (1,1) -- (1,1.4); \draw[-Latex] (0.5,1.2) -- (0.5,1); \draw[-Latex] (0,1.2) -- (0,1); \draw (0,1) -- (0,1.5); \draw (0.5,1) -- (0.5,1.5); \draw (1,1) -- (1,1.5); \draw[-Latex] (0.75,1.5) -- (0.75,2.5); \node at (0.75,2.7) {$\mathcal{L}$}; \node[gray] at (2.2,1.5) {$i_{U\otimes V}$};
        \end{tikzpicture} \;,\label{eq:LalgA} 
    \end{equation}
    \begin{equation}
        u = \begin{tikzpicture}[baseline={(0,0)}]
            \draw (0,0) circle[radius=0.1]; \draw[-Latex] (0,0.1) -- (0,1); \node at (0,1.2) {$\mathcal{L}$};
        \end{tikzpicture} 
        :=
        \begin{tikzpicture}[baseline={(0,0)}]
            \draw[gray] (0,0) -- (0.5,0); \node[gray] at (0.8,0) {$i_\1$}; \draw[-Latex] (0.25,0) -- (0.25,1); \node at (0.25,1.2) {$\mathcal{L}$};
        \end{tikzpicture} 
        \label{eq:LalgB}
    \end{equation}
\end{subequations}
make $(\mathcal{L},m,u)$ an algebra in $\cat{C}$.

\begin{definition}\label{def:QMM}
    Assume $\cat{C}$ is ribbon and contains the coend $\mathcal{L}$. Let $(A,m,u)$ be an algebra in $\cat{C}$. A \emph{quantum moment map} is an algebra homomorphism $\mu:\mathcal{L}\to A$ in $\cat{C}$ such that 
    \begin{equation}
        \begin{tikzcd}
            (V^\ast\otimes V)\otimes A && A\otimes A && A \\ V^\ast \otimes (A\otimes V) && (A\otimes V^\ast)\otimes V && A\otimes A
            \arrow["(\mu\circ i_V)\otimes A", from=1-1, to=1-3]\arrow["m", from=1-3, to=1-5]
            \arrow["(V^\ast\otimes\sigma_{V,A})\circ\alpha"', from=1-1, to=2-1]\arrow["(\sigma_{A,V^\ast}^{-1}\otimes V)\circ\alpha^{-1}"', from=2-1, to=2-3]\arrow["(A\otimes (\mu\circ i_V))\circ\alpha"', from=2-3, to=2-5]\arrow["m"', from=2-5, to=1-5]
        \end{tikzcd}
    \end{equation}
    for all $V\in\cat{C}$. This condition is depicted
    \begin{equation}\label{eq:QMM}
        \begin{tikzpicture}[baseline={(0,0.8)}]
            \node at (0,-0.2) {$V^\ast$}; \node at (0.5,-0.2) {$V$}; \draw[-Latex] (0,1) -- (0,0.3); \draw (0,1) -- (0,0); \draw[-Latex] (0.5,0) -- (0.5,0.8); \draw (0.5,1) -- (0.5,0); \draw[gray, ultra thick] (-0.2,1) -- (0.7,1); \draw (0.25,1) -- (0.25,1.5);
            \centerarc[green](0.75,1.5)(0:180:0.5); \filldraw[red] (0.25,1.5) circle[radius=0.1]; \node[red] at (0.55,1.5) {$\mu$}; \draw[green] (1.25,0) -- (1.25,1.5); \node[green] at (1.25,-0.2) {$A$}; \draw[green,-Latex] (0.75,2) -- (0.75,2.75); \node[green] at (0.75,2.95) {$A$};
        \end{tikzpicture}  
        =
        \begin{tikzpicture}[baseline={(0,0.8)}]
            \node at (0,-0.2) {$V^\ast$}; \node at (0.5,-0.2) {$V$}; \node[green] at (1,-0.2) {$A$}; \pic[rightstr = green] at (0.5,0) {braid}; \pic[rightstr = green] at (0,0.5) {invbraid}; \draw (0,0) -- (0,0.5); \draw[-Latex] (0,0.5) -- (0,0.2); \draw[-Latex] (1,0.5) -- (1,0.8); \draw (1,0.5) -- (1,1); \draw[gray] (0.3,1) -- (1.2,1); \draw (0.75,1) -- (0.75,1.5); \centerarc[green](0.375,1.5)(0:180:0.375); \filldraw[red] (0.75,1.5) circle[radius=0.1]; \node[red] at (1.05,1.5) {$\mu$}; \draw[green] (0,1) -- (0,1.5); \draw[green,-Latex] (0.375,1.875) -- (0.375,2.625); \node[green] at (0.375,2.825) {$A$};
        \end{tikzpicture}
    \end{equation}
\end{definition}
\begin{proposition}[{\cite[Sec.\,2.5]{JR21}}]\label{prop:QMMiffHalfBr}
Assume $\cat{C}$ is ribbon and contains the coend $\mathcal{L}$. Given an algebra $A\in\cat{C}$, there is a bijective correspondence between half-braided structures for $A$ and quantum moment maps for $A$. 
\end{proposition}
\begin{proof}
First we show that a half-braided structure can be used to construct a quantum moment map. From a half-braided structure for an algebra $A\in\cat{C}$ we obtain a dinatural transformation 
$b:(\blank)^\ast\otimes\bblank\,\ddot{\to}\,A$ 
as follows:
    \begin{equation}\label{eq:DinatFromHB}
        b_{V} := 
        \begin{tikzpicture}[baseline={(0,0.8)}]
            \node at (0,-0.2) {$V^\ast$}; \node at (0.5,-0.2) {$V$}; \pic[rightstr = green] at (0.5,0.25) {braid}; \pic[leftstr = green] at (0.5,0.75) {braid}; \centerarc[](0.25,1.25)(0:180:0.25); \draw[-Latex] (0,1.25) -- (0,0); \draw (0.5,0) -- (0.5,0.25); \draw[green] (1,0.15) circle[radius=0.1]; \draw[green,-Latex] (1,1.25) -- (1,1.8); \node[green] at (1,2) {$A$}; \filldraw[gray] (0.75,0.5) circle[radius=0.1];
        \end{tikzpicture}
    \end{equation}
    The universal property of the coend \eqref{eq:Coend-UP} implies that there exists a unique morphism $\mu:
    \mathcal{L}\to A$ such that $b_V = \mu\circ i_V$. 
    We now show that $\mu$ obtained this way is a quantum moment map. First, we show that it is an algebra homomorphism:
    \begin{align}
        \begin{tikzpicture}[baseline={(0,0.8)}]
            \node at (0,-0.2) {$U^\ast$}; \node at (0.5,-0.2) {$U$}; \draw[-Latex] (0,1) -- (0,0.3); \draw (0,1) -- (0,0); \draw[-Latex] (0.5,0) -- (0.5,0.8); \draw (0.5,1) -- (0.5,0); \draw[gray] (-0.2,1) -- (0.7,1); \node[gray] at (1,1) {$i_U$};
            \node at (1.5,-0.2) {$V^\ast$}; \node at (2,-0.2) {$V$}; \draw[-Latex] (1.5,1) -- (1.5,0.3); \draw (1.5,1) -- (1.5,0); \draw[-Latex] (2,0) -- (2,0.8); \draw (2,1) -- (2,0); \draw[gray] (1.3,1) -- (2.2,1); \node[gray] at (2.5,1) {$i_V$};
            \draw (0.25,1) -- (0.25,1.5); \draw (1.75,1) -- (1.75,1.5); \centerarc[green](1,1.5)(0:180:0.75); \draw[green,-Latex] (1,2.25) -- (1,3); \node[green] at (1,3.2) {$A$}; \filldraw[red] (0.25,1.5) circle[radius=0.1]; \filldraw[red] (1.75,1.5) circle[radius=0.1];
        \end{tikzpicture}
        &\overset{\eqref{eq:DinatFromHB}}{=}
        \begin{tikzpicture}[baseline={(0,0.8)}]
            \node at (0,-0.2) {$U^\ast$}; \node at (0.5,-0.2) {$U$}; \pic[rightstr = green] at (0.5,0.25) {braid}; \pic[leftstr = green] at (0.5,0.75) {braid}; \centerarc[](0.25,1.25)(0:180:0.25); \draw[-Latex] (0,1.25) -- (0,0); \draw (0.5,0) -- (0.5,0.25); \draw[green] (1,0.15) circle[radius=0.1]; \filldraw[gray] (0.75,0.5) circle[radius=0.1];
            \centerarc[green](1.75,1.25)(0:180:0.75); \draw[green,-Latex] (1.75,2) -- (1.75,2.75); \node[green] at (1.75,2.95) {$A$}; 
            \node at (1.5,-0.2) {$V^\ast$}; \node at (2,-0.2) {$V$}; \pic[rightstr = green] at (2,0.25) {braid}; \pic[leftstr = green] at (2,0.75) {braid}; \centerarc[](1.75,1.25)(0:180:0.25); \draw[-Latex] (1.5,1.25) -- (1.5,0); \draw (2,0) -- (2,0.25); \draw[green] (2.5,0.15) circle[radius=0.1]; \filldraw[gray] (2.25,0.5) circle[radius=0.1];
        \end{tikzpicture}
        =
        \begin{tikzpicture}[baseline={(0,0.8)}]
            \pic at (0,0) {braidtso}; \pic[rightstr = green] at (2,0.5) {braid}; \filldraw[gray] (2.25,0.75) circle[radius=0.1]; \draw[green] (2.5,0.4) circle[radius=0.1]; \pic[leftstr = green] at (2,1) {braid};
            \pic[rightstr = green] at (1,1) {braid}; \filldraw[gray] (1.25,1.25) circle[radius=0.1]; \draw[green] (1.5,0.9) circle[radius=0.1]; \pic[leftstr = green] at (1,1.5) {braid};
            \pic[leftstr = green] at (1.5,2) {braid}; \draw (2,1.5) -- (2,2); \draw[smooth] plot coordinates {(1.5,0) (1.6,0.2) (1.9,0.3) (2,0.5) };
            \centerarc[green](2.25,2.5)(0:180:0.25); \draw[green] (2.5,1.5) -- (2.5,2.5); \draw[green,-Latex] (2.25,2.75) -- (2.25,3.5); \node[green] at (2.25,3.7) {$A$};
            \centerarc[](0.75,2)(0:180:0.25); \draw (0.5,1) -- (0.5,2);
            \centerarc[](0.75,2.5)(0:180:0.75); \draw (0,1) -- (0,2.5);
            \draw[-Latex] (0.5,1.8) -- (0.5,1.3); \draw[-Latex] (0,1.8) -- (0,1.3);
            \node at (0,-0.2) {$U^\ast$}; \node at (0.5,-0.2) {$U$}; \node at (1,-0.2) {$V^\ast$}; \node at (1.5,-0.2) {$V$};
        \end{tikzpicture}
        \overset{\eqref{eq:HBA}}{=}
        \begin{tikzpicture}[baseline={(0,0.8)}]
            \pic at (0,0) {braidtso}; \pic[rightstr = green] at (1.5,0.5) {braid}; \filldraw[gray] (1.75,0.75) circle[radius=0.1]; \draw[green] (2,0.4) circle[radius=0.1]; 
            \pic[leftstr = green] at (1.5,2) {braid}; \draw[green,-Latex] (2,2.5) -- (2,3.5); \node[green] at (2,3.7) {$A$};
            \centerarc[](0.75,2)(0:180:0.25); \draw (0.5,1) -- (0.5,2);
            \centerarc[](0.75,2.5)(0:180:0.75); \draw (0,1) -- (0,2.5);
            \pic[rightstr = green] at (1,1) {braid}; \filldraw[gray] (1.25,1.25) circle[radius=0.1]; \pic[leftstr = green] at (1,1.5) {braid}; \draw (2,1) -- (2,2);
            \centerarc[green](1.5,0.75)(90:180:0.25); \draw[green] (1.25,0.65) circle[radius=0.1];
            \draw[-Latex] (0.5,1.8) -- (0.5,1.3); \draw[-Latex] (0,1.8) -- (0,1.3);
            \node at (0,-0.2) {$U^\ast$}; \node at (0.5,-0.2) {$U$}; \node at (1,-0.2) {$V^\ast$}; \draw (1.5,0) -- (1.5,0.5); \node at (1.5,-0.2) {$V$};
        \end{tikzpicture} 
        \nonumber\\
    &\overset{\eqref{eq:half-hexagon}}{=}
        \begin{tikzpicture}[baseline={(0,0.8)}]
            \pic at (0,0) {braidtso}; 
            \pic[leftstr = green] at (1,2) {braidots}; \draw[green,-Latex] (2,3) -- (2,4); \node[green] at (2,4.2) {$A$}; \draw[green] (2,0.9) circle[radius=0.1];
            \centerarc[](0.75,3)(0:180:0.25); \draw (0.5,1) -- (0.5,3);
            \centerarc[](0.75,3)(0:180:0.75); \draw (0,1) -- (0,3);
            \pic[rightstr = green] at (1,1) {braidtso}; \filldraw[gray] (1.5,1.5) circle[radius=0.3];  \draw[-Latex] (0.5,1.8) -- (0.5,1.3); \draw[-Latex] (0,1.8) -- (0,1.3);
            \node at (0,-0.2) {$U^\ast$}; \node at (0.5,-0.2) {$U$}; \node at (1,-0.2) {$V^\ast$};\draw (1.5,0) -- (1.5,1); \node at (1.5,-0.2) {$V$};
        \end{tikzpicture} \overset{\eqref{eq:DinatFromHB}}{=} 
        \begin{tikzpicture}[baseline={(0,0.8)}]
        \node at (0,-0.2) {$U^\ast$}; \node at (0.5,-0.2) {$U$}; \node at (1,-0.2) {$V^\ast$}; \node at (1.5,-0.2) {$V$}; \pic at (0,0) {braidtso}; \draw[-Latex] (1.5,0) -- (1.5,1.4); \draw (1.5,0) -- (1.5,1.5); \draw[gray] (-0.2,1.5) -- (1.7,1.5); \draw[-Latex] (1,1) -- (1,1.4); \draw[-Latex] (0.5,1.2) -- (0.5,1); \draw[-Latex] (0,1.2) -- (0,1); \draw (0,1) -- (0,1.5); \draw (0.5,1) -- (0.5,1.5); \draw (1,1) -- (1,1.5); \draw (0.75,1.5) -- (0.75,2); \draw[green,-Latex] (0.75,2) -- (0.75,3); \node[green] at (0.75,3.2) {$A$}; \filldraw[red] (0.75,2) circle[radius = 0.1];  \node[gray] at (2.2,1.5) {$i_{U\otimes V}$};
    \end{tikzpicture}
   \end{align}
    Preservation of the unit is clear from $\gamma^A_\1 = \Id_A$, which follows from the hexagon for $\gamma^A$. Now we show that the quantum moment map condition \eqref{eq:QMM} is satisfied:
    \begin{align}
        \begin{tikzpicture}[baseline={(0,0.8)}]
            \node at (0,-0.2) {$V^\ast$}; \node at (0.5,-0.2) {$V$}; \draw[-Latex] (0,1) -- (0,0.3); \draw (0,1) -- (0,0); \draw[-Latex] (0.5,0) -- (0.5,0.8); \draw (0.5,1) -- (0.5,0); \draw[gray, ultra thick] (-0.2,1) -- (0.7,1); \draw (0.25,1) -- (0.25,1.5);
            \centerarc[green](0.75,1.5)(0:180:0.5); \filldraw[red] (0.25,1.5) circle[radius=0.1]; \node[red] at (0.55,1.5) {$\mu$}; \draw[green] (1.25,0) -- (1.25,1.5); \node[green] at (1.25,-0.2) {$A$}; \draw[green,-Latex] (0.75,2) -- (0.75,2.75); \node[green] at (0.75,2.95) {$A$};
        \end{tikzpicture}
        &=
        \begin{tikzpicture}[baseline={(0,0.8)}]
            \node at (0,-0.2) {$V^\ast$}; \node at (0.5,-0.2) {$V$}; \pic[rightstr = green] at (0.5,0.25) {braid}; \pic[leftstr = green] at (0.5,0.75) {braid}; \centerarc[](0.25,1.25)(0:180:0.25); \draw[-Latex] (0,1.25) -- (0,0); \draw (0.5,0) -- (0.5,0.25); \draw[green] (1,0.15) circle[radius=0.1]; \filldraw[gray] (0.75,0.5) circle[radius=0.1];
            \draw[green] (1,1.25) -- (1,1.5);
            \centerarc[green](1.25,1.5)(0:180:0.25); \draw[green] (1.5,0) -- (1.5,1.5); \node[green] at (1.5,-0.2) {$A$}; \draw[green,-Latex] (1.25,1.75) -- (1.25,2.5); \node[green] at (1.25,2.7) {$A$};
        \end{tikzpicture}
        \overset{\eqref{eq:HBA}}{=}
        \begin{tikzpicture}[baseline={(0,0.8)}]
            \node at (0,-0.2) {$V^\ast$}; \node at (0.5,-0.2) {$V$}; \pic[rightstr = green] at (0.5,0.75) {braid}; \filldraw[gray] (0.75,1) circle[radius=0.1]; \pic[leftstr = green] at (0.5,1.25) {braid}; \centerarc[](0.25,1.75)(0:180:0.25); \draw[-Latex] (0,1.75) -- (0,0); \draw (0.5,0) -- (0.5,0.75);
            \centerarc[green](1,0.5)(0:180:0.25); \draw[green] (0.75,0.4) circle[radius=0.1]; \draw[green] (1.25,0) -- (1.25,0.5); \node[green] at (1.25,-0.2) {$A$};
            \draw[green,-Latex] (1,1.75) -- (1,2.5); \node[green] at (1,2.7) {$A$};            
        \end{tikzpicture}
        =
        \begin{tikzpicture}[baseline={(0,0.8)}]
            \node at (0,-0.2) {$V^\ast$}; \node at (0.5,-0.2) {$V$}; \pic[rightstr = green] at (0.5,1.25) {braid}; \filldraw[gray] (0.75,1.5) circle[radius=0.1]; \pic[leftstr = green] at (0.5,1.75) {braid}; \centerarc[](0.25,2.25)(0:180:0.25); \draw[-Latex] (0,2.25) -- (0,0); \draw (0.5,0) -- (0.5,1.25);
            \draw[green] (1.25,0.75) -- (1,1.25);
            \centerarc[green](1.25,0.5)(0:180:0.25); \draw[green] (1.5,0.4) circle[radius=0.1]; \draw[green] (1,0) -- (1,0.5); \node[green] at (1,-0.2) {$A$};
            \draw[green,-Latex] (1,2.25) -- (1,3); \node[green] at (1,3.2) {$A$};            
        \end{tikzpicture}
        \overset{\eqref{eq:HBA}}{=}
        \begin{tikzpicture}[baseline={(0,0.8)}]
            \node at (0,-0.2) {$V^\ast$}; \node at (0.5,-0.2) {$V$}; \node[green] at (1,-0.2) {$A$}; \pic[rightstr = green] at (0.5,0) {braid}; \pic[rightstr = green] at (0,0.5) {invbraid}; \draw (0,0) -- (0,0.5); \draw[-Latex] (0,0.5) -- (0,0.2);
            \pic[rightstr = green] at (1,0.5) {braid}; \filldraw[gray] (1.25,0.75) circle[radius=0.1]; \pic[leftstr = green] at (1,1) {braid};
            \centerarc[](0.75,1.5)(0:180:0.25); \draw (0.5,1) -- (0.5,1.5);
            \centerarc[green](0.75,1.5)(0:180:0.75); \draw[green] (1.5,0.4) circle[radius=0.1]; \draw[green] (0,1) -- (0,1.5); \draw[green,-Latex] (0.75,2.25) -- (0.75,3); \node[green] at (0.75,3.2) {$A$};
        \end{tikzpicture} \nonumber\\
        &= 
        \begin{tikzpicture}[baseline={(0,0.8)}]
            \node at (0,-0.2) {$V^\ast$}; \node at (0.5,-0.2) {$V$}; \node[green] at (1,-0.2) {$A$}; \pic[rightstr = green] at (0.5,0) {braid}; \pic[rightstr = green] at (0,0.5) {invbraid}; \draw (0,0) -- (0,0.5); \draw[-Latex] (0,0.5) -- (0,0.2); \draw[-Latex] (1,0.5) -- (1,0.8); \draw (1,0.5) -- (1,1); \draw[gray] (0.3,1) -- (1.2,1); \draw (0.75,1) -- (0.75,1.5); \centerarc[green](0.375,1.5)(0:180:0.375); \filldraw[red] (0.75,1.5) circle[radius=0.1]; \node[red] at (1.05,1.5) {$\mu$}; \draw[green] (0,1) -- (0,1.5); \draw[green,-Latex] (0.375,1.875) -- (0.375,2.625); \node[green] at (0.375,2.825) {$A$};
        \end{tikzpicture}
    \end{align}    
    Now we consider the other direction. Given a quantum moment map $\mu:\mathcal{L}\to A$, we define a natural transformation $\gamma^A:\blank\otimes A\Rightarrow A\otimes\blank$ via
    \begin{equation}\label{eq:HBfromDinat}
        \gamma^A_V := 
        \begin{tikzpicture}[baseline={(0,0.8)}]
            \node at (0.5,-0.2) {$V$}; \draw[-Latex] (0.5,0) -- (0.5,0.8); \draw (0.5,1) -- (0.5,0); \draw[gray, ultra thick] (-0.2,1) -- (0.7,1); \draw (0.25,1) -- (0.25,1.5);
            \centerarc[green](0.75,1.5)(0:180:0.5); \filldraw[red] (0.25,1.5) circle[radius=0.1]; \draw[green] (1.25,0) -- (1.25,1.5); \node[green] at (1.25,-0.2) {$A$};
            \centerarc[](-0.25,1)(180:360:0.25); \draw[smooth] plot coordinates {(-0.5,1) (-0.45,1.2) (0.2,1.8) (0.25,2)};
            \pic[rightstr = green] at (0.25,2) {invbraid};
            \draw[green,-Latex] (0.25,2.5) -- (0.25,3); \node[green] at (0.25,3.2) {$A$}; \draw[-Latex] (0.75,2.5) -- (0.75,3); \node at (0.75,3.2) {$V$};
        \end{tikzpicture}
    \end{equation}
    Dinaturality of $i:(\blank)^\ast\otimes\bblank\,\ddot{\to}\,\mathcal{L}$ implies that
$\gamma^A$ is natural. A straightforward string diagram computation shows that $\gamma^A$ satisfies the hexagon \eqref{eq:half-hexagon} and conditions \eqref{eq:HBA}. A computation using \eqref{eq:LalgA}, \eqref{eq:LalgB}, dinaturality of $i$ and the zig-zag identities shows that
    \begin{equation}
        (\gamma^A_V)^{-1} := \begin{tikzpicture}[baseline={(0,0.8)}]
            \centerarc[](0.75,0)(180:360:0.25);
            \pic at (1,0) {invbraid};
            \centerarc[](1.25,0.5)(0:180:0.25);
            \draw (1.5,-0.5) -- (1.5,0);
            \draw[green] (0,-0.5) -- (0,0);
            \pic[leftstr = green] at (0,0) {braid}; 
            \node[green] at (0,-0.7) {$A$}; \node at (1.5,-0.7) {$V$};
            \draw[gray, ultra thick] (-0.7,1.5) -- (0.2,1.5);
            \draw (0,0.5) -- (0,1); \pic at (-0.5,1) {braid}; \draw[-Latex] (1.5,-0.5) -- (1.5,0); \centerarc[](-0.75,1)(180:360:0.25);
            \draw[green] (0.5,0.5) -- (0.5,2);
            \draw (-0.25,1.5) -- (-0.25,2);
            \draw[green,-Latex] (0.125,2.375) -- (0.125,3);
            \node[green] at (0.125,3.2) {$A$};
            \draw[-Latex] (-1,1) -- (-1,3); \node at (-1,3.2) {$V$}; 
            \centerarc[green](0.125,2)(0:180:0.375);
            \filldraw[red] (-0.25,2) circle[radius=0.1];
        \end{tikzpicture}
    \end{equation} 
is inverse to $\gamma^A_V$.
    
    Finally, we show that these constructions are mutually inverse. Inserting the quantum moment map obtained from the half-braiding \eqref{eq:DinatFromHB} into \eqref{eq:HBfromDinat} results in the original half-braiding:
    \begin{equation}
        \begin{tikzpicture}[baseline={(0,0.8)}]
            \node at (0.5,-0.2) {$V$}; \draw[-Latex] (0.5,0) -- (0.5,0.8); \draw (0.5,1) -- (0.5,0); \draw[gray, ultra thick] (-0.2,1) -- (0.7,1); \draw (0.25,1) -- (0.25,1.5);
            \centerarc[green](0.75,1.5)(0:180:0.5); \filldraw[red] (0.25,1.5) circle[radius=0.1];; \draw[green] (1.25,0) -- (1.25,1.5); \node[green] at (1.25,-0.2) {$A$};
            \centerarc[](-0.25,1)(180:360:0.25); \draw[smooth] plot coordinates {(-0.5,1) (-0.45,1.2) (0.2,1.8) (0.25,2)};
            \pic[rightstr = green] at (0.25,2) {invbraid};
            \draw[green,-Latex] (0.25,2.5) -- (0.25,3); \node[green] at (0.25,3.2) {$A$}; \draw[-Latex] (0.75,2.5) -- (0.75,3); \node at (0.75,3.2) {$V$};
        \end{tikzpicture}
        =
        \begin{tikzpicture}[baseline={(0,0.8)}]
            \node at (0,-0.2) {$V$}; \draw (0,0) -- (0,0.25); \pic[rightstr = green] at (0,0.25) {braid}; \filldraw[gray] (0.25,0.5) circle[radius=0.1]; \pic[leftstr = green] at (0,0.75) {braid}; \draw[green] (0.5,0.15) circle[radius=0.1];
            \centerarc[green](0.75,1.25)(0:180:0.25); \draw[green] (1,0) -- (1,1.25); \node[green] at (1,-0.2) {$A$};
            \pic[rightstr = green] at (0.25,1.5) {invbraid}; 
            \centerarc[](0.125,1.5)(270:360:0.125); \centerarc[](0.125,1.25)(90:180:0.125);
            \draw[green,-Latex] (0.25,2) -- (0.25,2.75); \node[green] at (0.25,2.95) {$A$}; \draw[-Latex] (0.75,2) -- (0.75,2.75);  \node at (0.75,2.95) {$V$};
        \end{tikzpicture}
        \overset{\eqref{eq:HBA}}{=}
        \begin{tikzpicture}[baseline={(0,0.8)}]
            \node at (0,-0.2) {$V$}; \pic[rightstr = green] at (0,0.5) {braid}; \filldraw[gray] (0.25,0.75) circle[radius=0.1]; \centerarc[green](0.5,0.35)(0:180:0.15); \draw[green] (0.35,0.25) circle[radius=0.1]; \draw (0,0) -- (0,0.5); \draw[green] (0.65,0) -- (0.65,0.35); \node[green] at (0.65,-0.2) {$A$};
            \draw[green,-Latex] (0,1) -- (0,1.3); \node[green] at (0,1.5) {$A$};
            \draw[-Latex] (0.5,1) -- (0.5,1.3); \node at (0.5,1.5) {$V$};
        \end{tikzpicture}
        =
        \begin{tikzpicture}[baseline={(0,0.25)}]
            \pic[rightstr=green] at (0,0) {braid}; \filldraw[gray] (0.25,0.25) circle[radius=0.1]; \draw[green,-Latex] (0,0.5) -- (0,0.75); \draw[-Latex] (0.5,0.5) -- (0.5,0.75); \node at (0,-0.2) {$V$}; \node[green] at (0.5,-0.2) {$A$}; \node[green] at (0,0.95) {$A$}; \node at (0.5,0.95) {$V$};
        \end{tikzpicture} \;.
    \end{equation}
A similar computation shows that inserting the half-braiding \eqref{eq:HBfromDinat} into the defining relation \eqref{eq:DinatFromHB} the quantum moment map results in the defining relation for the original quantum moment map.
\end{proof}

Combining Propositions~\ref{Prop:BrModCatiffHalfBr} and~\ref{prop:QMMiffHalfBr}, we arrive at the following result. 
\begin{corollary}[{\cite[Cor.\,4.11]{BBJ16}}]\label{Cor:QMMModuleBr}
Assume $\cat{C}$ is ribbon and contains the coend $\mathcal{L}$. Let $(A,m,u)$ be an algebra in $\cat{C}$. There is a bijective correspondence between module braidings for $\cat{C}_A$ and quantum moment maps $\mu:\mathcal{L}\to A$. In this correspondence, a module braiding is constructed from a quantum moment map via
    \begin{equation}
        e_{V,M} = 
        \begin{tikzpicture}[baseline={(0,0.8)}]
            \node at (0,-0.2) {$V$}; \node[blue] at (0.5,-0.2) {$M$};
            \pic[rightstr = blue] at (0,0) {braid};
            \draw[gray] (0.8,1.5) -- (1.7,1.5);
            \draw[smooth] plot coordinates {(0.5,0.5) (0.6,0.7) (1.4,1.3) (1.5,1.5)};
            \centerarc[](0.75,1.5)(180:360:0.25);
            \draw[blue] (0,0.5) -- (0,1.5);
            \pic[leftstr = blue] at (0,1.5) {braid};
            \centerarc[](0.5,1.5)(0:45:0.75);
            \draw[blue,-Latex] (0.5,2) -- (0.5,3); \node[blue] at (0.5,3.2) {$M$};
            \centerarc[green](0.5,1.5)(45:90:0.75);
            \filldraw[red] (1.03,2.03) circle[radius=0.1];
            \draw[-Latex] (0,2) -- (0,3); \node at (0,3.2) {$V$};
        \end{tikzpicture} \;.
    \end{equation}
\end{corollary}

\begin{remark}
    The pair $(A,\gamma^A)$ is an object in the Drinfeld centre $Z(\cat{C})$ of $\cat{C}$. However, equipping the Drinfeld centre with the canonical monoidal structure which turns it into a braided monoidal category (see \cite[Sec.\,8.5]{EGNO15}), the conditions \eqref{eq:HBA} do not allow one to lift the algebra structure to an algebra structure in the Drinfeld centre. In \cite[Cor.\,4.6]{BBJ16} (see also \cite[Def.\,2.47]{JR21}) it was observed that if $\cat{C}$ is ribbon, one can use the equivalence between $\cat{C}_\mathcal{L}$ and $Z(\cat{C})$ to define an alternative monoidal structure on $Z(\cat{C})$ via the relative tensor product over $\cat{L}$, such that a half-braided algebra in $\cat{C}$ lifts to an algebra in $Z(\cat{C})$. The category $\cat{C}_\mathcal{L}$ with this monoidal structure is sometimes called the \emph{Harish-Chandra category} \cite[Sec.\,2.5]{GJS19}.
\end{remark}

\subsection{Small braidings}\label{ssec:SmallBraidings}

Let $\cat{C}$ and $(\cat{D},\otimes,\1,\alpha,l,r)$ be monoidal categories, 
$\cat{C}\xhookrightarrow{\iota}\cat{D}$ a full monoidal subcategory. Denote by $\iota_{U,V}:\iota(U\otimes V) \to \iota U\otimes \iota V$ the monoidal structure of $\iota$.
\begin{definition}\label{def:SmallBraiding}
    A \emph{small braiding for $\cat{C}\xhookrightarrow{\iota}\cat{D}$} is the data of two natural isomorphisms\footnote{
One way to remember the convention for the arrow on $\sigma$ is that for the source objects, it points the same way as $\iota$, namely $\cat{C}\to\cat{D}$. For example, $\lc_{U,X}$ has $U \in \cat{C}$ and $X \in \cat{D}$.
    }
    \begin{equation}
        \begin{tikzcd}
        	{\cat{D}\times\cat{C}} && \cat{D}
        	\arrow["\blank\otimes\iota(\bblank)"{name=U, anchor=center, above}, bend left=50, from=1-1, to=1-3]
            \arrow["\iota(\bblank)\otimes\blank"{name=D, anchor=center, below}, bend right=50, from=1-1, to=1-3]
        	\arrow["\rc", Rightarrow, from=U, to=D]
        \end{tikzcd} \;,\quad 
        \begin{tikzcd}
        	{\cat{C}\times\cat{D}} && \cat{D}
        	\arrow["\iota(\blank)\otimes \bblank"{name=U, anchor=center, above}, bend left=50, from=1-1, to=1-3]
            \arrow["\bblank\otimes\iota(\blank)"{name=D, anchor=center, below}, bend right=50, from=1-1, to=1-3]
        	\arrow["\lc", Rightarrow, from=U, to=D]
        \end{tikzcd}
    \end{equation}
    which agree on $\cat{C}$, (for all $U,V\in\cat{C}$ one has $\rc_{\iota U,V} = \lc_{U,\iota V}$) and satisfy the hexagon constraints, namely, for all $U,V\in\cat{C}$ and $X,Y\in\cat{D}$ the following diagrams commute 
    \begin{subequations}
    \begin{equation}\label{eq:SmallHexA}
        \begin{tikzcd}
        	& {\scriptstyle (\iota U\otimes X)\otimes \iota V} & {\scriptstyle \iota U\otimes(X\otimes \iota V)} \\
        	{\scriptstyle (X\otimes \iota U)\otimes \iota V} &&& {\scriptstyle \iota U\otimes(\iota V\otimes X)} \\
        	& {\scriptstyle X\otimes \iota(U\otimes V)} & {\scriptstyle \iota(U\otimes V)\otimes X}
        	\arrow["{\alpha}", from=1-2, to=1-3]
        	\arrow["{\iota U\otimes\rc_{X,V}}", from=1-3, to=2-4]
        	\arrow["{\rc_{X,U}\otimes \iota V}", from=2-1, to=1-2]
        	\arrow["{(X\otimes \iota_{U,V}^{-1})\circ \alpha}"', from=2-1, to=3-2]
        	\arrow["{\rc_{X,U\otimes V}}"', from=3-2, to=3-3]
        	\arrow["{\alpha\circ(\iota_{U,V}\otimes X)}"', from=3-3, to=2-4]
        \end{tikzcd}
    \end{equation}
    \begin{equation}\label{eq:SmallHexB}
        \begin{tikzcd}
        	& {\scriptstyle \iota U\otimes(X\otimes Y)} & {\scriptstyle (\iota U\otimes X)\otimes Y} & \\
        	{\scriptstyle (X\otimes Y)\otimes \iota U} &&& {\scriptstyle (X\otimes \iota U)\otimes Y} \\
        	& {\scriptstyle X\otimes(Y\otimes \iota U)} & {\scriptstyle X\otimes(\iota U\otimes Y)}
        	\arrow["{\alpha^{-1}}", from=1-2, to=1-3]
        	\arrow["{\rc_{X\otimes Y,U}}", from=2-1, to=1-2]
        	\arrow["\alpha"', from=2-1, to=3-2]
        	\arrow["{\rc_{X, U}\otimes Y}"', from=2-4, to=1-3]
        	\arrow["{X\otimes\rc_{Y,U}}"', from=3-2, to=3-3]
        	\arrow["{\alpha^{-1}}"', from=3-3, to=2-4]
        \end{tikzcd}
    \end{equation}
    \begin{equation}\label{eq:SmallHexC}
        \begin{tikzcd}
        	& {\scriptstyle \iota U\otimes(X\otimes \iota V)} & {\scriptstyle (\iota U\otimes X)\otimes \iota V} \\
        	{\scriptstyle \iota U\otimes(\iota V\otimes X)} &&& {\scriptstyle (X\otimes \iota U)\otimes \iota V} \\
        	& {\scriptstyle \iota(U\otimes V)\otimes X} & {\scriptstyle X\otimes\iota(U\otimes V)}
        	\arrow["{\alpha^{-1}}", from=1-2, to=1-3]
        	\arrow["{\lc_{U,X}\otimes \iota V}", from=1-3, to=2-4]
        	\arrow["{\iota U\otimes\lc_{V,X}}", from=2-1, to=1-2]
        	\arrow["{(\iota_{U,V}^{-1}\otimes X)\circ\alpha^{-1}}"', from=2-1, to=3-2]
        	\arrow["{\lc_{U\otimes V,X}}"', from=3-2, to=3-3]        	\arrow["{\alpha^{-1}\circ(X\otimes\iota_{U,V})}"', from=3-3, to=2-4]
        \end{tikzcd}
    \end{equation}
    \begin{equation}\label{eq:SmallHexD}
        \begin{tikzcd}
        	& {\scriptstyle (X\otimes Y)\otimes \iota U} & {\scriptstyle X\otimes(Y\otimes \iota U)} & \\
        	{\scriptstyle \iota U\otimes(X\otimes Y)} &&& {\scriptstyle X\otimes(\iota U\otimes Y)} \\
        	& {\scriptstyle (\iota U\otimes X)\otimes Y} & {\scriptstyle (X\otimes \iota U)\otimes Y}
        	\arrow["\alpha", from=1-2, to=1-3]
        	\arrow["{\lc_{U,X\otimes Y}}", from=2-1, to=1-2]
        	\arrow["\alpha^{-1}"', from=2-1, to=3-2]
        	\arrow["{X\otimes\lc_{U,Y}}"', from=2-4, to=1-3]
        	\arrow["{\lc_{U,X}\otimes Y}"', from=3-2, to=3-3]
        	\arrow["\alpha"', from=3-3, to=2-4]
        \end{tikzcd} \;.
    \end{equation}
    \end{subequations}
\end{definition}

We will think of $\cat{C}$ as a subcategory of $\cat{D}$ and so in many places we no longer write out the embedding $\iota$ explicitly.
We depict a small braiding similarly to a regular braiding by using a thicker line to indicate the object which comes from $\cat{D}$:
\begin{equation}
    \rc_{X,V} = \begin{tikzpicture}[baseline={(0,0.25)}]
        \pic[leftstr=thick black] at (0,0) {braid}; \draw[-Latex] (0,0.5) -- (0,0.75); \draw[thick black, -Latex] (0.5,0.5) -- (0.5,0.75); \node at (0,-0.2) {$X$}; \node at (0.5,-0.2) {$V$}; \node at (0,0.95) {$V$}; \node at (0.5,0.95) {$X$};
    \end{tikzpicture}\;,\quad \lc_{V,X} =  \begin{tikzpicture}[baseline={(0,0.25)}]
        \pic[rightstr=thick black] at (0,0) {braid}; \draw[thick black,-Latex] (0,0.5) -- (0,0.75); \draw[-Latex] (0.5,0.5) -- (0.5,0.75); \node at (0,-0.2) {$V$}; \node at (0.5,-0.2) {$X$}; \node at (0,0.95) {$X$}; \node at (0.5,0.95) {$V$};
    \end{tikzpicture} \;.
\end{equation}

\begin{lemma}\label{lemma:SmallBraidingGivesBraiding}
    A small braiding on $(\cat{C}\xhookrightarrow{\iota}\cat{D},\rc,\lc)$ endows $\cat{C}$ with the structure of a braided monoidal category, and $l_X\circ\rc_{X,\1} = r_X$ and $r_X\circ\lc_{\1,X} = l_X$ for all $X\in\cat{D}$.
\end{lemma}
\begin{proof}
    Since $\rc$ and $\lc$ agree on $\cat{C}$, their respective hexagons together define a braiding for $\cat{C}$. The unitality property $l_X\circ\rc_{X,\1} = r_X$ follows from setting $U = \1$ in \eqref{eq:SmallHexA} and inserting unitors in the interior as follows (cf.\ \cite[Prop.\,XIII.1.2]{Ka95}):
    \begin{equation}
        \begin{tikzcd}[column sep=3em, row sep=3em]
        	& {\scriptstyle (\iota \1\otimes X)\otimes\iota V} & {\scriptstyle \iota\1\otimes(X\otimes \iota V)} & \\
        	{\scriptstyle (X\otimes\iota\1)\otimes\iota V} & {\scriptstyle X\otimes\iota V} & {\scriptstyle \iota V\otimes X} & {\scriptstyle \iota \1\otimes(\iota V\otimes X)} \\
        	& {\scriptstyle X\otimes(\iota\1\otimes\iota V)} & {\scriptstyle (\iota\1\otimes\iota V)\otimes X}
        	\arrow["\alpha", from=1-2, to=1-3]
        	\arrow["{l_X\otimes \Id_{\iota V}}"', from=1-2, to=2-2]
        	\arrow["{l_{X\otimes \iota V}}"', from=1-3, to=2-2]
        	\arrow["{\iota\1\otimes\rc_{X,V}}", from=1-3, to=2-4]
        	\arrow["{\rc_{X,\1}\otimes \Id_{\iota V}}", from=2-1, to=1-2]
        	\arrow["{r_X\otimes \Id_{\iota V}}"', from=2-1, to=2-2]
        	\arrow["\alpha"', from=2-1, to=3-2]
        	\arrow["{\rc_{X,V}}"', from=2-2, to=2-3]
        	\arrow["{l_{\iota V\otimes X}}", from=2-4, to=2-3]
        	\arrow["{\Id_X\otimes l_{\iota V}}"', from=3-2, to=2-2]
        	\arrow["{\rc_{X,\1\otimes V}}"', from=3-2, to=3-3]
        	\arrow["{l_{\iota V}\otimes \Id_X}", from=3-3, to=2-3]
        	\arrow["\alpha"', from=3-3, to=2-4]
        \end{tikzcd}
    \end{equation}
    The square in the top half of the hexagon commutes via naturality of the unitor, and the square in the bottom half commutes from naturality of $\rc$. The two triangles adjacent to the top left triangle as well as the bottom right triangle commute from the definition of the unitor. Thus the top left triangle commutes, which gives us the result upon setting $V=\1$ and using the right unitor. The property $r_X\circ\lc_{\1,X} = l_X$ can be proven by filling in the hexagon \eqref{eq:SmallHexC} similarly.
\end{proof}
In all sections that follow we denote with $\sigma$ the underlying braiding in $\cat{C}$ coming from a small braiding $(\rc,\lc)$ on $\cat{C}\xhookrightarrow{\iota}\cat{D}$. 
We also suppress coherence morphisms.

\begin{proposition}\label{prop:ClosedSubcatBraidedModuleCat}
    Let $\cat{M}$ be a full subcategory of $\cat{D}$ such that $\iota V\otimes M$ is an object in $\cat{M}$ for all objects $V$ in $\cat{C}$ and $M$ in $\cat{M}$. Then $(\cat{M},\rc\circ\lc)$ is a braided left $\cat{C}$-module category with action functor $\iota\otimes\blank$.
\end{proposition} 
\begin{proof}
    Let $\cat{M}$ be as above. The $\cat{C}$-module structure is given by the monoidal structure for $\cat{D}$. Let $U,V\in\cat{C}$ and $M\in\cat{M}$. The identities for $\rc\circ\lc$ we need to verify are \eqref{eq:BrM1} and \eqref{eq:BrM2}. The identity \eqref{eq:BrM1} follows from
    \begin{equation}
        \begin{tikzpicture}[baseline={(0,0.65)}]
             \pic[moduleA = thick blue] at (0,0) {modbrots}; \draw (-0.5,-0.5) -- (-0.5,0); \node at (-0.5,-0.7) {$U$}; \draw (0,-0.5) -- (0,0); \node at (0,-0.7) {$V$}; \draw[thick blue] (0.5,-0.5) -- (0.5,0); \node[blue] at (0.5,-0.7) {$M$}; \draw[-Latex] (0,1.5) -- (0,2); \node at (-0.5,2.2) {$U$}; \draw[-Latex] (-0.5,1.5) -- (-0.5,2); \node at (0,2.2) {$V$}; \draw[-Latex, thick blue] (0.5,1.5) -- (0.5,2); \node[blue] at (0.5,2.2) {$M$};
        \end{tikzpicture} 
        = 
        \begin{tikzpicture}[baseline={(0,1)}]
            \pic[rightstrB = thick blue] at (0,0) {braidots}; \pic[leftstrB = thick blue] at (0,1) {braidtso}; \draw[-Latex] (0,2) -- (0,2.3); \draw[-Latex] (0.5,2) -- (0.5,2.3); \draw[-Latex, thick blue] (1,2) -- (1,2.3); \node at (0,2.5) {$U$}; \node at (0.5,2.5) {$V$}; \node[blue] at (1,2.5) {$M$};
            \node at (0,-0.2) {$U$}; \node at (0.5,-0.2) {$V$}; \node[blue] at (1,-0.2) {$M$};
        \end{tikzpicture}
        =
        \begin{tikzpicture}[baseline={(0,0.4)}]
            \pic[rightstr = thick blue] at (0,0) {braid}; \pic[leftstr = thick blue] at (0,0.5) {braid}; \pic at (-0.5,1) {braid}; \pic at (-0.5,-0.5) {braid}; \draw (-0.5,0) -- (-0.5,1); \draw[-Latex] (-0.5,1.5) -- (-0.5,1.75); \draw[-Latex] (0,1.5) -- (0,1.75); \draw[thick blue] (0.5,-0.5) -- (0.5,0); \draw[thick blue,-Latex] (0.5,1) -- (0.5,1.75); \node at (-0.5,-0.7) {$U$}; \node at (0,-0.7) {$V$}; \node[blue] at (0.5,-0.7) {$M$}; \node at (-0.5,1.95) {$U$}; \node at (0,1.95) {$V$}; \node[blue] at (0.5,1.95) {$M$};
        \end{tikzpicture}
        =
        \begin{tikzpicture}[baseline={(0,0.65)}]
            \pic[module = thick blue] at (0,0) {modbr}; \pic at (-0.5,1.5) {braid}; \pic at (-0.5,-0.5) {braid}; \draw (-0.5,0) -- (-0.5,1.5); \draw[-Latex] (-0.5,2) -- (-0.5,2.25); \draw[-Latex] (0,2) -- (0,2.25); \draw[thick blue] (0.5,-0.5) -- (0.5,0); \draw[thick blue,-Latex] (0.5,1.5) -- (0.5,2.25); \node at (-0.5,-0.7) {$U$}; \node at (0,-0.7) {$V$}; \node[blue] at (0.5,-0.7) {$M$}; \node at (-0.5,2.45) {$U$}; \node at (0,2.45) {$V$}; \node[blue] at (0.5,2.45) {$M$};
        \end{tikzpicture} \;,
    \end{equation}
    where in the first equality we inserted the definition of the module braiding $\rc\circ\lc$, and we used the hexagons in Definition~\ref{def:SmallBraiding} in the second equality. The identity \eqref{eq:BrM2} follows from
    \begin{equation}
        \begin{tikzpicture}[baseline={(0,0.65)}]
            \pic[module = thick blue] at (0,0) {modbrtso}; \draw (-0.5,-0.3) -- (-0.5,0); \node at (-0.5,-0.5) {$U$}; \draw (0,-0.3) -- (0,0); \node at (0,-0.5) {$V$}; \draw[thick blue] (0.5,-0.3) -- (0.5,0); \node[blue] at (0.5,-0.5) {$M$};
            \draw[-Latex] (-0.5,1.5) -- (-0.5,2); \node at (-0.5,2.2) {$U$}; \draw[-Latex] (0,1.5) -- (0,2); \node at (0,2.2) {$V$}; \draw[thick blue,-Latex] (0.5,1.5) -- (0.5,2); \node[blue] at (0.5,2.2) {$M$};
        \end{tikzpicture}
        =
        \begin{tikzpicture}[baseline={(0,1)}]
            \pic[rightstr = thick blue] at (0,0) {braidtso}; \pic[leftstr = thick blue] at (0,1) {braidots}; \draw[-Latex] (0,2) -- (0,2.3); \draw[-Latex] (0.5,2) -- (0.5,2.3); \draw[-Latex, thick blue] (1,2) -- (1,2.3); \node at (0,2.5) {$U$}; \node at (0.5,2.5) {$V$}; \node[blue] at (1,2.5) {$M$};
            \node at (0,-0.2) {$U$}; \node at (0.5,-0.2) {$V$}; \node[blue] at (1,-0.2) {$M$};
        \end{tikzpicture}
        =
        \begin{tikzpicture}[baseline={(0,1.5)}]
            \node at (0,-0.2) {$U$}; \node at (0.5,-0.2) {$V$}; \node[blue] at (1,-0.2) {$M$};
            \draw[thick blue] (1,0) -- (1,0.5);
            \pic at (0,0) {invbraid}; \pic[rightstr = thick blue] at (0.5,0.5) {braid}; \pic[leftstr = thick blue] at (0.5,1) {braid};
            \draw (0,0.5) -- (0,1.5);
            \pic at (0,1.5) {braid};
            \pic[rightstr = thick blue] at (0.5,2) {braid}; \pic[leftstr = thick blue] at (0.5,2.5) {braid};
            \draw[thick blue] (1,1.5) -- (1,2);
            \draw[-Latex] (0,2) -- (0,3.3); \node at (0,3.5) {$U$}; \draw[-Latex] (0.5,3) -- (0.5,3.3); \node at (0.5,3.5) {$V$}; \draw[-Latex, thick blue] (1,3) -- (1,3.3); \node[blue] at (1,3.5) {$M$};
        \end{tikzpicture}
        =
        \begin{tikzpicture}[baseline={(0,1.6)}]
            \pic at (-0.5,-0.5) {invbraid}; \pic[module = thick blue] at (0,0) {modbr};  \pic at (-0.5,1.5) {braid}; \pic[module = thick blue] at (0,2) {modbr}; \node at (-0.5,-0.7) {$U$}; \draw (-0.5,0) -- (-0.5,1.5); \draw[-Latex] (-0.5,2) -- (-0.5,4); \node at (-0.5,4.2) {$U$}; \draw[-Latex] (0,3.5) -- (0,4); \node at (0,4.2) {$V$}; \node at (0,-0.7) {$V$}; \draw[thick blue] (0.5,-0.5) -- (0.5,0); \draw[thick blue] (0.5,1.5) -- (0.5,2); \draw[thick blue,-Latex] (0.5,3.5) -- (0.5,4); \node[blue] at (0.5,-0.7) {$M$}; \node[blue] at (0.5,4.2) {$M$};
        \end{tikzpicture}\;,
    \end{equation}
    with the first equality being the definition of the module braiding and the second following from the hexagons and naturality of the small braidings.
\end{proof}

Let $(A,m,u)$ be an algebra in $\cat{D}$. 
By virtue of $\cat{C}$ being a monoidal subcategory of $\cat{D}$, the category $\cat{D}_A$ is a $\cat{C}$-module category. The following proposition adapts one direction of Proposition~\ref{Prop:BrModCatiffHalfBr} in order to define a module braiding for $\cat{D}_A$ from a half-braided structure on $A$. 
The proof is exactly the same as in Proposition~\ref{Prop:BrModCatiffHalfBr}. For $A\in\cat{D}$, we define a half-braiding $\widehat{\gamma}^A$ to be a natural isomorphism $\iota\otimes A \Rightarrow A\otimes \iota$ which satisfies the hexagon
\eqref{eq:half-hexagon}
with $\iota$ and its monoidal structure inserted where appropriate.

\begin{proposition}\label{prop:SmallModuleBrIfHalfBr}
    Let $(A,m,u)$ be an algebra in $\cat{D}$. A half-braiding $\widehat{\gamma}^A:\iota\otimes A \Rightarrow A\otimes \iota$ satisfying
    \begin{equation}\label{eq:hHBA}   
        \begin{tikzpicture}[baseline={(0,0.65)}]
            \centerarc[thick green](0,0.3)(0:180:0.5); \draw[thick green] (-0.5,0) -- (-0.5,0.3); \draw[thick green] (0.5,0) -- (0.5,0.3);
            \draw[thick green] (0,0.8) -- (0,1); \pic[rightstr = thick green] at (-0.5,1) {braid};             \filldraw[gray] (-0.25,1.25) circle[radius=0.1]; \draw[thick green, -Latex] (-0.5,1.5) -- (-0.5,1.8); \draw[-Latex] (0,1.5) -- (0,1.8);
            \centerarc[](-0.75,1)(270:360:0.25); \centerarc[](-0.75,0.5)(90:180:0.25); \draw (-1,0) -- (-1,0.5);
            \node at (-1,-0.2) {$V$}; \node[green] at (-0.5,-0.2) {$A$}; \node[green] at (0.5,-0.2) {$A$};
            \node at (0,2) {$V$}; \node[green] at (-0.5,2) {$A$};
        \end{tikzpicture}
        =
        \begin{tikzpicture}[baseline={(0,0.65)}]
            \pic[rightstr = thick green] at (-0.5,0) {braid}; \filldraw[gray] (-0.25,0.25) circle[radius=0.1]; \pic[rightstr = thick green] at (0,0.5) {invbraid};
            \node at (-0.5,-0.2) {$V$}; \node[green] at (0,-0.2) {$A$}; \node[green] at (0.5,-0.2) {$A$};
            \draw[-Latex] (0.5,1) -- (0.5,1.8); \node at (0.5,2) {$V$}; \draw[thick green] (-0.5,0.5) -- (-0.5,1); \draw[thick green] (0.5,0) -- (0.5,0.5); \centerarc[thick green](-0.25,1)(0:180:0.25); \draw[thick green, -Latex] (-0.25,1.25) -- (-0.25,1.8); \node[green] at (-0.25,2) {$A$};
        \end{tikzpicture}
        =
        \begin{tikzpicture}[baseline={(0,0.65)}]
            \pic[rightstr = thick green] at (-0.5,0) {braid}; \pic[rightstr = thick green] at (0,0.5) {braid}; \filldraw[gray] (0.25,0.75) circle[radius=0.1];
            \node at (-0.5,-0.2) {$V$}; \node[green] at (0,-0.2) {$A$}; \node[green] at (0.5,-0.2) {$A$};
            \draw[-Latex] (0.5,1) -- (0.5,1.8); \node at (0.5,2) {$V$}; \draw[thick green] (-0.5,0.5) -- (-0.5,1); \draw[thick green] (0.5,0) -- (0.5,0.5); \centerarc[thick green](-0.25,1)(0:180:0.25); \draw[thick green, -Latex] (-0.25,1.25) -- (-0.25,1.8); \node[green] at (-0.25,2) {$A$};
        \end{tikzpicture}
    \end{equation}
    and
    \begin{equation}\label{eq:small-half-hexagon}
        \begin{tikzpicture}[baseline={(0,0.8)}]
            \draw (0,0) -- (0,0.2);
            \draw (0.5,0) -- (0.5,0.2);
            \draw[thick green] (1,0) -- (1,0.2);
            \pic[rightstr=thick green] at (0,0.2) {braidtso};
            \draw[thick green,-Latex] (0,1.2) -- (0,1.7);
            \draw[-Latex] (0.5,1.2) -- (0.5,1.7);
            \draw[-Latex] (1,1.2) -- (1,1.7);
            \node[green] at (1,-0.2) {$A$}; \node[green] at (0,1.9) {$A$};
            \node at (0,-0.2) {$U$}; \node at (0.5,1.9) {$U$};
            \node at (0.5,-0.2) {$V$}; \node at (1,1.9) {$V$};
            \filldraw[gray] (0.5,0.7) circle[radius=0.3];
        \end{tikzpicture}
        =
        \begin{tikzpicture}[baseline={(0,0.8)}]
            \draw (0,0) -- (0,0.7);
            \draw (0.5,0) -- (0.5,0.2);
            \draw[thick green] (1,0) -- (1,0.2);
            \pic[rightstr = thick green] at (0.5,0.2) {braid}; \filldraw[gray] (0.75,0.45) circle[radius=0.1];
            \pic[rightstr = thick green] at (0,0.7) {braid}; \filldraw[gray] (0.25,0.95) circle[radius=0.1];
            \draw[thick green,-Latex] (0,1.2) -- (0,1.7);
            \draw[-Latex] (0.5,1.2) -- (0.5,1.7);
            \draw[-Latex] (1,0.7) -- (1,1.7);
            \node[green] at (1,-0.2) {$A$}; \node[green] at (0,1.9) {$A$};
            \node at (0,-0.2) {$U$}; \node at (0.5,1.9) {$U$};
            \node at (0.5,-0.2) {$V$}; \node at (1,1.9) {$V$};
        \end{tikzpicture} \;.
    \end{equation}
    for all $U,V\in\cat{C}$, equips the $\cat{C}$-module category $\cat{D}_A$ with a module braiding given by
    \begin{equation}\label{eq:HCmodulebraiding}
        e_{V,M} := \begin{tikzpicture}[baseline={(0,0.9)}]
            \pic[module = thick blue, algebra = thick green] at (0,0) {HBAmodbr}; 
            \draw (0,-0.3) -- (0,0); \node at (0,-0.5) {$V$}; \draw[-Latex] (0,2) -- (0,2.5); \node at (0,2.7) {$V$};
            \draw[thick blue] (0.5,-0.3) -- (0.5,0); \node[blue] at (0.5,-0.5) {$M$}; \draw[thick blue,-Latex] (0.5,2) -- (0.5,2.5); \node[blue] at (0.5,2.7) {$M$};
        \end{tikzpicture}
    \end{equation}  
\end{proposition}
As in the previous subsection, we call the pair consisting of an algebra $(A,m,u)\in\cat{D}$ and a half-braiding $\widehat{\gamma}^A$ satisfying \eqref{eq:hHBA} a half-braided algebra in $\cat{D}$.
\begin{remark}
    It is not clear to us if one can construct a half-braided structure for $A$ given a module braiding for $\cat{D}_A$: The inverse for the construction out of which one would wish to obtain a half-braided structure is
    \begin{equation}
        \widehat{\gamma}^A_V = \begin{tikzpicture}[baseline={(0,0.9)}]
            \pic[module = thick green] at (0,0) {modbr}; \pic[rightstr = thick green] at (0,1.5) {invbraid}; \draw (0,-0.3) -- (0,0); \node at (0,-0.5) {$V$}; \draw[thick green] (0.5,-0.3) -- (0.5,0); \node[green] at (0.5,-0.5) {$A$}; \draw[thick green,-Latex] (0,2) -- (0,2.5); \node[green] at (0,2.7) {$A$}; \draw[-Latex] (0.5,2) -- (0.5,2.5); \node at (0.5,2.7) {$V$};
        \end{tikzpicture} \;.
    \end{equation}
    The first equality in \eqref{eq:hHBA} follows from each $e_{V,A}$ being a morphism in $\cat{D}_{A}$ and the fact that the small braidings are natural. We do not know how to prove the second equality from the basic properties of module braidings. 
One would want to use \eqref{eq:BrMpicturesA} with $U,V$ set to $V,A$, but that does not apply as $A$ is not in $\cat{C}$.
    That being said, we do not have an example of a module braiding for $\cat{D}_A$ that does not give $A$ a half-braided structure.
\end{remark}

\medskip

Assume that $\cat{C}$ is a ribbon category and that the coend 
\begin{equation}
    \mathcal{L} = \int^{V\in\cat{C}} \iota(V^\ast\otimes V) \;,
\end{equation}
of the functor $\iota((\blank)^\ast\otimes\bblank):\cat{C}^{\mathrm{op}}\times\cat{C} \to \cat{D}$ exists in $\cat{D}$. 
Assume further that $\blank\otimes V$ and $V\otimes\blank$ preserve $\mathcal{L}$ for all $V\in\cat{D}$. The algebra structure for $\mathcal{L}$ is the same as \eqref{eq:LalgA}, \eqref{eq:LalgB} with $\iota$ inserted where appropriate.
Our assumption on preservation under $\blank\otimes V$ and $V\otimes\blank$ ensures that the algebra structure is well-defined.
\begin{remark}
    Preservation of $\mathcal{L}$ under $\blank\otimes V$ and $V\otimes\blank$ is not automatic in general. The assumption holds if the tensor product preserves all colimits, or $\cat{D}$ is a finite tensor category in the sense of \cite[Ch.\,4]{EGNO15}. Our example in Section~\ref{ssec:QMMs from internal-homs}
does satisfy the additional assumption by virtue of being closed monoidal and braided.
\end{remark}
As we did in Proposition~\ref{prop:SmallModuleBrIfHalfBr}, we obtain one direction of the results in Proposition~\ref{prop:QMMiffHalfBr} and Corollary~\ref{Cor:QMMModuleBr}. The proof is exactly the same as in the appropriate directions of the previous subsection, with braidings being replaced by small braidings, so we omit it.
\begin{proposition}\label{Prop:hCQMM}
    Let $\cat{C}$ be ribbon, $\cat{D}$ contain the coend $\mathcal{L}$, and assume that $\blank\otimes V$ and $V\otimes\blank$ preserve it for all $V\in\cat{D}$. Let $(A,m,u)$ be an algebra in $\cat{D}$, and let $\mu:\mathcal{L}\to A$ be an algebra homomorphism in $\cat{D}$ satisfying 
    \begin{equation}\label{eq:DQMM}
        \begin{tikzpicture}[baseline={(0,0.8)}]
            \node at (0,-0.2) {$V^\ast$}; \node at (0.5,-0.2) {$V$}; \draw[-Latex] (0,1) -- (0,0.3); \draw (0,1) -- (0,0); \draw[-Latex] (0.5,0) -- (0.5,0.8); \draw (0.5,1) -- (0.5,0); \draw[gray] (-0.2,1) -- (0.7,1); \draw[thick black] (0.25,1) -- (0.25,1.5);
            \centerarc[thick green](0.75,1.5)(0:180:0.5); \filldraw[red] (0.25,1.5) circle[radius=0.1]; \node[red] at (0.55,1.5) {$\mu$}; \draw[thick green] (1.25,0) -- (1.25,1.5); \node[green] at (1.25,-0.2) {$A$}; \draw[thick green,-Latex] (0.75,2) -- (0.75,2.75); \node[green] at (0.75,2.95) {$A$};
        \end{tikzpicture}  
        =
        \begin{tikzpicture}[baseline={(0,0.8)}]
            \node at (0,-0.2) {$V^\ast$}; \node at (0.5,-0.2) {$V$}; \node[green] at (1,-0.2) {$A$}; \pic[rightstr = thick green] at (0.5,0) {braid}; \pic[rightstr = thick green] at (0,0.5) {invbraid}; \draw (0,0) -- (0,0.5); \draw[-Latex] (0,0.5) -- (0,0.2); \draw[-Latex] (1,0.5) -- (1,0.8); \draw (1,0.5) -- (1,1); \draw[gray] (0.3,1) -- (1.2,1); \draw[thick black] (0.75,1) -- (0.75,1.5); \centerarc[thick green](0.375,1.5)(0:180:0.375); \filldraw[red] (0.75,1.5) circle[radius=0.1]; \node[red] at (1.05,1.5) {$\mu$}; \draw[thick green] (0,1) -- (0,1.5); \draw[thick green,-Latex] (0.375,1.875) -- (0.375,2.625); \node[green] at (0.375,2.825) {$A$};
        \end{tikzpicture}
    \end{equation}
    for all $V\in\cat{C}$. Then the morphisms
    \begin{equation}
        e_{V,M} = 
        \begin{tikzpicture}[baseline={(0,0.8)}]
            \node at (0,-0.2) {$V$}; \node[blue] at (0.5,-0.2) {$M$};
            \pic[rightstr =  thick blue] at (0,0) {braid};
            \draw[gray] (0.8,1.5) -- (1.7,1.5);
            \draw[smooth] plot coordinates {(0.5,0.5) (0.6,0.7) (1.4,1.3) (1.5,1.5)};
            \centerarc[](0.75,1.5)(180:360:0.25);
            \draw[thick blue] (0,0.5) -- (0,1.5);
            \pic[leftstr = thick blue] at (0,1.5) {braid};
            \centerarc[thick black](0.5,1.5)(0:45:0.75);
            \draw[thick blue,-Latex] (0.5,2) -- (0.5,3); \node[blue] at (0.5,3.2) {$M$};
            \centerarc[thick green](0.5,1.5)(45:90:0.75);
            \filldraw[red] (1.03,2.03) circle[radius=0.1];
            \draw[-Latex] (0,2) -- (0,3); \node at (0,3.2) {$V$};
        \end{tikzpicture} 
    \end{equation}
    form a module braiding for $\cat{D}_A$. There is a bijective correspondence between such maps $\mu$ and half-braided structures for $A$ in the sense of Proposition~\ref{prop:SmallModuleBrIfHalfBr}.
\end{proposition}

As in the previous subsection, we refer to algebra maps $\mu:\mathcal{L}\to A$ in $\cat{D}$ satisfying \eqref{eq:DQMM} as quantum moment maps for $A$.

\section{Braided module categories for quantum groups}\label{sec:BrModCatsforQuantumGroups}

In this section we construct examples of braided module categories, half-braided algebra structures and quantum moment maps coming from quantum groups. First we discuss small braidings associated to quantum groups for quantum parameter not a root of unity and introduce several concrete examples. We then change settings, and discuss small braidings coming from the $\hbar$-adic quantum group $U_\hbar\mathfrak{g}$. In the $\hbar$-adic setting we apply the machinery of Appendix~\ref{app:closed} and the previous section to construct quantum moment maps from internal-hom algebras.

Throughout this section, fix a semisimple Lie algebra $\mathfrak{g}$ of rank $n$. Let $\{\alpha_i\}$ be the simple roots of $\mathfrak{g}$, $\{H_i\}$ the simple coroots and $P$ the weight lattice. Let $(\blank,\bblank)$ be the standard symmetric bilinear form on $\mathfrak{h}^\ast$ proportional to the Killing form on $\mathfrak{g}$, normalised so that $(\alpha,\alpha)=2$ for short roots, and $d_i = (\alpha_i,\alpha_i)/2$. Let $A = ((\alpha_i,\alpha_j)/d_i)$ be the Cartan matrix for $\mathfrak{g}$.

\subsection{Small braidings from unrolled quantum groups}\label{ssec:Small braidings for Uqg}

We work over $\C$ and with quantum parameter $q\in\C^\times$, for which we fix a logarithm $\log q \in \mathbb{C}$.
We denote $e^{a\log q}$ by $q^a$, where $a$ can be a number or an operator. We denote with $U_q^H\mathfrak{g}$ the unrolled quantum group based on $\mathfrak{g}$, which was originally defined for $\mathfrak{g} = \mathfrak{sl}_2$ in \cite{GPT07} and for an arbitrary semisimple Lie algebra in \cite{GP10}. Concretely, 
we use the definition of $U_q^H\mathfrak{g}$ from \cite[Def.\,3.1]{Ru19}, which allows for $q$ not a root of unity.

Let us briefly recall a few basic facts and notation. 
Denote by $U_q\mathfrak{h}$, $U_q\mathfrak{n}_+$ and $U_q\mathfrak{n}_-$ the subalgebras of $U_q^H\mathfrak{g}$ generated by $\{K_i^{\pm1}\}_{i=1}^n$, $\{E_i\}_{i=1}^n$ and $\{F_i\}_{i=1}^n$, respectively, and denote the generators of $U\mathfrak{h}$ with $\{H_i\}_{i=1}^n$. The translation between the notation of \cite{Ru19} and our notation is $K_i = K_i$, $K_{-i} = K_i^{-1}$, $X_i = E_i$ and $X_{-i}=  F_i$. The generators coming from $U_q\mathfrak{n}_{\pm}$ act nilpotently in finite-dimensional representations. Let $U_q\mathfrak{g}$ be the subalgebra of $U_q^H\mathfrak{g}$ generated by $U_q\mathfrak{h}$ and $U_q\mathfrak{n}_\pm$. This is the De Concini--Kac specialisation of $U_q\mathfrak{g}$ \cite{DK90} (see also \cite[p.\,6]{Ru19}), which agrees with the definition of $U_q\mathfrak{g}$ given in \cite[Sec.\,6.1.2]{KS97}.

Recall that a $U_q^H\mathfrak{g}$-module $V$ is a \emph{weight module} if there exist $\lambda\in\mathfrak{h}^\ast$ such that $V$ has a direct sum decomposition \cite[Sec.\,5.3]{GP10}
\begin{equation}\label{eq:weight module}
    V \cong \bigoplus_{\lambda\in\mathfrak{h}^\ast} V_\lambda \;,\quad V_\lambda = \{v: H_i\cdot v = \lambda(H_i) v \;,\quad K_i\cdot v = q^{d_i\lambda(H_i)}\cdot v\}\;.
\end{equation}
We do not demand the $V_\lambda$ to be finite-dimensional.
    
For applications in skein theory, we are interested in defining a small braiding for the category of so-called finite-dimensional type 1 modules for $U_q\mathfrak{g}$. Recall that a $U_q\mathfrak{g}$-module $V$ is of type 1 if there exist $\lambda\in P$ such that $V$ has a direct sum decomposition \cite[p.\,314]{CP94} (see also \cite[p.\,4]{Cha94})
\begin{equation}\label{eq: type 1}
        V \cong \bigoplus_{\lambda\in P} V_\lambda \;,\quad V_\lambda = \{v: K_i\cdot v = q^{d_i\lambda(H_i)} v \}\;.
\end{equation}
Let $\cat{D}^H$ be the full monoidal subcategory of $U_q^H\mathfrak{g}$-modules consisting of weight modules, $\cat{C}^H$ the full monoidal subcategory of $\cat{D}^H$ consisting of finite-dimensional modules, and $\cat{C}$ be the category of finite-dimensional $U_q\mathfrak{g}$-modules of type 1.

It is well-known that finite-dimensional $U_q\mathfrak{g}$-modules of type 1 form a braided monoidal category when $q$ is not a root of unity, see for example \cite[Prop.\,8.19]{KS97} (but there is no braiding for all finite-dimensional $U_q\mathfrak{g}$-modules; this is shown for $\mathfrak{g}=\mathfrak{sl}_2$ in \cite{WW26}). The braiding is constructed in \cite[Sec.\,8.3.3]{KS97} by defining an endomorphism of the tensor product of $U_q\mathfrak{g}$-modules which implements the action of the universal $R$-matrix for the $\hbar$-adic quantum group with the formal parameter being specialised to a complex number. Replacing the formal parameter with a complex number does not recover an element of the tensor square of $U_q\mathfrak{g}$, so one cannot obtain a universal $R$-matrix for $U_q\mathfrak{g}$ this way. To get around this, one defines an element $R$ of a completion $\widehat{U_q\mathfrak{g}\otimes U_q\mathfrak{g}}$ of the tensor square of $U_q\mathfrak{g}$ \cite[Thm.\,8.18]{KS97}, which contains infinite sums over powers of $U_q\mathfrak{n}_-\otimes U_q\mathfrak{n}_+\subset\widehat{U_q\mathfrak{g}\otimes U_q\mathfrak{g}}$ generators, and an endomorphism $\Phi$ of the tensor product of $U,V\in\cat{C}$ which acts as $u\otimes v \mapsto q^{(\lambda,\mu)} u\otimes v$ for $u\in U_\lambda$ and $v\in V_\mu$ \cite[Prop.\,8.19]{KS97}, such that
\begin{equation}\label{eq:Uqg braiding}
    u\otimes v \mapsto \Phi R(u\otimes v)
\end{equation}
implements the action of the $\hbar$-adic $R$-matrix with the formal variable specialised to a complex number. We call $R$ the $R$-matrix for $U_q\mathfrak{g}$ (taking care not to refer to it as a universal $R$-matrix). Composing this with the swap operation of vector spaces then implements a braiding.

We now outline how this braiding can be extended to finite-dimensional $U_q^H\mathfrak{g}$-weight modules. Using the duality between $\mathfrak{h}^\ast$ and $\mathfrak{h}$ the endomorphism $u\otimes v \mapsto q^{(\lambda,\mu)} u\otimes v$ can be realised in terms of the generators of $U\mathfrak{h}\subset U_q^H\mathfrak{g}$ \cite[Lem.\,40]{GP10},
\begin{equation}\label{eq:UqHg automorphism}
    u\otimes v \mapsto q^{\sum_{i,j=1}^nd_i(A^{-1})_{ij}H_i\otimes H_j}(u\otimes v) \;.
\end{equation}
This endomorphism is well-defined on the tensor product of any two weight modules for $U_q^H\mathfrak{g}$, regardless of whether or not their weight supports are integral. Since $U_q\mathfrak{g}$ is a subalgebra of $U_q^H\mathfrak{g}$, one can ask whether the map
\begin{equation}
    u\otimes v \mapsto q^{\sum_{i,j=1}^nd_i(A^{-1})_{ij}H_i\otimes H_j}R(u\otimes v)
\end{equation}
composed with the swap operation on $U\otimes V$ defines a braiding between finite-dimensional $U_q^H\mathfrak{g}$-weight modules. Because
\begin{equation}
    q^{\sum_{i,j=1}^nd_i(A^{-1})_{ij}H_i\otimes H_j}\Delta(H_k) = \Delta(H_k)q^{\sum_{i,j=1}^nd_i(A^{-1})_{ij}H_i\otimes H_j}
\end{equation}
for all $k = 1,\dots,n$, it suffices to check $R\Delta(H_i) = \Delta^\mathrm{op}(H_i)R$. A straightforward computation shows this to be true.
\begin{remark}
    The truncated $R$-matrix of \cite{GP10} can only be used to implement a braiding when $q$ is a root of unity and one restricts to a certain subcategory of $\cat{C}^H$. Alternatively, in the $q$ a root of unity case one can take the quotient of $U_q^H\mathfrak{g}$ by ideals generated by certain powers of the generators coming from $U_q\mathfrak{n}_\pm$ to define the restricted quantum group \cite{GPT07,CGP14,Ru19}, the modules over which admit a braiding in terms of a truncated $R$-matrix. 
\end{remark}

Removing the assumptions on the weight support of finite-dimensional $U_q^H\mathfrak{g}$-weight modules $U$ and $V$, the proof of \cite[Prop.\,8.19]{KS97} can be adapted to show that finite-dimensional $U_q^H\mathfrak{g}$-weight modules form a braided monoidal category when $q$ is not a root of unity: One only needs to replace $B_{U,V}$ of \cite[Prop.\,8.19]{KS97} with \eqref{eq:UqHg automorphism}. The obstruction to braiding $U_q^H\mathfrak{g}$-modules of arbitrary dimension is that the sum over powers of the $U_q\mathfrak{n}_-\otimes U_q\mathfrak{n}_+$ generators in the $R$-matrix will not truncate when acting on a tensor product of modules which are not locally finite as $U_q\mathfrak{n}_+$-  or $U_q\mathfrak{n}_-$-modules, see \cite[Thm.\,8.18]{KS97} for the explicit form of the $R$-matrix. The torsion-free modules of \cite{FHW14} are an example of such non-locally finite modules. However, the sum will truncate if at least one of the modules is finite-dimensional, so one can braid a finite-dimensional module with an infinite-dimensional module. We then arrive at our first example of a small braiding. Note that it follows from \cite[Thm.\,8.18]{KS97} that the hexagons \eqref{eq:SmallHexA}--\eqref{eq:SmallHexD} are satisfied and the small braiding morphisms are invertible.
\begin{proposition}\label{prop:Unrolled small braiding}
    For all $U\in\cat{C}^H$ and $X\in\cat{D}^H$, if $q$ is not a root of unity, then the morphisms
    \begin{subequations}
        \begin{equation}
            \lc_{U,X}(u\otimes x) = s_{U,X}\left(q^{\sum_{i,j=1}^nd_i(A^{-1})_{ij}H_i\otimes H_j}R(u\otimes x) \right)
        \end{equation}
        \begin{equation}
            \rc_{X,U}(x\otimes u) = s_{X,U}\left(q^{\sum_{i,j=1}^nd_i(A^{-1})_{ij}H_i\otimes H_j}R(x\otimes u)\right) \;,
        \end{equation}
    \end{subequations}
    where $u\in U$, $x\in X$, $R$ is the $R$-matrix for $U_q\mathfrak{g}$ and $s$ is the swap operation on the tensor product of vector spaces, define a small braiding on $\cat{C}^H\hookrightarrow\cat{D}^H$.
\end{proposition}

The following proposition allows us to define a small braiding between finite-dimensional $U_q\mathfrak{g}$-modules of type 1 and $U_q^H\mathfrak{g}$-weight modules.

\begin{proposition}\label{prop: rolled type 1 includes into unrolled}
    Assume $q\in\C^\times$ is not a root of unity. Let $\cat{C}$ be the category of finite-dimensional $U_q\mathfrak{g}$-modules of type 1. There is a full braided monoidal inclusion of categories $\iota:\cat{C}\hookrightarrow\cat{C}^H$.
\end{proposition}
\begin{proof}
    If $q$ is not a root of unity the decomposition \eqref{eq: type 1} of $V\in\cat{C}$ is unique as each $\lambda(H_i)\in\Z$ when $\lambda\in P$. Define the functor $\iota:\cat{C}\to\cat{C}^H$ 
    to leave the underlying vector spaces of the objects and linear maps of the morphisms the same, and define an action of the generators of $U\mathfrak{h}$ to be $H_i\cdot v = \lambda(H_i)v$ for $v\in V_\lambda$. The uniqueness of the decomposition \eqref{eq: type 1} implies that linear maps which intertwine the $U_q\mathfrak{h}$ generators also intertwine the generators of $U\mathfrak{h}$, so $\iota$ is a well-defined functor. Furthermore the uniqueness of the decomposition implies such maps are in bijection, so $\iota$ is fully faithful. Let us now show $\iota$ is monoidal. Let $u\in U_\mu$ and $v\in V_\lambda$ for $U,V\in\cat{C}$. Consider $u\otimes v \in (U\otimes V)_{\mu+\lambda}$, on which the action of $K_i\in U_q\mathfrak{h}$ is given by
    \begin{equation}
        K_i\cdot(u\otimes v) = q^{d_i\mu(H_i) + d_i\lambda(H_i)}u\otimes v = q^{d_i(\mu+\lambda)(H_i)} u\otimes v \;.
    \end{equation}
    The action of $H_i\in U\mathfrak{h}$ on $u\otimes v$ is
    \begin{equation}
        H_i\cdot(u\otimes v) = H_i\cdot u\otimes v + u\otimes H_i\cdot v = (\mu+\lambda)(H_i)u\otimes v \;,
    \end{equation}
    as desired. We see that this inclusion is braided because of the identification of the endomorphism $\Phi$ in \eqref{eq:Uqg braiding} with $q^{\sum_{i,j=1}^nd_i(A^{-1})_{ij}H_i\otimes H_j}$ \cite[Lem.\,40]{GP10}.
\end{proof}

Propositions~\ref{prop:Unrolled small braiding} and~\ref{prop: rolled type 1 includes into unrolled} imply the main theorem of this section:

\begin{theorem}\label{thm:SmallBraidingUqg}
    Let $\cat{C}$ be the category of finite-dimensional $U_q\mathfrak{g}$-modules of type 1, $\cat{D}^H$ be the category of $U_q^H\mathfrak{g}$-weight modules, and $\cat{C}^H$ the full subcategory of $\cat{D}^H$ of finite-dimensional modules. For $q\in\C^\times$ not a root of unity, there is a small braiding on the monoidal inclusion $\cat{C}\xhookrightarrow{\iota}\cat{C}^H\hookrightarrow\cat{D}^H$ given by
    \begin{subequations}
        \begin{equation}
            \lc_{U,X}(u\otimes x) = s_{U,X}\Bigl( q^{\sum_{i,j=1}^nd_i(A^{-1})_{ij}H_i\otimes H_j}R(u\otimes x) \Bigr)
        \end{equation}
        \begin{equation}
            \rc_{X,U}(x\otimes u) = s_{X,U}\Bigl(q^{\sum_{i,j=1}^nd_i(A^{-1})_{ij}H_i\otimes H_j}R(x\otimes u)\Bigr) \;,
        \end{equation}
    \end{subequations}
    for all $U\in\cat{C}$, $X\in\cat{D}^H$, $u\in U$ and $x\in X$, where $R$ is the $R$-matrix for $U_q\mathfrak{g}$ and $s$ is the swap operation on the tensor product of vector spaces.
\end{theorem}

\subsection{Examples for quantum parameter not a root of unity}\label{ssec:Examples for generic q}

We fix $q\in\C^\times$ not a root of unity and define $\cat{C}$ and $\cat{D}^H$ as below \eqref{eq: type 1}. We now study subcategories of $\cat{D}^H$ closed under tensoring with objects in $\cat{C}$ in order to construct examples of braided module categories from the small braiding on $\cat{C}\hookrightarrow\cat{D}^H$ in Theorem~\ref{thm:SmallBraidingUqg}.

In the examples that follow, the module braiding is given by the monodromy $\rc\circ\lc$, as per the prescription in Proposition~\ref{prop:ClosedSubcatBraidedModuleCat}. We understand the following example to be folklore and have not seen it written down anywhere in the literature. It can be viewed as an unrolled version of the quantum category $\cat{O}$ of \cite{AM11}.
\begin{example}\label{ex:O}
    Define $\cat{O}^H_q$ to be the full subcategory of $\cat{D}^H$ consisting of modules which are finitely generated and in which $U_q\mathfrak{n}_+$ acts locally nilpotently. Tensoring with a finite-dimensional module preserves finite generation and local $U_q\mathfrak{n}_+$-nilpotency, so $\cat{O}_q^H$ is closed under tensoring with $\cat{C}$ and thus a braided $\cat{C}$-module category.
\end{example}

For each pair $\lambda,\mu\in\mathfrak{h}^\ast$, say $\lambda\sim\mu$ if $\lambda-\mu\in P$. Let $W\in\cat{C}$ and $V\in\cat{D}^H$. Because the weight support of $\iota W$ is in $P$, for all weights $\mu$ in the support of $\iota W\otimes V$ there exists $\lambda$ in the support of $V$ such that $\mu\sim\lambda$. We then arrive at the following example of a subcategory of $\cat{D}^H$ closed under tensoring with objects in $\cat{C}$.
\begin{example}\label{ex:WeylOrbit}
    For each $\Lambda\in\mathfrak{h}^\ast/\sim$, the full subcategory $\cat{D}^H_\Lambda$ of $\cat{D}^H$ consisting of modules with weight support in $\Lambda$ is closed under tensoring with objects in $\cat{C}$ and thus a braided $\cat{C}$-module category.
\end{example}
We call a weight $\lambda\in\mathfrak{h}^\ast$ \emph{$q$-generic} if  $q^{2(\lambda,\alpha)-m(\alpha,\alpha)} \neq 1$  for all positive roots $\alpha$ and $m\in\Z$. Given $p\in P$, it follows from the fact that $\frac{2(p,\alpha)}{(\alpha,\alpha)}\in\Z$ for all positive roots $\alpha$ that if $\lambda$ is $q$-generic then so is $\lambda+p$ for all $p\in P$. We then call $\Lambda\in\mathfrak{h}^\ast/\sim$ $q$-generic if its representatives $\lambda\in\Lambda$ are $q$-generic. Given $\lambda\in\mathfrak{h}^\ast$ one can define a Verma module $V^H(\lambda)$ for the unrolled quantum group analogously to the classical case, see \cite[Sec.\,5.5]{GP10} for a definition. 
\begin{lemma}\label{lem: q-generic implies irreducible}
    Let $q\in\C^\times$ not be a root of unity. A Verma module $V^H(\lambda)$ is irreducible if $\lambda\in\mathfrak{h}^\ast$ is $q$-generic.
\end{lemma}
\begin{proof}
    By forgetting the action of $U\mathfrak{h}$, a Verma module $V^H(\lambda)$ defines a Verma module $V(\lambda)$ for $U_q\mathfrak{g}$ in the sense of \cite[Sec.\,1.9]{DK90}. 
It follows from \cite[Cor.\,1.9]{DK90} that if $\lambda$ is $q$-generic, then $V^H(\lambda)$ is irreducible as a $U_q\mathfrak{g}$-module, hence also as a $U_q^H\mathfrak{g}$-module.
\end{proof}

Using Lemma~\ref{lem: q-generic implies irreducible}, we define module subcategories of $\cat{D}^H_\Lambda$ which contain only finite direct sums of irreducible Verma modules.  First we show how the tensor product of an irreducible Verma module with a finite-dimensional module decomposes as a direct sum of irreducible Verma modules.

\begin{proposition}\label{prop: Verma tensor finite}
    Let $\lambda\in\mathfrak{h}^\ast$ be $q$-generic and $V^H(\lambda)$ be a Verma module for $U_q^H\mathfrak{g}$. The tensor product $\iota W\otimes V^H(\lambda)$ decomposes as
    \begin{equation}\label{eq:VermaDecomposition}
        \iota W\otimes V^H(\lambda)\cong\bigoplus_{\mu\in\mathfrak{h}^\ast} V^H(\lambda+\mu)^{\bigoplus\dim(W_\mu)}
    \end{equation}
    for all modules $W\in\cat{C}$, where $\iota W \cong \bigoplus_{\mu\in\mathfrak{h}^\ast} W_\mu$ is the weight space decomposition of $\iota W$.
\end{proposition}
\begin{proof}
    Given $\lambda\in\mathfrak{h}^\ast$ we denote with $v_\lambda\in V^H(\lambda)$ the highest weight vector.
    Because $U\mathfrak{h}$ acts as a scalar on each weight space and $V^H(\lambda)$ is irreducible, the result \cite[Prop.\,2.2]{ES99} applies, meaning that for $\lambda\in\mathfrak{h}^\ast$ $q$-generic, there is for every $w\in W_\mu$ a unique intertwiner $\Phi_{\lambda+\mu}^w:V^H(\lambda+\mu)\to \iota W \otimes V^H(\lambda)$ such that
    \begin{equation}\label{eq:VermaDecompA}
        \Phi_{\lambda+\mu}^w(v_{\lambda+\mu}) = w\otimes v_\lambda + u \;,
    \end{equation}
    where $u\in \bigoplus_{\nu<\lambda}(\iota W)\otimes (V^H(\lambda)_\nu)$.
    Combining these, we obtain an intertwiner 
    \begin{equation}\label{eq:VermaDecompB}
        \eta_W^\lambda:\bigoplus_{\mu\in\mathfrak{h}^\ast} W_\mu\otimes V^H(\lambda+\mu)  \to \iota W\otimes V^H(\lambda) 
        ~~,\quad
        w\otimes v \mapsto \Phi^w_{\lambda+\mu}(v)
        \;,
    \end{equation}
where the action of $U_q^H\mathfrak{g}$ on $W_\mu\otimes V^H(\lambda+\mu)$ is only on the second factor.
As a $U_q^H\mathfrak{g}$-module, the domain of $\eta_W^\lambda$ is isomorphic to 
$\bigoplus_{\mu\in\mathfrak{h}^\ast} V^H(\lambda+\mu)^{\bigoplus\dim(W_\mu)}$. To establish \eqref{eq:VermaDecomposition}, it thus remains to show that $\eta_W^\lambda$ is an isomorphism. 
    
We begin with injectivity of $\eta^\lambda_W$. 
First note that by construction the linear map     
    \begin{equation}
        \iota W = \bigoplus_{\mu\in\mathfrak{h}^\ast} W_\mu \to \iota W\otimes V^H(\lambda) \;, \quad w \mapsto \Phi_{\lambda+\mu}^w(v_{\lambda+\mu})
        \quad \text{(for $w \in W_\mu$)}
    \end{equation}
is injective. Because $W\in\cat{C}$, the weight support of $\iota W$ is contained in the weight lattice $P$, hence each $V^H(\lambda+\mu)$ is irreducible. It then follows that if the kernel of $\eta^\lambda_W$ is non-zero, it must also contain an element of the form $w \otimes v_{\lambda+\mu}$, which by the above observation cannot happen. 
    
    Now we show $\eta_W^\lambda$ is surjective. The restriction of $\eta_W^\lambda$ to each weight space is also an injective map. Because the weight spaces of the left- and right-hand side of \eqref{eq:VermaDecomposition} are finite-dimensional, if they have equal dimension, $\eta_W^\lambda$ is surjective. 
    This in turn follows from
    \begin{align}
        \dim(\iota W\otimes V^H(\lambda))_\rho &= \sum_\mu\dim(W_\mu)\dim(V^H(\lambda)_{\rho-\mu})
        \nonumber\\
        &= \sum_\mu \dim(W_\mu)\dim(V^H(\lambda+\mu)_\rho) \nonumber\\
        &= \dim\Bigl(\bigoplus_\mu W_\mu\otimes V^H(\lambda+\mu)\Bigr)_\rho \;.
    \end{align}
\end{proof}

Now we can define an example containing only finite direct sums of irreducible Verma modules.
\begin{example}\label{ex:Verma}
    For each $q$-generic $\Lambda\in\mathfrak{h}^\ast/\sim$ the full subcategory $\cat{M}^H_\Lambda$ of $\cat{D}^H_\Lambda$ of modules that are isomorphic to a finite direct sum of Verma modules of highest weight $\lambda\in\Lambda$ is closed under tensoring with $\cat{C}$, and thus a braided $\cat{C}$-module category. 
\end{example}

\begin{remark}
In the case $\mathfrak{g} = \mathfrak{sl}_2$, facets of the categories $\cat{M}^H_\Lambda$ have appeared in the literature before. Fix the pivotal structure for $\cat{C}$ which recovers the Kauffman bracket skein relations (see \cite{Ti10} for more details) and let $V$ be the fundamental representation of $U_q\mathfrak{sl}_2$. In this setting the Temperley--Lieb algebra on $n$ strands is isomorphic to the endomorphism space of $V^{\otimes n}\in\cat{C}$. Furthermore, the hom-space from $m$ to $n$ in the Temperley--Lieb category is isomorphic to $\cat{C}(V^{\otimes m},V^{\otimes n})$. From the results of \cite{ILZ18}, the blob algebra of \cite{MS93} on $n$ strands is isomorphic to the endomorphism algebra of $V^{\otimes n}\triangleright V^H(\lambda)$ in $\cat{M}^H_\Lambda$. In fact, the hom-spaces of the blob category defined in \cite[Sec.\,2.2]{LNV20} are given by hom-spaces of $\cat{M}^H_\Lambda$. The blob category can be viewed as a module category over the Temperley--Lieb category, with the action given by concatenation. Furthermore, from the identification of the hom-spaces in the blob category with those of $\cat{M}_\Lambda^H$ we observe that the blob category is a braided module category over the Temperley--Lieb category. See also \cite{CGS22} for the relation between endomorphism spaces of $\cat{M}^H_\Lambda$, the blob algebra, and boundary conditions of the XXZ spin chain.
\end{remark}

\subsection{Quantum moment maps from internal-homs}\label{ssec:QMMs from internal-homs}

We now construct examples of quantum moment maps. In order to do this we change settings and work with the $\hbar$-adic Hopf algebra $U_\hbar\mathfrak{g}$ over $\C\llb\hbar\rrb$, see \cite[Sec.\,6.1.3]{KS97} for details regarding the definition. First we introduce some terminology which will let us define categories of $U_\hbar\mathfrak{g}$-modules from which we can obtain quantum moment maps via small braidings. 
Let $W$ be a complex vector space. Define $W\llb\hbar\rrb$ to be the $\C\llb\hbar\rrb$-module containing all formal series $\sum_{n\geq0}w_n\hbar^n$, where $(w_0,w_1,\dots)$ is a sequence of elements in $W$ (see \cite[Eq.\,(XVI.1.2)]{Ka95} for the $\C\llb\hbar\rrb$-module structure). A $\C\llb\hbar\rrb$-module $V$ is called \emph{topologically free} if it is isomorphic as a $\C\llb\hbar\rrb$-module to $W\llb\hbar\rrb$, where $W$ is a complex vector space \cite[Sec.\,XVI.2]{Ka95}. A topologically free module $V\cong W\llb\hbar\rrb$ is said to be of \emph{finite rank} if $W$ is a finite-dimensional vector space \cite[Sec.\,XVII.3]{Ka95}. The \emph{completed tensor product} of $\C\llb\hbar\rrb$-modules, $\widehat{\otimes}$, is the $\hbar$-adic completion of $\otimes_{\C\llb\hbar\rrb}$, see \cite[Def.\,XVI.3.1]{Ka95}. For topologically free modules $X\llb\hbar\rrb$ and $Y\llb\hbar\rrb$ one has $X\llb\hbar\rrb\widehat{\otimes}Y\llb\hbar\rrb \cong (X\otimes Y)\llb\hbar\rrb$ as $\C\llb\hbar\rrb$-modules, so topologically free modules form a monoidal subcategory of $\C\llb\hbar\rrb$-modules \cite[Prop.\,XVI.3.2]{Ka95}. The Hopf algebra $U_\hbar\mathfrak{g}$ is itself topologically free, and isomorphic to $U\mathfrak{g}\llb\hbar\rrb$ as an algebra \cite[Thm.\,XVIII.4.1]{Ka95}. Hence given a topologically free $U_\hbar\mathfrak{g}$-module $W\llb\hbar\rrb$, the vector space $W$ has the canonical structure of a $U\mathfrak{g}$-module.

\begin{remark}
    The category of topologically free $\C\llb\hbar\rrb$-modules, and hence topologically free $U_\hbar\mathfrak{g}$-modules, is not abelian. Consider the endomorphism of $\C\llb\hbar\rrb$ given by $f: r\mapsto r\hbar$. $f$ is injective, hence a monomorphism, and it follows from $\C\llb\hbar\rrb$-linearity and the fact that the target of $f$ is $\hbar$-torsion free that $f$ is an epimorphism. However, $f$ is not surjective and thus not an isomorphism, which implies topologically free modules do not form an abelian category. One workaround for this is to consider the category of $L_0^{\hbar\C\llb\hbar\rrb}$-complete $\C\llb\hbar\rrb$-modules of \cite{PW20}, which is abelian and contains topologically free modules as a monoidal subcategory \cite[Prop.\,3.9]{PW20}. We do not do this in this manuscript, in favour of imposing stronger constraints on our module categories which are sufficient to prove the desired results in our non-abelian setting.
\end{remark}

The Hopf algebra $U_\hbar\mathfrak{g}$ is a quasitriangular $\hbar$-adic Hopf algebra \cite[Thm.\,8.17]{KS97}, meaning the $R$-matrix $\mathcal{R} = \mathcal{R}_1\widehat{\otimes}\mathcal{R}_2$ for $U_\hbar\mathfrak{g}$ is an element of $U_\hbar\mathfrak{g}\widehat{\otimes}U_\hbar\mathfrak{g}$. Hence acting with $\mathcal{R}$ is well-defined on any tensor product of topologically free $U_\hbar\mathfrak{g}$-modules, not just the finite rank ones. We then have that the category of topologically free $U_\hbar\mathfrak{g}$-modules is a braided monoidal category. From \cite[Exercise\,XVI.2]{Ka95} we also have that the category of topologically free $U_\hbar\mathfrak{g}$-modules is a closed monoidal category, with the internal-hom objects $\underline{\cat{D}}(U,V)$ given by $\Hom_{\C\llb\hbar\rrb}(U,V)$. The action of $U_\hbar\mathfrak{g}$ on $\underline{\cat{D}}(U,V)$ is
\begin{equation}
    (x\triangleright f)(\blank) = x_{(1)}f(S(x_{(2)})\blank)
\end{equation}
for all $x\in U_\hbar\mathfrak{g}$ and $f\in \underline{\cat{D}}(U,V)$. Note that this implies the completed tensor product for $\cat{D}$ has a right adjoint, hence preserves colimits.

Now we fix our setting: Let $\cat{C}$ be the category of finite rank topologically free $U_\hbar\mathfrak{g}$-modules and $\cat{D}$ be the category of all topologically free $U_\hbar\mathfrak{g}$-modules. When we write the tensor product $\otimes$ we always mean $\widehat{\otimes}$. The small braiding on $\cat{C}\hookrightarrow\cat{D}$ is simply given by the braiding on $\cat{D}$. We leave the inclusion functor $\cat{C}\hookrightarrow\cat{D}$ implicit, since it is just the identification of finite rank topologically free modules within all topologically free modules. We indicate successive applications of $\mathcal{R}$ with a prime superscript. For example, the monodromy $\rc\circ\lc$ is given by
\begin{equation}\label{eq:SmallHexHopf}
   \rc_{V,U}\circ\lc_{U,V}(u\otimes v) = \mathcal{R}_2^\prime\mathcal{R}_1u\otimes \mathcal{R}_1^\prime\mathcal{R}_2v \;,
\end{equation}
for all $U\in\cat{C}$, $V\in\cat{D}$ and $u\in U$, $v\in V$.

We frequently use the fact that $\cat{D}$ is a closed category, and $\cat{C}$ is rigid \cite[Sec.\,XVII.3]{Ka95}, which allows us to apply the results of Appendix~\ref{app:closed}. Denote by $\underline{\cat{D}}(V)$ the internal endomorphism algebra $(\underline{\cat{D}}(V,V),m_{V,V}^V,u_V)$ of Corollary~\ref{cor:AlgFromClosed}.

\begin{proposition}\label{prop:HalfBraidingClosed}
    Let $\Phi$ be the natural transformation \eqref{eq:Phi}, and denote $\Phi^V_{\blank,\bblank}$ by $\Phi_{\blank,\bblank}$. For all $V\in\cat{D}$, the natural isomorphism $\gamma^{\underline{\cat{D}}(V)}_\blank:\blank\otimes\underline{\cat{D}}(V) \Rightarrow \underline{\cat{D}}(V)\otimes\blank$ 
    defined via
    \begin{equation}\label{eq:InHomHB}
        \begin{tikzcd}
            U\otimes\underline{\cat{D}}(V) && \underline{\cat{D}}(V,U\otimes V) \\
            \underline{\cat{D}}(V) \otimes U & U\otimes \underline{\cat{D}}(V) & \underline{\cat{D}}(V,U\otimes V)
            \arrow["{\Phi_{U,V}}", from=1-1, to=1-3]\arrow["{\underline{\cat{D}}(V,\rc_{V,U}\lc_{U,V})}", from=1-3, to=2-3]\arrow["{\Phi_{U,V}}", from=2-2, to=2-3]\arrow["{\rc_{\underline{\cat{D}}(V),U}}", from=2-1, to=2-2]
            \arrow[dashed, "{\gamma^{\underline{\cat{D}}(V)}_U}", from=1-1, to=2-1]
        \end{tikzcd}
    \end{equation}
    for all $U\in\cat{C}$ is a half-braided structure for $\underline{\cat{D}}(V)$. Let $\cat{M}$ be a full subcategory of $\cat{D}$ which is closed under tensoring with objects in $\cat{C}$. Given $M\in\cat{M}$, consider $\cat{M}$ and $\cat{D}_{\underline{\cat{D}}(M)}$ as braided module categories, with module braiding coming from  the small braiding and half-braided structure respectively. The functor $\underline{\cat{D}}(M,\blank):\cat{M}\to\cat{D}_{\underline{\cat{D}}(M)}$ is a braided $\cat{C}$-module functor.
\end{proposition}
\begin{proof}
That $\gamma^{\underline{\cat{D}}(V)}_\blank$ is an isomorphism follows from invertibility of the small braidings and $\Phi_{U,V}$ for $U\in\cat{C}$ (see Lemma~\ref{lem:PhisoIffRigid}). Naturality of $\gamma^{\underline{\cat{D}}(V)}_\blank$ follows from naturality of the coherence morphisms and the small braidings. It remains to show that the identities \eqref{eq:hHBA} and \eqref{eq:small-half-hexagon} are satisfied.
    
Let $X,Y\in\cat{C}$. Conjugating the image of a morphism $f:X\otimes V\to Y\otimes V$ under $\underline{\cat{D}}(V,\blank)$ with $\Phi$ defines a morphism
    \begin{equation}
        \underline{f} := \Phi_{Y,V}^{-1} \underline{\cat{D}}(V,f) \Phi_{X,V}:X\otimes\underline{\cat{D}}(V)\to Y\otimes\underline{\cat{D}}(V) \;.
    \end{equation}
    We depict this as follows.
    \begin{equation}
        \underline{f} = 
        \begin{tikzpicture}[baseline={(0,0.8)}]
            \node at (0,-0.2) {$X$}; \node[green] at (0.7,-0.2) {$\underline{\cat{D}}(V)$};
            \draw (0,0) -- (0,0.5); \draw[thick green] (0.5,0) -- (0.5,0.5);
            \draw[-Latex] (0,1.5) -- (0,2.25); \draw[-Latex,thick green] (0.5,1.5) -- (0.5,2.25);
            \node at (0.25,1) {$f$};
            \draw (-0.2,0.75) -- (0.7,0.75) -- (0.7,1.25) -- (-0.2,1.25) -- (-0.2,0.75);
            \draw (0,0.5) -- (0,0.75); \draw[thick black] (0.5,0.5) -- (0.5,0.75);
            \draw (0,1.25) -- (0,1.5); \draw[thick black] (0.5,1.25) -- (0.5,1.5);
            \draw[red] (-0.2,0.5) -- (0.7,0.5);
            \draw[red] (-0.2,1.5) -- (0.7,1.5);
            \node at (0,2.45) {$Y$}; \node[green] at (0.7,2.45) {$\underline{\cat{D}}(V)$};
        \end{tikzpicture}
    \end{equation}
    Note that $U\otimes \underline{f} = \underline{U\otimes f}$, which is depicted
    \begin{equation}
        \begin{tikzpicture}[baseline={(0,0.8)}]
            \node at (0,-0.2) {$X$}; \node[green] at (0.7,-0.2) {$\underline{\cat{D}}(V)$};
            \draw (0,0) -- (0,0.5); \draw[thick green] (0.5,0) -- (0.5,0.5);
            \draw[-Latex] (0,1.5) -- (0,2.25); \draw[-Latex,thick green] (0.5,1.5) -- (0.5,2.25);
            \node at (0.25,1) {$f$};
            \draw (-0.2,0.75) -- (0.7,0.75) -- (0.7,1.25) -- (-0.2,1.25) -- (-0.2,0.75);
            \draw (0,0.5) -- (0,0.75); \draw[thick black] (0.5,0.5) -- (0.5,0.75);
            \draw (0,1.25) -- (0,1.5); \draw[thick black] (0.5,1.25) -- (0.5,1.5);
            \draw[red] (-0.2,0.5) -- (0.7,0.5);
            \draw[red] (-0.2,1.5) -- (0.7,1.5);
            \node at (0,2.45) {$Y$}; \node[green] at (0.7,2.45) {$\underline{\cat{D}}(V)$};
            \node at (-0.5,-0.2) {$U$}; \draw[-Latex] (-0.5,0) -- (-0.5,2.25); \node at (-0.5,2.45) {$U$};
        \end{tikzpicture}
        =
        \begin{tikzpicture}[baseline={(0,0.8)}]
            \node at (0,-0.2) {$X$}; \node[green] at (0.7,-0.2) {$\underline{\cat{D}}(V)$};
            \draw (0,0) -- (0,0.5); \draw[thick green] (0.5,0) -- (0.5,0.5);
            \draw[-Latex] (0,1.5) -- (0,2.25); \draw[-Latex,thick green] (0.5,1.5) -- (0.5,2.25);
            \node at (0.25,1) {$f$};
            \draw (-0.2,0.75) -- (0.7,0.75) -- (0.7,1.25) -- (-0.2,1.25) -- (-0.2,0.75);
            \draw (0,0.5) -- (0,0.75); \draw[thick black] (0.5,0.5) -- (0.5,0.75);
            \draw (0,1.25) -- (0,1.5); \draw[thick black] (0.5,1.25) -- (0.5,1.5);
            \draw[red] (-0.7,0.5) -- (0.7,0.5);
            \draw[red] (-0.7,1.5) -- (0.7,1.5);
            \node at (0,2.45) {$Y$}; \node[green] at (0.7,2.45) {$\underline{\cat{D}}(V)$};
            \node at (-0.5,-0.2) {$U$}; \draw[-Latex] (-0.5,0) -- (-0.5,2.25); \node at (-0.5,2.45) {$U$};
        \end{tikzpicture}
    \end{equation}
    Using this and the hexagons \eqref{eq:SmallHexA} and \eqref{eq:SmallHexD}, it is a straightforward computation to show that $\gamma^{\underline{\cat{D}}(V)}_\blank$ satisfies the hexagon \eqref{eq:small-half-hexagon}.
    
    The candidate \eqref{eq:InHomHB} for the half-braided structure is then depicted
    \begin{equation}
        \gamma_U^{\underline{\cat{D}}(V)} = 
        \begin{tikzpicture}[baseline={(0,0.8)}]
            \node at (0,-0.2) {$U$}; \node[green] at (0.65,-0.2) {$\underline{\cat{D}}(V)$};
            \draw (0,0) -- (0,0.5); \draw[thick green] (0.5,0) -- (0.5,0.5);
            \draw[-Latex] (0.5,2) -- (0.5,2.75); \draw[-Latex,thick green] (0,2) -- (0,2.75);
            \pic[rightstr = thick black] at (0,0.5) {braid}; \pic[leftstr = thick black] at (0,1) {braid}; \pic[rightstr = thick green] at (0,1.5) {invbraid};
            \draw[red] (-0.2,0.5) -- (0.7,0.5);
            \draw[red] (-0.2,1.5) -- (0.7,1.5);
            \node at (0.5,2.95) {$U$}; \node[green] at (-0.2,2.95) {$\underline{\cat{D}}(V)$};
        \end{tikzpicture}
    \end{equation}
    The first equality in \eqref{eq:hHBA} follows from the naturality of the small braidings and the fact that $\underline{\cat{D}}(V,\blank):\cat{D}\to\cat{D}_{\underline{\cat{D}}(V)}$ is a $\cat{C}$-module functor, see Corollary~\ref{cor:C-module-structure-from-rigid}, which implies $\underline{f}$ is a morphism in $\cat{D}_{\underline{\cat{D}}(V)}$:
    \begin{equation}
        \begin{tikzpicture}[baseline={(0,0.8)}]
            \centerarc[thick green](1,0.5)(0:180:0.5); \draw[thick green] (0.5,0) -- (0.5,0.5); \draw[thick green] (1.5,0) -- (1.5,0.5); \node[green] at (0.7,-0.2) {$\underline{\cat{D}}(V)$}; \node[green] at (1.7,-0.2) {$\underline{\cat{D}}(V)$};
            \draw (0,0) -- (0,1.5); \node at (0,-0.2) {$U$};
            \centerarc[thick green](0.75,1)(0:90:0.25); \centerarc[thick green](0.75,1.5)(180:270:0.25);
            \draw[-Latex] (0.5,3) -- (0.5,3.5); \node at (0.5,3.7) {$U$};
            \draw[thick green, -Latex] (0,3) -- (0,3.5); \node[green] at (-0.2,3.7) {$\underline{\cat{D}}(V)$};
            \pic[rightstr = thick black] at (0,1.5) {braid}; \pic[leftstr = thick black] at (0,2) {braid}; \pic[rightstr = thick green] at (0,2.5) {invbraid};
            \draw[red] (-0.2,1.5) -- (0.7,1.5);
            \draw[red] (-0.2,2.5) -- (0.7,2.5);
        \end{tikzpicture}
        =
        \begin{tikzpicture}[baseline={(0,0.8)}]
            \node at (0,-0.2) {$U$}; \node[green] at (0.7,-0.2) {$\underline{\cat{D}}(V)$};
            \draw (0,0) -- (0,0.5); \draw[thick green] (0.5,0) -- (0.5,0.5);
            \draw[-Latex] (0.5,2.5) -- (0.5,3); \draw[-Latex,thick green] (0,2.5) -- (0,3);
            \pic[rightstr = thick black] at (0,0.5) {braid}; \pic[leftstr = thick black] at (0,1) {braid}; \pic[rightstr = thick green] at (0,2) {invbraid};
            \draw[thick green] (0.5,1.5) -- (0.5,2); \draw (0,1.5) -- (0,2); \centerarc[thick green](0.5,1.25)(0:90:0.75);
            \draw[thick green] (1.25,0) -- (1.25,1.25); \node[green] at (1.65,-0.2) {$\underline{\cat{D}}(V)$};
            \draw[red] (-0.2,0.5) -- (0.7,0.5);
            \draw[red] (-0.2,1.5) -- (0.7,1.5);
            \node at (0.5,3.2) {$U$}; \node[green] at (-0.2,3.2) {$\underline{\cat{D}}(V)$};
        \end{tikzpicture}
        =
        \begin{tikzpicture}[baseline={(0,0.8)}]
            \node at (0,-0.2) {$U$}; \node[green] at (0.7,-0.2) {$\underline{\cat{D}}(V)$};
            \draw (0,0) -- (0,0.5); \draw[thick green] (0.5,0) -- (0.5,0.5);
            \draw[-Latex] (1,2.5) -- (1,3.5); \draw[-Latex,thick green] (0.25,2.75) -- (0.25,3.5);
            \centerarc[thick green](0.25,2.5)(0:180:0.25);
            \pic[rightstr = thick black] at (0,0.5) {braid}; \pic[leftstr = thick black] at (0,1) {braid};
            \pic[rightstrA = thick green, rightstrB = thick green] at (0,1.5) {braiduts};
            \draw[thick green,smooth] plot coordinates {(1.25,0) (1.25,1) (1.05,1.3) (1,1.5)}; \node[green] at (1.65,-0.2) {$\underline{\cat{D}}(V)$};
            \draw[red] (-0.2,0.5) -- (0.7,0.5);
            \draw[red] (-0.2,1.5) -- (0.7,1.5);
            \node at (1,3.7) {$U$}; \node[green] at (0.15,3.7) {$\underline{\cat{D}}(V)$};
        \end{tikzpicture}
    \end{equation}
    The second equality in \eqref{eq:hHBA} is depicted
    \begin{equation}\label{eq:ClosedHBPentagon2}
        \begin{tikzpicture}[baseline={(0,0.8)}]
            \node at (-0.6,-0.2) {$U$}; \node[green] at (0,-0.2) {$\scriptstyle \underline{\cat{D}}(V)$}; \node[green] at (0.8,-0.2) {$\scriptstyle \underline{\cat{D}}(V)$};
            \draw (0,0.5) -- (0,1); \draw[thick green] (-0.5,0.5) -- (-0.5,2.5); \pic[rightstr = thick green] at (-0.5,0) {braid};
            \pic[rightstr = thick black] at (0,1) {braid}; \pic[leftstr = thick black] at (0,1.5) {braid}; \pic[rightstr = thick green] at (0,2) {invbraid};
            \draw[thick green] (0.5,0) -- (0.5,1);
            \draw[red] (-0.2,1) -- (0.7,1);
            \draw[red] (-0.2,2) -- (0.7,2);
            \centerarc[thick green](-0.25,2.5)(0:180:0.25); \draw[thick green,-Latex] (-0.25,2.75) -- (-0.25,3.5);
            \draw[-Latex] (0.5,2.5) -- (0.5,3.5);
            \node at (0.5,3.7) {$U$}; \node[green] at (-0.5,3.7) {$\underline{\cat{D}}(V)$};
        \end{tikzpicture}
        =
        \begin{tikzpicture}[baseline={(0,0.8)}]
            \centerarc[thick green](1,0.5)(0:180:0.5); \draw[thick green] (0.5,0) -- (0.5,0.5); \draw[thick green] (1.5,0) -- (1.5,0.5); \node[green] at (0.7,-0.2) {$\underline{\cat{D}}(V)$}; \node[green] at (1.7,-0.2) {$\underline{\cat{D}}(V)$};
            \draw (0,0) -- (0,1.5); \node at (0,-0.2) {$U$};
            \centerarc[thick green](0.75,1)(0:90:0.25); \centerarc[thick green](0.75,1.5)(180:270:0.25);
            \draw[-Latex] (0.5,3) -- (0.5,3.5); \node at (0.5,3.7) {$U$};
            \draw[thick green, -Latex] (0,3) -- (0,3.5); \node[green] at (-0.2,3.7) {$\underline{\cat{D}}(V)$};
            \pic[rightstr = thick black] at (0,1.5) {braid}; \pic[leftstr = thick black] at (0,2) {braid}; \pic[rightstr = thick green] at (0,2.5) {invbraid};
            \draw[red] (-0.2,1.5) -- (0.7,1.5);
            \draw[red] (-0.2,2.5) -- (0.7,2.5);
        \end{tikzpicture} \;.
    \end{equation}
    We show that this holds via a direct calculation.

    First it is useful to express the half-braiding in components. Denote by $\Delta$, $\epsilon$ and $S$ the coproduct, counit and antipode for $U_\hbar\mathfrak{g}$, respectively. Let $\mathcal{R}^{-1} = \mathcal{S}_1\otimes\mathcal{S}_2$. For all $U\in\cat{C}$, $u\in U$ and $a\in\underline{\cat{D}}(V)$, we have
    \begin{align}\label{eq:ClosedHBPentagon3}
        \gamma_U^{\underline{\cat{D}}(V)}(u\otimes a) &= [ v\mapsto \mathcal{S}_1\triangleright(\mathcal{R}_1^\prime\mathcal{R}_2a)(v)]\otimes \mathcal{S}_2\mathcal{R}_2^\prime\mathcal{R}_1u \nonumber\\
        &= [ v\mapsto \mathcal{S}_{1(1)}\mathcal{R}_1^\prime\mathcal{R}_2a(S(\mathcal{S}_{1(2)})v)]\otimes \mathcal{S}_2\mathcal{R}_2^\prime\mathcal{R}_1u \;,
    \end{align}
    where $xa$ for $x\in U_\hbar\mathfrak{g}$ and $a\in\underline{\cat{D}}(V)$ denotes post-composition (that is, $xa = (v\mapsto xa(v))$ for all $v\in V$). Using the identity 
    \begin{equation}
        (\Delta\otimes\Id)(\mathcal{R}^{-1})(\mathcal{R}_1\otimes1\otimes\mathcal{R}_2) = (1\otimes \mathcal{R}^{-1})(\mathcal{S}_1\mathcal{R}_1\otimes1\otimes\mathcal{S}_2\mathcal{R}_2) = 1\otimes\mathcal{R}^{-1} \;,
    \end{equation}
    we write
    \begin{align}
        &\mathcal{S}_{1(1)}\mathcal{R}_1^\prime\mathcal{R}_2 a \otimes S(\mathcal{S}_{1(2)})v \otimes \mathcal{S}_2\mathcal{R}_2^\prime\mathcal{R}_1 u 
        \nonumber \\
        &= (\Id\otimes S\otimes\Id)(\Delta\otimes\Id)(\mathcal{R}^{-1})(\mathcal{R}_1^\prime\mathcal{R}_2 a \otimes v \otimes \mathcal{R}_2^\prime\mathcal{R}_1 u) \nonumber\\
        &= \mathcal{R}_2 a \otimes S(\mathcal{S}_1)v \otimes \mathcal{S}_2\mathcal{R}_1 u \;.
    \end{align}
    Hence \eqref{eq:ClosedHBPentagon3} simplifies to
    \begin{equation}
        \gamma_U^{\underline{\cat{D}}(V)}(u\otimes a) = [v\mapsto \mathcal{R}_2 a(S(\mathcal{S}_1)v)]\otimes \mathcal{S}_2\mathcal{R}_1 u \;.
    \end{equation}
    Evaluating the left-hand side of \eqref{eq:ClosedHBPentagon2} on any $u\otimes a\otimes b\in U\otimes\underline{\cat{D}}(V)\otimes\underline{\cat{D}}(V)$ results in
    \begin{align}\label{eq:ClosedHBPentagon4}
        u\otimes a\otimes b &\xmapsto{\lc_{U,\underline{\cat{D}}(V)}} \mathcal{R}_2\triangleright a\otimes \mathcal{R}_1 u \otimes b \nonumber\\
        &\xmapsto{\underline{\cat{D}}(V)\otimes\gamma_U^{\underline{\cat{D}}(V)}}  \mathcal{R}_2\triangleright a\otimes [v\mapsto\mathcal{R}_2^\prime b(S(\mathcal{S}_1)v)]\otimes \mathcal{S}_2\mathcal{R}_1^\prime\mathcal{R}_1 u \nonumber\\
        &\xmapsto{m\otimes U} [v\mapsto (\mathcal{R}_2\triangleright a)(\mathcal{R}_2^\prime b(S(\mathcal{S}_1)v))]\otimes \mathcal{S}_2\mathcal{R}_1^\prime\mathcal{R}_1 u \;.
    \end{align}
    Using $(\Id\otimes\Delta)(\mathcal{R}) = \mathcal{R}_1^\prime\mathcal{R}_1\otimes\mathcal{R}_2\otimes\mathcal{R}_2^\prime$, the right-hand side of \eqref{eq:ClosedHBPentagon4} becomes
    \begin{equation}
        \begin{aligned}
        &[v \mapsto(\mathcal{R}_{2(1)} \triangleright a)(\mathcal{R}_{2(2)} b(S(\mathcal{S}_1)v))]\otimes \mathcal{S}_2\mathcal{R}_1 u
        \\&\quad =[v \mapsto\mathcal{R}_{2(1)}a(S(\mathcal{R}_{2(2)})\mathcal{R}_{2(3)}b(S(\mathcal{S}_1)v))]\otimes \mathcal{S}_2\mathcal{R}_1 u \;.
        \end{aligned}
    \end{equation}
    Because $S(\mathcal{R}_{2(2)})\mathcal{R}_{2(3)} = \epsilon(\mathcal{R}_{2(2)})$, this becomes
    \begin{equation}
        [v\mapsto \mathcal{R}_2 (a\circ b)(S(\mathcal{S}_1)v)]\otimes \mathcal{S}_2\mathcal{R}_1 u \;.
    \end{equation}
    Evaluating the right-hand side of \eqref{eq:ClosedHBPentagon2}, we find
    \begin{equation}
        \gamma^{\underline{\cat{D}}(V)}_U\circ(U\otimes m)(u\otimes a\otimes b) = [v\mapsto \mathcal{R}_2 (a\circ b)(S(\mathcal{S}_1)v)]\otimes \mathcal{S}_2\mathcal{R}_1 u \;.
    \end{equation}
    So \eqref{eq:ClosedHBPentagon2} holds, as desired.
    
    That $\underline{\cat{D}}(M,\blank)$ is a $\cat{C}$-module functor follows from Corollary~\ref{cor:C-module-structure-from-rigid}. The proof that $\underline{\cat{D}}(M,\blank)$ is a braided module functor follows from an almost identical computation to the one proving \eqref{eq:ClosedHBPentagon2}. To check that \eqref{eq:BrModFunct} is satisfied all one needs to do is to replace the left $\underline{\cat{D}}(M)$ strand in \eqref{eq:ClosedHBPentagon2} with a $\underline{\cat{D}}(M)$-module, and apply the unit to the right strand.
\end{proof}

\begin{remark}
    Let $H$ be a quasitriangular Hopf algebra over a commutative ring $K$. Proposition~\ref{prop:HalfBraidingClosed} also holds for $\cat{D}$ being the category of $H$-modules and $\cat{C}$ the full subcategory of 
    rigid objects in $\cat{D}$, with the monoidal structure given by the tensor product over $K$. The proof is exactly the same.
\end{remark}

It is well-known that $\cat{C}$ is a ribbon category, see for example \cite[Sec.\,XVII.3]{Ka95}. Furthermore, the coend $\mathcal{L} = \int^{V\in\cat{C}}V^\ast\otimes V$ exists in $\cat{D}$, and is given by the so-called \emph{braided dual} \cite[p.\,26]{Ke23} (see also \cite[Sec.\,6]{BBJ15} for the braided dual in a related setting). Because $\cat{D}$ is closed monoidal, the functors $\blank\otimes V$ are left adjoint to $\underline{\cat{D}}(V,\blank)$ for all $V\in\cat{D}$, hence they preserve colimits and thus coends, namely $\mathcal{L}$. $\cat{D}$ is braided, so $\blank\otimes V\cong V\otimes\blank$ for all $V\in\cat{D}$, which implies $V\otimes\blank$ also preserves $\cat{L}$. Applying Proposition~\ref{Prop:hCQMM} to the half-braiding \eqref{eq:InHomHB}, we arrive at the following source of quantum moment maps coming from $U_\hbar\mathfrak{g}$-modules.
\begin{theorem}\label{thm:QMMHmod}
    Let $\cat{D}$ be the category of topologically free $U_\hbar\mathfrak{g}$-modules, $\cat{C}$ the category of finite rank topologically free $U_\hbar\mathfrak{g}$-modules, and $\underline{\cat{D}}(V)\in\cat{D}$ the internal endomorphism algebra of $V\in\cat{D}$. For all $V\in\cat{D}$, the morphism $\mu:\mathcal{L}\to\underline{\cat{D}}(V)$ defined via
    \begin{equation}
        \begin{tikzcd}
            U^\ast\otimes U && (U^\ast\otimes U)\otimes\underline{\cat{D}}(V) & U^\ast\otimes\underline{\cat{D}}(V,U\otimes V) \\
            \mathcal{L} &&&\\
            \underline{\cat{D}}(V) && (U^\ast\otimes U)\otimes\underline{\cat{D}}(V) & U^\ast\otimes\underline{\cat{D}}(V,U\otimes V)
            \arrow["(U^\ast\otimes U)\otimes u_V", from=1-1, to=1-3]\arrow["\sim", from=1-3, to=1-4]\arrow["{U^\ast\otimes\underline{\cat{D}}(V,\rc_{V,U}\lc_{U,V})}"', from=1-4, to=3-4]\arrow["\sim", from=3-4, to=3-3]\arrow["{\ev_U\otimes\underline{\cat{D}}(V)}", from=3-3, to=3-1]
            \arrow["i_U", from=1-1, to=2-1]\arrow[dashed, "\mu", from=2-1, to=3-1]
        \end{tikzcd}
    \end{equation}
    is a quantum moment map.
\end{theorem}
\begin{remark}\label{rem:Mudrov}
    The quantum moment maps $\mu:\mathcal{L}\to\underline{\cat{D}}(V)$ are related to the quantised conjugacy classes of Lie groups defined in \cite{Mu04}. To see this, we begin by writing the explicit form of $\mu \circ i_U(f\otimes u)$ as an endomorphism of $V$. This is
    \begin{equation}\label{eq:QMMexplicit}
        \mu \circ i_U(f\otimes u) = (v\mapsto f(\mathcal{R}_2^\prime\mathcal{R}_1 u)\mathcal{R}_1^\prime\mathcal{R}_2 v )
    \end{equation}
    for all $v\in V$. The map
    \begin{equation}
        \varphi:\mathcal{L} \to U_\hbar\mathfrak{g} \;,\quad i_U(f\otimes u) \mapsto f(\mathcal{R}_2^\prime\mathcal{R}_1 u)\mathcal{R}_1^\prime\mathcal{R}_2
    \end{equation}
    is precisely \cite[Eq.\,(21)]{Mu04}. Let $\rho_V:U_\hbar\mathfrak{g}\to\mathrm{End}_{\C\llb\hbar\rrb}(V)$ be the $U_\hbar\mathfrak{g}$-module structure of $V$. The map \cite[Eq.\,(23)]{Mu04} is not available in our setting as the module $M_{\mathfrak{p},\lambda/2\hbar}$ in \cite{Mu04} is not in $\cat{D}$, see \cite[Rem.\,3.3]{Mu04}. The composition $\rho_V\circ\varphi$ coincides with $\mu$ defined via \eqref{eq:QMMexplicit}, and closely resembles \cite[Eq.\,(23)]{Mu04}, suggesting that the quantised conjugacy classes of \cite{Mu04} may conceptually be thought of as internal-end algebras in an appropriate category of modules for the $\hbar$-adic quantum group. 
\end{remark}

The examples in Section~\ref{ssec:Examples for generic q} make sense here as well.  Moreover, we can now define quantum moment maps for these examples. Denote by $V_\hbar(\lambda)$ a Verma module for $U_\hbar\mathfrak{g}$ of highest weight $\lambda\in\mathfrak{h}^\ast$, see for example \cite[Sec.\,3.1]{Mu04} for a definition of Verma modules for the $\hbar$-adic quantum group. The irreducibility criterion for Verma modules is different however, namely it is the same as for $U\mathfrak{g}$-modules. Replacing $q$-generic with generic in the sense of \cite[Sec.\,3.2]{Mu04} (which is also invariant under shifts by weights in $P$), the tensor product decomposition of Proposition~\ref{prop: Verma tensor finite} follows from \cite[Prop.\,3.7]{Mu04}. Hence we define the $\hbar$-adic analogue of Example~\ref{ex:Verma}, $\cat{M}_\Lambda$, to be the full subcategory of $\cat{D}$ whose objects are isomorphic to finite direct sums of Verma modules with highest weight $\lambda\in\Lambda$ where $\Lambda\in\mathfrak{h}^\ast/\sim$ is generic. Because of the decomposition in \eqref{eq:VermaDecomposition}, for all $M\in\cat{M}_\Lambda$ and $\lambda\in\Lambda$ there exists $U\in\cat{C}$ such that $M$ is a direct summand (hence retract) of $U\otimes V_\hbar(\lambda)$, which implies $\cat{M}_\Lambda$ satisfies the conditions of $\cat{M}$ in Proposition~\ref{prop:FFMCEmbedding}. Upon applying Proposition~\ref{prop:FFMCEmbedding} and Theorem~\ref{thm:QMMHmod} we obtain inclusions of braided module categories from $\cat{M}_\Lambda$:

\begin{proposition}\label{prop:VermaQMM}
    For all $\lambda\in\Lambda$, consider $\cat{D}_{A_\lambda}$ where $A_\lambda := \underline{\cat{D}}(V_\hbar(\lambda))$, as a braided $\cat{C}$-module category, with the module braiding coming from the quantum moment map of Theorem~\ref{thm:QMMHmod} via the prescription in Proposition~\ref{Prop:hCQMM}. There is a full inclusion of braided module categories $\underline{\cat{D}}(V_\hbar(\lambda),\blank):\cat{M}_\Lambda\hookrightarrow\cat{D}_{A_\lambda}$.
\end{proposition}
\begin{remark}\label{rem:Not an equivalence}
    The functor $\underline{\cat{D}}(V_\hbar(\lambda),\blank)$ is not essentially surjective: For example, the full dual $\underline{\cat{D}}(V_\hbar(\lambda),\C\llb\hbar\rrb)$ of the Verma module $V_\hbar(\lambda)$ is not isomorphic to $\underline{\cat{D}}(V_\hbar(\lambda),N)$ for any $N\in\cat{M}_\Lambda$. This follows from the fact that the full dual has a lowest weight $-\lambda$, whereas $\underline{\cat{D}}(V_\hbar(\lambda),N)$ for non-zero $N$ does not have a lowest weight. Hence $\underline{\cat{D}}(V_\hbar(\lambda),\blank)$ is not an equivalence of braided module categories.
\end{remark}

\appendix

\section{Closed monoidal categories}\label{app:closed}

We recall the definition and some basic facts about closed monoidal categories. We refer the reader to \cite{EK66} for proofs of the general statements.
\begin{definition}[{\cite[Sec.\,II.2 \& Thm.\,II.5.3]{EK66}}]
    A monoidal category $\cat{D}$ is \emph{closed} if the functor $\blank\otimes V:\cat{D}\to\cat{D}$ has a right adjoint $\underline{\cat{D}}(V,\blank)$ for all $V\in\cat{D}$. That is, for all $U,V,W\in\cat{D}$ there is a natural isomorphism $\cat{D}(U\otimes V,W) \cong \cat{D}(U,\underline{\cat{D}}(V,W))$.
\end{definition}
Let $\cat{D}$ be a closed monoidal category in the remainder of the appendix. Given $f\in\cat{D}(U\otimes V,W)$ we denote with $f^\prime\in\cat{D}(U,\underline{\cat{D}}(V,W))$ the image of $f$ under the isomorphism $\cat{D}(U\otimes V,W) \to \cat{D}(U,\underline{\cat{D}}(V,W))$. The internal-hom objects $\underline{\cat{D}}(V,W)$ come equipped with extra structure, which we now describe. 

For each $V\in\cat{D}$, define the canonical natural transformations
\begin{equation}\label{eq:eviso}
    \ev_{V,\blank}:\underline{\cat{D}}(V,\blank)\otimes V \Rightarrow \blank \quad\text{and}\quad \coev_{\blank,V}:\blank\Rightarrow \underline{\cat{D}}(V,\blank\otimes V)
\end{equation}
to be the images of the identity morphisms under the isomorphisms 
\begin{equation}\label{eq:coeviso}
    \cat{D}(\underline{\cat{D}}(V,\blank),\underline{\cat{D}}(V,\blank))\Rightarrow\cat{D}(\underline{\cat{D}}(V,\blank)\otimes V,\blank) \;,
\end{equation}
and
\begin{equation}
    \cat{D}(\blank\otimes V,\blank\otimes V) \Rightarrow \cat{D}(\blank,\underline{\cat{D}}(V,\blank\otimes V)) \;,
\end{equation}
respectively.

The proofs of the two lemmas which follow can be found in \cite[Ch.\,II]{EK66}.
\begin{lemma}[{\cite[Eqs.\,(II.3.4) \& (II.3.5)]{EK66}}] \label{lemma:ClosedStructureA}
    The natural transformations $\ev_{V,\blank}$ and $\coev_{\blank,V}$ are such that the diagrams
    \begin{equation}
        \begin{tikzcd}
            U & \underline{\cat{D}}(V,W) \\
             & \underline{\cat{D}}(V,U\otimes V)
            \arrow["{f^\prime}", from=1-1, to=1-2]
            \arrow["{\coev_{U,V}}"', from=1-1, to=2-2]\arrow["{\underline{\cat{D}}(V,f)}"', from=2-2, to=1-2]
        \end{tikzcd} \;,\quad \begin{tikzcd}
            U\otimes V & W \\
             & \underline{\cat{D}}(V,W)\otimes V
            \arrow["{f}", from=1-1, to=1-2]
            \arrow["{f^\prime\otimes V}"', from=1-1, to=2-2]\arrow["{\ev_{V,W}}"', from=2-2, to=1-2]
        \end{tikzcd}
    \end{equation}
    commute for all $U,V,W\in\cat{D}$ and $f\in\cat{D}(U\otimes V,W)$.
\end{lemma}
For all $U,V,W\in\cat{D}$ there are canonical morphisms
\begin{equation}
    m_{U,V}^W:\underline{\cat{D}}(U,V)\otimes\underline{\cat{D}}(W,U) \to \underline{\cat{D}}(W,V) \;,
\end{equation}
given by the image of 
\begin{equation}
    \ev_{U,V}\circ(\underline{\cat{D}}(U,V)\otimes\ev_{W,U}): \underline{\cat{D}}(U,V)\otimes(\underline{\cat{D}}(W,U)\otimes W) \to V
\end{equation}
under the isomorphism
\begin{equation}
    \cat{D}(\underline{\cat{D}}(U,V)\otimes(\underline{\cat{D}}(W,U)\otimes W), V) \to \cat{D}(\underline{\cat{D}}(U,V)\otimes\underline{\cat{D}}(W,U), \underline{\cat{D}}(W,V)) \;.
\end{equation}

\begin{lemma}[{\cite[Thm.\,6.4, proof of (iii)]{EK66}}]\label{lemma:ClosedStructureB}
For all $U,V,W,X\in\cat{D}$, the following diagram commutes:
    \begin{equation}
        \begin{tikzcd}[column sep=0.5em]
            & {\scriptstyle \underline{\cat{D}}(U,V)\otimes\underline{\cat{D}}(X,U)} & \\
            {\scriptstyle \underline{\cat{D}}(U,V)\otimes(\underline{\cat{D}}(W,U)\otimes\underline{\cat{D}}(X,W))} && {\scriptstyle \underline{\cat{D}}(X,V) } \\
            {\scriptstyle (\underline{\cat{D}}(U,V)\otimes\underline{\cat{D}}(W,U))\otimes\underline{\cat{D}}(X,W) } && {\scriptstyle \underline{\cat{D}}(W,V)\otimes\underline{\cat{D}}(X,W)}
            \arrow["{\underline{\cat{D}}(U,V)\otimes m_{W,U}^X}", from=2-1, to=1-2]\arrow["{m_{U,V}^X}", from=1-2, to=2-3]
            \arrow["{\alpha}"', from=3-1, to=2-1]\arrow["{m_{U,V}^W\otimes\underline{\cat{D}}(X,W)}"', from=3-1, to=3-3]\arrow["{m_{W,V}^X}"', from=3-3, to=2-3]
        \end{tikzcd}
    \end{equation}
\end{lemma}

Denote by $\underline{\cat{D}}(V)$ the internal endomorphism object $\underline{\cat{D}}(V,V)$. Define the composition
\begin{equation}
    u_V:\1\xrightarrow{\coev_{\1,V}}\underline{\cat{D}}(V,\1\otimes V)\xrightarrow{\underline{\cat{D}}(V,l_V)} \underline{\cat{D}}(V,V)\;.
\end{equation}
Using Lemmas~\ref{lemma:ClosedStructureA} and~\ref{lemma:ClosedStructureB} we construct algebras and bimodules in $\cat{D}$ from the internal-hom objects:
\begin{corollary}\label{cor:AlgFromClosed}
     Each tuple $(\underline{\cat{D}}(V),m_{V,V}^V,u_V)$ is an algebra in $\cat{D}$, which we refer to simply as $\underline{\cat{D}}(V)$. Each tuple $(\underline{\cat{D}}(U,V),m_{V,V}^U,m_{U,V}^U)$ is a $(\underline{\cat{D}}(V),\underline{\cat{D}}(U))$-bimodule in $\cat{D}$.
\end{corollary}

Define the natural transformation $\Phi_{\blank,\bblank}^\bbblank:\blank\otimes\underline{\cat{D}}(\bbblank,\bblank) \to \underline{\cat{D}}(\bbblank,\blank\otimes\bblank)$ via the composition
\begin{equation}\label{eq:Phi}
    \Phi_{V,W}^U = m_{W,V\otimes W}^U\circ(\coev_{V,W}\otimes \underline{\cat{D}}(U,W)): V\otimes\underline{\cat{D}}(U,W) \to \underline{\cat{D}}(U,V\otimes W) \;.
\end{equation}
\begin{lemma}\label{lem:PhisoIffRigid}
    The natural transformation $\Phi_{V,\blank}^\bblank$ is an isomorphism if and only if $V$ has a left dual in $\cat{D}$.
\end{lemma}
\begin{proof}
    Let $U,V,W\in\cat{D}$ and assume $V$ has a left dual. For all $X\in\cat{D}$, we have
    \begin{align}
        \cat{D}(X\otimes U,V\otimes W) &\cong \cat{D}(V^\ast\otimes(X\otimes U),W) \nonumber\\
        &\cong  \cat{D}((V^\ast\otimes X)\otimes U,W) \nonumber\\
        &\cong  \cat{D}(V^\ast\otimes X,\underline{\cat{D}}(U,W)) \nonumber\\
        &\cong  \cat{D}(X,V\otimes\underline{\cat{D}}(U,W)) \;.
    \end{align}
    By definition of the internal-hom, we also have
    \begin{equation}
        \cat{D}(X\otimes U,V\otimes W)\cong\cat{D}(X,\underline{\cat{D}}(U,V\otimes W)) \;.
    \end{equation}
    Hence, by the Yoneda Lemma, we have $V\otimes\underline{\cat{D}}(U,W) \cong \underline{\cat{D}}(U,V\otimes W)$. Furthermore, this isomorphism exhibits $V\otimes\underline{\cat{D}}(U,W)$ as representing $\cat{D}(\blank\otimes U,V\otimes W)$, making it exactly $\Phi_{V,W}^U$.

    Assume $\Phi_{V,\blank}^\bblank$ is an isomorphism. Let $V\in\cat{D}$ and consider the object $\underline{\cat{D}}(V,\1)$. We have evaluation and coevaluation morphisms given by
    \begin{equation}
        \ev_V := \ev_{V,\1}:\underline{\cat{D}}(V,\1)\otimes V \to \1
    \end{equation}
    and
    \begin{equation}
        \coev_V : \1 \xrightarrow{u_V} \underline{\cat{D}}(V) \xrightarrow{\underline{\cat{D}}(V,r_V^{-1})} \underline{\cat{D}}(V,V\otimes\1) \xrightarrow{(\Phi^V_{V,\1})^{-1}} V\otimes \underline{\cat{D}}(V,\1) \;.
    \end{equation}
    We need to check the zig-zag identities. For brevity we suppress the associators in what follows. We check the identity
    \begin{equation}\label{eq:zigzagA}
        \begin{tikzcd}
            V && V \\ \1\otimes V && V\otimes\1 \\ & (V\otimes \underline{\cat{D}}(V,\1))\otimes V
            \arrow["{l_V^{-1}}"', from=1-1, to=2-1]\arrow["{\coev_V\otimes V}"', from=2-1, to=3-2]\arrow["{V\otimes\ev_V}"', from=3-2, to=2-3]\arrow["{r_V}"', from=2-3, to=1-3]
            \arrow[equals, from=1-1, to=1-3]
        \end{tikzcd} \;.
    \end{equation}
    The other identity follows via a similar computation. The diagram \eqref{eq:zigzagA} is the outside pentagon of
    \begin{equation}\label{eq:zigzagB}
        \begin{tikzcd}
            V && V \\
            \1\otimes V && V\otimes \1 \\ \underline{\cat{D}}(V)\otimes V && (V\otimes\underline{\cat{D}}(V,\1))\otimes V \\ & \underline{\cat{D}}(V,V\otimes\1)\otimes V &
            \arrow["{l_V^{-1}}"', from=1-1, to=2-1]\arrow["{u_V\otimes V}"', from=2-1, to=3-1]\arrow["{\underline{\cat{D}}(V,r_V^{-1})\otimes V}"', from=3-1, to=4-2]\arrow["{(\Phi_{V,\1}^V)^{-1} \otimes V}"', from=4-2, to=3-3]\arrow["{V\otimes\ev_{V,\1}}"', from=3-3, to=2-3]\arrow["{r_V}"', from=2-3, to=1-3]
            \arrow[equals, from=1-1, to=1-3]
            \arrow["{\ev_{V,V}}", from=3-1, to=1-3]
            \arrow["{\ev_{V,V\otimes\1}}",bend left=30, from=4-2, to=2-3]
        \end{tikzcd} \;.
    \end{equation}
    The square in the top left of \eqref{eq:zigzagB} commutes from the second diagram in Lemma~\ref{lemma:ClosedStructureA} and the definition of $u_V$, and the square in the middle from the naturality of $\ev_{V,\blank}$. It remains to show that the bottom right triangle in \eqref{eq:zigzagB} commutes. This triangle is the outside of the diagram
    \begin{equation}\label{eq:zigzagC}
        \begin{tikzcd}
            \underline{\cat{D}}(V,V\otimes\1)\otimes V && V\otimes\1 \\
            & \underline{\cat{D}}(\1,V\otimes\1)\otimes\1 & \\
            & \underline{\cat{D}}(\1,V\otimes\1)\otimes\underline{\cat{D}}(V,\1)\otimes V & \\ && \\
            & V\otimes\underline{\cat{D}}(V,\1)&
            \arrow["{\ev_{V,V\otimes\1}}\otimes V", from=1-1, to=1-3]
            \arrow["{(\Phi_{V,\1}^V)^{-1}}\otimes V"', bend right=30, from=1-1, to=5-2]\arrow["{V\otimes\ev_{V,\1}}"', bend right=30, from=5-2, to=1-3]
            \arrow["{\coev_{V,\1}\otimes\underline{\cat{D}}(V,\1)\otimes V}" description, from=5-2, to=3-2]\arrow["{m_{\1,V\otimes\1}^V\otimes V}", from=3-2, to=1-1]
            \arrow["{\underline{\cat{D}}(\1,V\otimes\1)\otimes\ev_{V,\1}}"', from=3-2, to=2-2]\arrow["{\ev_{\1,V\otimes\1}}"', from=2-2, to=1-3]
        \end{tikzcd}
    \end{equation}
    The triangle on the bottom left of \eqref{eq:zigzagC} commutes by the definition of $\Phi_{V,\1}^V$, and the upper left square of \eqref{eq:zigzagC} commutes from the definition of $m_{\1,V\otimes\1}^V$ and the second diagram in Lemma~\ref{lemma:ClosedStructureA}. The remaining square in \eqref{eq:zigzagC} can be reshaped into
    \begin{equation}\label{eq:zigzagD}
        \begin{tikzcd}
            V\otimes\underline{\cat{D}}(V,\1)\otimes V & V\otimes\1 & \underline{\cat{D}}(\1,V\otimes \1) \otimes\1
            \arrow["{V\otimes\ev_{V,\1}}", from=1-1, to=1-2]
            \arrow["{\coev_{V,\1}\otimes \1}", bend left=30, from=1-2, to=1-3]\arrow["{\ev_{\1,V\otimes\1}}", bend left=30, from=1-3, to=1-2]
        \end{tikzcd} \;.
    \end{equation}
    Upon setting $U = V$, $V = \1$, $W = V\otimes\1$ and $f = \Id_{V\otimes\1}$ in the second diagram in Lemma~\ref{lemma:ClosedStructureA} we obtain the loop in \eqref{eq:zigzagD}, so \eqref{eq:zigzagD} commutes. Thus \eqref{eq:zigzagC} commutes.
\end{proof}

Let $\cat{C}\subset\cat{D}$ be a full rigid monoidal subcategory and $\cat{M}\subset\cat{D}$ a full subcategory which is closed under tensoring with objects in $\cat{C}$, making $\cat{M}$ a $\cat{C}$-module category. Given $M\in\cat{M}$, we consider the algebra $\underline{\cat{D}}(M)$
and the $\cat{D}$-module category $\cat{D}_{\underline{\cat{D}}(M)}$. By virtue of the fact that $\cat{C}$ is a monoidal subcategory of $\cat{D}$, the category $\cat{D}_{\underline{\cat{D}}(M)}$ is also a $\cat{C}$-module category. Consider the functor $\underline{\cat{D}}(M,\blank):\cat{M} \to \cat{D}_{\underline{\cat{D}}(M)}$. The following is a corollary of Lemma~\ref{lem:PhisoIffRigid}.

\begin{corollary}\label{cor:C-module-structure-from-rigid}
    Let $\cat{C}\subset\cat{D}$ be a full rigid monoidal subcategory and $\cat{M}\subset\cat{D}$ a full subcategory which is closed under tensoring with objects in $\cat{C}$. Given $M\in\cat{M}$, the natural isomorphism $\Phi^M_{\blank,\bblank}$ gives the functor
    \begin{equation}
        \underline{\cat{D}}(M,\blank):\cat{M} \to \cat{D}_{\underline{\cat{D}}(M)} \;,\quad N \mapsto \underline{\cat{D}}(M,N)
    \end{equation}
     the structure of a $\cat{C}$-module functor.
\end{corollary}

An object $N\in\cat{M}$ is said to be a retract of an object $M\in\cat{M}$ if there exist morphisms $s:N\to M$ and $p:M\to N$ such that $p\circ s = \Id_N$. The following proposition is in the spirit of the characterisation given in \cite[Thm.\,7.10.1]{EGNO15} in the non-finite and non-rigid setting.
\begin{proposition}\label{prop:FFMCEmbedding}
    Let $\cat{C}\subset\cat{D}$ be a full rigid monoidal subcategory and $\cat{M}\subset\cat{D}$ a full subcategory which is closed under tensoring with objects in $\cat{C}$. Assume there exists an $M\in\cat{M}$ such that every $N\in\cat{M}$ is a retract of $U\otimes M$ for some $U\in\cat{C}$. The $\cat{C}$-module functor $\underline{\cat{D}}(M,\blank):\cat{M}\to\cat{D}_{\underline{\cat{D}}(M)}$ is fully faithful.
\end{proposition}
\begin{proof}
    Say $N\in\cat{M}$ is of the form $V\otimes M$ for some $V\in\cat{C}$. Given $L\in\cat{M}$, we have
    \begin{align}\label{eq:FFonInduced}
        \cat{D}_{\underline{\cat{D}}(M)}(\underline{\cat{D}}(M,V\otimes M),\underline{\cat{D}}(M,L)) &\cong \cat{D}_{\underline{\cat{D}}(M)}(V\otimes\underline{\cat{D}}(M),\underline{\cat{D}}(M,L)) \nonumber\\
        &\cong \cat{D}(V,\underline{\cat{D}}(M,L)) \nonumber\\
        &\cong \cat{M}(V\otimes M,L) \nonumber\\
        &= \cat{M}(N,L) \;,
    \end{align}
    where the second line follows from the fact that $\cat{D}_A(X\otimes A,B) \cong \cat{D}(X,B)$ for all $X,B\in\cat{D}$, and any algebra $A\in\cat{D}$. As the $\cat{C}$-module structure $\Phi^M_{\blank,\bblank}$ is defined via the adjunction data for $\blank\otimes M$ and $\underline{\cat{D}}(M,\blank)$, the isomorphism \eqref{eq:FFonInduced} is the map induced by $\underline{\cat{D}}(M,\blank)$.
    
    Define $F:=\underline{\cat{D}}(M,\blank)$ and denote by $\phi_{N,L}$ the canonical maps $\cat{M}(N,L)\to \cat{D}_{\underline{\cat{D}}(M)}(F(N),F(L))$, which per \eqref{eq:FFonInduced} are isomorphisms when $N\cong V\otimes M$ for some $V\in\cat{C}$. Choose $U\in\cat{C}$ such that $N\in\cat{M}$ is a retract of $U\otimes M$. Consider the diagram
    \begin{equation}\label{eq:Retracts preserve isos}
        \begin{tikzcd}
            {\scriptstyle \cat{M}(N,L)} & {\scriptstyle \cat{M}(U\otimes M,L)} & {\scriptstyle \cat{M}(N,L)} \\
            {\scriptstyle \cat{D}_{\underline{\cat{D}}(M)}(F(N),F(L))} & {\scriptstyle \cat{D}_{\underline{\cat{D}}(M)}(F(U\otimes M),F(L))} & {\scriptstyle \cat{D}_{\underline{\cat{D}}(M)}(F(N),F(L))}
            \arrow["{p^\ast}", from=1-1, to=1-2]\arrow["{s^\ast}", from=1-2, to=1-3]
            \arrow["{\phi_{N,L}}", from=1-1, to=2-1]\arrow["{\phi_{U\otimes M,L}}", from=1-2, to=2-2]\arrow["{\phi_{N,L}}", from=1-3, to=2-3]
            \arrow["{F(p)^\ast}", from=2-1, to=2-2]\arrow["{F(s)^\ast}", from=2-2, to=2-3]
        \end{tikzcd} \;.
    \end{equation}
    Functoriality of $F$ implies that the diagram commutes. Since $p\circ s = \Id_N$, composition along both rows is the identity. Define
    \begin{equation}
        \rho_{N,L} := s^\ast \circ \phi_{U\otimes M,L}^{-1} \circ F(p)^\ast: \cat{D}_{\underline{\cat{D}}(M)}(F(N),F(L))\to\cat{M}(N,L) \;.
    \end{equation}
    Commutativity of \eqref{eq:Retracts preserve isos} implies $\phi_{N,L}\circ\rho_{N,L} = \Id$ and $\rho_{N,L}\circ\phi_{N,L} = \Id$. Thus $\underline{\cat{D}}(M,\blank)$ is fully faithful.
\end{proof}


\begin{thebibliography}{AAA99}
\small
\setlength{\itemsep}{-1pt}

\bibitem[AGS94a]{AGS94a}
A.~Alekseev, H.~Grosse, V.~Schomerus,
\textit{Combinatorial quantization of the Hamiltonian Chern--Simons theory~I},
\doi{10.1007/BF02099431}{Communications in Mathematical Physics, \textbf{172} (1995): 317–358},
\arxiv{hep-th/9403066}{[hep-th]}.
 
\bibitem[AGS94b]{AGS94b}
A.~Alekseev, H.~Grosse, V.~Schomerus,
\textit{Combinatorial quantization of the Hamiltonian Chern--Simons theory~II},
\doi{10.1007/BF02101528}{Communications in Mathematical Physics, \textbf{174} (1996): 561–604},
\arxiv{hep-th/9408097}{[hep-th]}.

\bibitem[Al93]{Al93}
A.~Alekseev,
\textit{Integrability in the Hamiltonian Chern--Simons theory},
St.~Petersburg Mathematical Journal, \textbf{6} (1995): 241–253,
\arxiv{hep-th/9311074}{[hep-th]}.
 
\bibitem[AM11]{AM11}
H. H.~Andersen, V.~Mazorchuk,
\textit{Category {$\mathcal{O}$} for quantum groups},
\doi{10.4171/JEMS/506}{Journal of the European Mathematical Society, \textbf{17} (2015): 405–431},
\arxiv{1105.5500}{[math.RT]}.
 
\bibitem[AMM97]{AMM97}
A.~Alekseev, A.~Malkin, E.~Meinrenken,
\textit{Lie group valued moment maps},
\doi{10.4310/jdg/1214460860}{Journal of Differential Geometry, \textbf{48} (1998): 445–495},
\arxiv{dg-ga/9707021}{[math.DG]}.
 
\bibitem[ARS99]{ARS99}
A.~Alekseev, A.~Recknagel, V.~Schomerus,
\textit{Non-commutative world-volume geometries: branes on SU(2) and fuzzy spheres},
\doi{10.1088/1126-6708/1999/09/023}{Journal of High Energy Physics, \textbf{09} (1999): 023},
\arxiv{hep-th/9908040}{[hep-th]}.
 
\bibitem[AS95]{AS95}
A.~Alekseev, V.~Schomerus,
\textit{Representation theory of Chern--Simons observables},
\doi{10.1215/S0012-7094-96-08519-1}{Duke Mathematical Journal, \textbf{85} (1996): 447–510},
\arxiv{q-alg/9503016}{[math.QA]}.
 
\bibitem[BBJ15]{BBJ15}
D.~Ben-Zvi, A.~Brochier, D.~Jordan,
\textit{Integrating quantum
groups over surfaces},
\doi{10.1112/topo.12072}{Journal of Topology, \textbf{11} (2018):874–917},
\arxiv{1501.04652}{[math.QA]}.
 
\bibitem[BBJ16]{BBJ16}
D.~Ben-Zvi, A.~Brochier, D.~Jordan,
\textit{Quantum character varieties and braided module categories},
\doi{10.1007/s00029-018-0426-y}{Selecta Mathematica, \textbf{24} (2018):4711–4748},
\arxiv{1606.04769}{[math.QA]}.

\bibitem[BGJV26]{BGJV26}
J.~Brown, J. R. G\'omez García, D.~Jordan, M.~Vancraeynest,
\textit{Weyl defects in skein theory and quantum cluster charts},
\arxiv{2609.09308}{[math.QA]}.
 
\bibitem[BJ25]{BJ25}
J.~Brown, D.~Jordan,
\textit{Parabolic skein modules},
\arxiv{2505.14836}{[math.QA]}.
 
\bibitem[BJS18]{BJS18}
A.~Brochier, D.~Jordan, N.~Snyder,
\textit{On dualizability of braided tensor categories},
\doi{10.1112/S0010437X20007630}{Compositio Mathematica, \textbf{157} (2021):435-483},
\arxiv{1804.07538}{[math.QA]}.

\bibitem[BR94]{BR94}
E.~Buffenoir, Ph.~Roche,
\textit{Two dimensional lattice gauge theory based on a quantum group},
\doi{10.1007/BF02099153}{Communications in Mathematical Physics, \textbf{170} (1995): 669–698},
\arxiv{hep-th/9405126}{[hep-th]}.
 
\bibitem[Br12]{Br12}
A.~Brochier,
\textit{Cyclotomic associators and finite type invariants for tangles in the solid torus},
\doi{10.2140/agt.2013.13.3365}{Algebraic \& Geometric Topology, \textbf{13} (2013): 3365–3409},
\arxiv{1209.0417}{[math.QA]}.

\bibitem[CGS22]{CGS22}
D.~Chernyak, A. M.~Gainutdinov, H.~Saleur,
\textit{{$U_q\mathfrak{sl}_2$}-invariant non-compact boundary conditions for the XXZ spin chain},
\doi{10.1007/JHEP11(2022)016}{Journal of High Energy Physics, \textbf{11} (2022): 16},
\arxiv{2207.12772}{[hep-th]}.
 
\bibitem[CGP14]{CGP14}
F.~Costantino, N.~Geer, B.~Patureau-Mirand,
\textit{Some remarks on the unrolled quantum group of $\mathfrak{sl}(2)$},
\doi{10.1016/j.jpaa.2014.10.012}{Journal of Pure and Applied Algebra, \textbf{219} (2015):3238-3262},
\arxiv{1406.0410}{[math.QA]}.

\bibitem[Cha94]{Cha94}
V.~Chari,
\textit{Minimal affinizations of representations of quantum groups: the $U_q(\mathfrak{g})$-module structure},
\arxiv{hep-th/9411144}{[hep-th]}.

\bibitem[Che84]{Ch84}
I.~V.~Cherednik,
\textit{Factorizing particles on a half-line and root systems},
\doi{10.1007/BF01038545}{Theoretical and Mathematical Physics, \textbf{61} (1984): 977–983}.

\bibitem[CJR26]{CJR26}
E. Y.-Y.~Chen, D.~Jordan, I.~Romaidis,
\textit{Skein theory, line defects, and quantum symmetric pairs},
\arxiv{2609.18902}{[math.QA]}.

\bibitem[CP94]{CP94}
V.~Chari, A.~Pressley,
\textit{A Guide to Quantum Groups},
Cambridge University Press (1994), Cambridge.
 
\bibitem[DK90]{DK90}
C.~De~Concini, V.~G.~Kac,
\textit{Representations of Quantum Groups at Roots of 1},
\href{https://math.mit.edu/~kac/not-easily-available/representations.pdf}{Operator Algebras, Unitary Representations, Enveloping Algebras, and Invariant Theory,
Progress in Mathematics, \textbf{92} (1990): 471-506, Birkh\"auser}.

\bibitem[DN20]{DN20}
A.~Davydov, D.~Nikshych,
\textit{Braided Picard groups and graded extensions of braided tensor categories},
\doi{10.1007/s00029-021-00670-1}{Selecta Mathematica, \textbf{27} (2021): 65},
\arxiv{2006.08022}{[math.QA]}.
 
\bibitem[EGNO15]{EGNO15}
P.~Etingof, S.~Gelaki, D.~Nikshych, V.~Ostrik,
\textit{Tensor Categories},
\doi{10.1090/surv/205}{Mathematical Surveys and Monographs, \textbf{205} (2015), American Mathematical Society, Providence, RI}.

\bibitem[EK66]{EK66}
S.~Eilenberg, G. M.~Kelly,
\textit{Closed Categories},
\doi{10.1007/978-3-642-99902-4_22}{Proceedings of the Conference on Categorical Algebra -- La Jolla 1965, (1966): 421-562}.

\bibitem[En04]{En04}
B.~Enriquez,
\textit{Quasi-reflection algebras and cyclotomic associators},
\doi{10.1007/s00029-007-0048-2}{Selecta Mathematica, \textbf{13} (2008): 391-463},
\arxiv{math/0408035}{[math.QA]}.

\bibitem[ES99]{ES99}
P.~Etingof, O.~Schiffmann,
\textit{Lectures on the dynamical Yang--Baxter equations},
\doi{10.1017/CBO9780511542848.007}{Quantum groups and Lie theory (Durham 1999), LMS Lecture Note Series \textbf{290}
(2001): 89-129},
\arxiv{math/9908064}{[math.QA]}.

\bibitem[FHW14]{FHW14}
V.~Futorny, J.~Hartwig, E.~Wilson,
\textit{Irreducible completely pointed modules of quantum groups of type A},
\doi{10.1016/j.jalgebra.2015.03.006}{Journal of Algebra, \textbf{432} (2015): 252-279},
\arxiv{1404.0305}{[math.RT]}.

\bibitem[FrRS]{paper2-skeins}
T.~Franke, I.~Runkel, V.~Schomerus, 
\textit{Skein modules with line and point defects}, in preparation.
 
\bibitem[FuRS02]{FRS02}
J.~Fuchs, I.~Runkel, C.~Schweigert,
\textit{TFT construction of RCFT correlators I: partition functions},
\doi{10.1016/S0550-3213(02)00744-7}{Nuclear Physics B, \textbf{646} (2002): 353–497},
\arxiv{hep-th/0204148}{[hep-th]}.
 
\bibitem[GJS19]{GJS19}
S.~Gunningham, D.~Jordan, P.~Safronov,
\textit{The finiteness conjecture for skein modules},
\doi{10.1007/s00222-022-01167-0}{Inventiones Mathematicae, \textbf{232} (2023): 301-363},
\arxiv{1908.05233}{[math.QA]}.

\bibitem[Go26]{G26}
J. R. G\'omez García,
\textit{HOMFLY parabolic restriction, defect skein theory and the Turaev coproduct},
\arxiv{2601.03196}{[math.QA]}.

\bibitem[GP10]{GP10}
N.~Geer, B.~Patureau-Mirand,
\textit{Topological invariants from non-restricted quantum groups},
\doi{10.2140/agt.2013.13.3305}{Algebraic \& Geometric Topology, \textbf{13} (2013): 3305–3363},
\arxiv{1009.4120}{[math.GT]}.

\bibitem[GPT07]{GPT07}
N.~Geer, B.~Patureau-Mirand, V.~Turaev,
\textit{Modified quantum dimensions and re-normalized link invariants},
\doi{10.1112/S0010437X08003795}{Compositio Mathematica, \textbf{145} (2009): 196-212},
\arxiv{0711.4229}{[math.QA]}.
 
\bibitem[GW06]{GW06}
S.~Gukov, E.~Witten,
\textit{Gauge theory, ramification, and the geometric Langlands program},
\doi{10.4310/CDM.2006.v2006.n1.a2}{Current Developments in Mathematics 2006, International Press (2008): 35–180},
\arxiv{hep-th/0612073}{[hep-th]}.
 
\bibitem[HS26]{HS26}
A.~Hofstetter, C.~Schweigert,
\textit{The 2-categorical S-matrix of a braided fusion 1-category is a character table},
\arxiv{2601.14168}{[math.QA]}.
 
\bibitem[ILZ18]{ILZ18}
K.~Iohara, G.~Lehrer, R.~Zhang,
\textit{Schur--Weyl duality for certain infinite dimensional $U_q(\mathfrak{sl}_2)$-modules},
\arxiv{1811.01325}{[math.RT]}.
 
\bibitem[Jo23]{Jo23}
D.~Jordan,
\textit{Langlands duality for skein modules of 3-manifolds},
String-Math 2022, Proceedings of Symposia in Pure Mathematics, AMS, \textbf{107} (2024): 127-150,
\arxiv{2302.14734}{[math.QA]}.

\bibitem[JR21]{JR21}
T.~Johnson-Freyd, D.~Reutter,
\textit{Minimal nondegenerate extensions},
\doi{10.1090/jams/1023}{Journal of the American Mathematical Society, \textbf{37} (2024): 81-150},
\arxiv{2105.15167}{[math.QA]}.
 
\bibitem[JS93]{JS93}
A.~Joyal, R.~Street,
\textit{Braided tensor categories}, \doi{10.1006/aima.1993.1055}{Advances in Mathematics, \textbf{102} (1993): 20-78}.
 
\bibitem[Ka95]{Ka95}
C.~Kassel,
\textit{Quantum Groups},
\doi{10.1007/978-1-4612-0783-2}{Graduate Texts in Mathematics \textbf{155}, Springer-Verlag, New York, 1995}.
 
\bibitem[Ke23]{Ke23}
C.~Keller,
\textit{Generalized Character Varieties and Quantization via Factorization Homology},
\href{https://imag.umontpellier.fr/~calaque/students-Corina_Keller_PhD_Thesis.pdf}{Ph.D.~thesis, Université de Montpellier (2023)}.

\bibitem[KL01]{KL01}
T.~Kerler, V. V.~Lyubashenko,
\textit{Non-Semisimple Topological Quantum Field Theories for 3-Manifolds with Corners},
\doi{10.1007/b82618}{Lecture Notes in Mathematics \textbf{1765} (2001), Springer-Verlag, Berlin}.
 
\bibitem[Ko17]{Ko17}
S.~Kolb,
\textit{Braided module categories via quantum symmetric pairs},
\doi{10.1112/plms.12303}{Proceedings of the London Mathematical Society, \textbf{121} (2020):1-31},
\arxiv{1705.04238}{[math.QA]}.
 
\bibitem[KR26]{KR26}
P.~Kinnear, I.~Runkel,
\textit{Defects in skein theory and TQFT},
\arxiv{2606.07432}{[math.QA]}.
 
\bibitem[KSc97]{KS97}
A. U.~Klimyk, K.~Schm\"udgen,
\textit{Quantum Groups and Their Representations},
\doi{10.1007/978-3-642-60896-4}{Texts and Monographs in Physics, Springer, Berlin (1997)}.

\bibitem[KSk92]{KS92}
P.~P.~Kulish, E.~K.~Sklyanin,
\textit{Algebraic structures related to reflection equations},
\doi{10.1088/0305-4470/25/22/022}{Journal of Physics A: Mathematical and General, \textbf{25} (1992): 5963–5975},
\arxiv{hep-th/9209054}{[hep-th]}.
 
\bibitem[LNV20]{LNV20}
A.~Lacabanne, G.~Naisse, P.~Vaz,
\textit{Tensor product categorifications, Verma modules and the blob 2-category},
\doi{10.4171/qt/156}{Quantum Topology, \textbf{12} (2021):705-812},
\arxiv{2005.06257}{[math.QA]}.
 
\bibitem[LWY23]{LWY23}
R.~Laugwitz, C.~Walton, M.~Yakimov,
\textit{Reflective centers of module categories and quantum K-matrices},
\doi{10.1017/fms.2025.10055}{Forum of Mathematics, Sigma, \textbf{13} (2025):e95},
\arxiv{2307.14764}{[math.QA]}.

\bibitem[MN92]{MN92}
L.~Mezincescu, R.~I.~Nepomechie,
\textit{Fusion procedure for open chains},
\doi{10.1088/0305-4470/25/9/024}{Journal of Physics A: Mathematical and General, \textbf{25} (1992): 2533–2543}.
 
\bibitem[MSa93]{MS93}
P.~Martin, H.~Saleur,
\textit{The blob algebra and the periodic Temperley--Lieb algebra},
\doi{10.1007/BF00805852}{Letters in Mathematical Physics, \textbf{30} (1994):189-206},
\arxiv{hep-th/9302094}{[hep-th]}.
 
\bibitem[MSe89]{MS89}
G.~Moore, N.~Seiberg,
\textit{Classical and quantum conformal field theory},
\doi{10.1007/BF01238857}{Communications in Mathematical Physics, \textbf{123} (1989): 177–254}.
 
\bibitem[Mu04]{Mu04}
A.~Mudrov,
\textit{Quantum Conjugacy Classes of Simple Matrix Groups},
\doi{10.1007/s00220-007-0222-6}{Communications in Mathematical Physics, \textbf{272} (2007):635-660},
\arxiv{math/0412538}{[math.QA]}.

\bibitem[PS90]{PS90}
V.~Pasquier, H.~Saleur,
\textit{Common structures between finite systems and conformal field theories through quantum groups},
\doi{10.1016/0550-3213(90)90122-T}{Nuclear Physics B, \textbf{330} (1990): 523–556}.
 
\bibitem[PW20]{PW20}
L.~Pol, J.~Williamson,
\textit{The Homotopy Theory of Complete Modules},
\doi{10.1016/j.jalgebra.2021.11.030}{Journal of Algebra, \textbf{594} (2022): 74-100},
\arxiv{2011.06989}{[math.AC]}.
 
\bibitem[Ru19]{Ru19}
M.~Rupert,
\textit{Categories of Weight Modules for Unrolled Restricted Quantum Groups at Roots of Unity},
\doi{10.4213/im9234e}{Izvestiya: Mathematics, \textbf{86} (2022):1204-1224},
\arxiv{1910.05922}{[math.RT]}.
 
\bibitem[Sk88]{Sk88}
E.~K.~Sklyanin,
\textit{Boundary conditions for integrable quantum systems},
\doi{10.1088/0305-4470/21/10/015}{Journal of Physics A: Mathematical and General, \textbf{21} (1988): 2375–2389}.
 
\bibitem[SZ24]{SZ24}
M.~Stroiński, T.~Zorman,
\textit{Reconstruction of module categories in the infinite and non-rigid settings},
\arxiv{2409.00793}{[math.QA]}.
 
\bibitem[Te95]{Te95}
J.~Teschner,
\textit{On the Liouville three-point function},
\doi{10.1016/0370-2693(95)01200-A}{Physics Letters B, \textbf{363} (1995): 65–70},
\arxiv{hep-th/9507109}{[hep-th]}.
 
\bibitem[Ti10]{Ti10}
P.~Tingley,
\textit{A minus sign that used to annoy me but now I know why it is there},
\arxiv{1002.0555}{[math.GT]}.

\bibitem[Wi89]{Wi89}
E.~Witten,
\textit{Quantum field theory and the Jones polynomial},
\doi{10.1007/BF01217730}{Communications in Mathematical Physics, \textbf{121} (1989): 351–399}.
 
\bibitem[Wi11]{Wi11}
E.~Witten,
\textit{Fivebranes and knots},
\doi{10.4171/QT/26}{Quantum Topology, \textbf{3} (2012): 1–137},
\arxiv{1101.3216}{[hep-th]}.
 
\bibitem[WW26]{WW26}
L.~Wang, Z.~Wojciechowski,
\textit{No braiding for quantum $\mathfrak{sl}_2$ reps -- A short proof},
\arxiv{2610.29748}{[math.QA]}.
 
\end{thebibliography}
\end{document}